%% file: main.tex
\documentclass{imsart}
\RequirePackage{amsthm,amsmath,amsfonts,amssymb}
\RequirePackage[numbers]{natbib}
\RequirePackage[colorlinks,citecolor=blue,urlcolor=blue]{hyperref}
\RequirePackage{graphicx}

\startlocaldefs
\theoremstyle{plain}

\newtheorem{theorem}{Theorem}[]
\newtheorem{lemma}[]{Lemma}
\newtheorem{corollary}{Corollary}
\theoremstyle{definition}

\newtheorem{assumption}{Assumption}

\theoremstyle{remark}

\usepackage{appendix}
\input{preamble.tex}

\endlocaldefs

\begin{document}
\begin{frontmatter}
\title{Learning cross-layer dependence structure in multilayer networks}
\runtitle{Learning cross-layer dependence structure in multilayer networks}

\begin{aug}
\author[A]{\fnms{Jiaheng}~\snm{Li}\ead[label=e1]{jl20gx@fsu.edu}},
\author[A]{\fnms{Jonathan}~\snm{Stewart}\ead[label=e2]{jrstewart@fsu.edu}}

\address[A]{Department of Statistics,
Florida State University\printead[presep={,\ }]{e1,e2}}
\runauthor{J. Li et al.}
\end{aug}


\begin{abstract}
We propose a novel class of separable multilayer network models to capture cross-layer dependencies in multilayer networks,
enabling the analysis of how interactions in one or more layers may influence interactions in other layers.
Our approach separates the network formation process from the layer formation process,
and is able to extend existing single-layer network models to multilayer network models
that accommodate cross-layer dependence.
We establish non-asymptotic and minimax-optimal error bounds for maximum likelihood estimators
and demonstrate the convergence rate
in scenarios of increasing parameter dimension.
Additionally, we establish non-asymptotic error bounds for multivariate normal approximations and propose a model selection method that controls the false discovery rate. 
Simulation studies and an application to the Lazega lawyers network show that our framework and method perform well in realistic settings.
\end{abstract}


\end{frontmatter}

\section{Introduction} 

\input{introduction_ml.tex}

The rest of the paper is organized as follows. 
Section \ref{sec2} introduces our modeling framework and includes illustrative examples.
The consistency and the minimax optimal results are contained in Section \ref{sec3}.
The multivariate normal approximation theory is presented in Section \ref{sec:normal}. 
The results of simulation studies are provided in Section \ref{sec:sim}, 
together with different testing procedures for model selection
which control the false discovery rate.
An application of our developed framework and methodology is given in Section \ref{sec:app},
concluding with a discussion presented in Section \ref{sec:disc}.
The code and data to reproduce the simulations and analyses can be found in our package online.\footnote{
\url{https://github.com/jiaheng-li/mlyrnetwork}}

\section{Modeling cross-layer dependence in multilayer networks} 
\label{sec2}

A multilayer network can be represented as 
a sequence of $1 \leq K < \infty$ random graphs $\bX^{(1)}, \ldots, \bX^{(K)}$
each defined on a common set of $N \geq 3$ nodes, 
which we take without loss to be the set $\mN = \{1,\ldots,N\}$. 
We call the graphs $\bX^{(1)}, \ldots, \bX^{(K)}$ the {\it layers} of the network,
and represent the multilayer network as the quantity $\bX = (\bX^{(1)}, \ldots, \bX^{(K)})$. 

Connections between pairs of nodes $\{i,j\} \subset \mN$ in each layer $k \in \{1, \ldots, K\}$
are modeled by random variables  
\beno
X_{i,j}^{(k)}
\= \begin{cases}
1 & \mbox{nodes } i \mbox{ and } j \mbox{ are connected in layer } k \\
0 & \mbox{otherwise}
\end{cases}. 
\ee
We refer to all connections of a pair of nodes $\{i,j\} \subset \mN$ across the $K$ layers 
as a {\it  dyad} which we denote by  
$\bX_{i,j} = (X_{i,j}^{(1)}, \ldots, X_{i,j}^{(K)})\in \{0,1\}^{K}$. 
A multilayer network can be represented by a collection of dyads as $\bX = (\bX_{i,j})_{\{i,j\} \subset \mN}$ alternatively.

For notational ease, 
we will consider undirected multilayer networks, 
which imply that the network layers $\bX^{(1)}, \ldots, \bX^{(K)}$ are undirected random graphs; 
extensions to directed multilayer networks or mixed multilayer networks with both directed and undirected 
layers will typically be straightforward, 
involving only notational adaptations in subscripts in most cases. 
We adopt the usual conventions for undirected networks, 
i.e., 
we assume that $X_{i,j}^{(k)} = X_{j,i}^{(k)}$ (all $\{i,j\} \subset \mN$, $1 \leq k \leq K$) 
and $X_{i,i}^{(k)} = 0$ (all $i \in \mN$, $1 \leq k \leq K$). 
The sample space of each layer $\bX^{(k)}$ is therefore 
the product space $\mbX^{(k)} \coloneqq \{0, 1\}^{\binom{N}{2}}$ ($k = 1, \ldots, K$), 
and the sample space $\mbX$ of $\bX$ is the product space of the sample spaces of the individual layers, 
i.e.,
$\mbX \coloneqq \mbX^{(1)} \times \cdots \times \mbX^{(K)}$. 
The sample space of dyad $\{i,j\} \subset \mN$ is the product space $\mbX_{i,j} \coloneqq \{0,1\}^{K}$. 

A challenge in the statistical modeling of network data lies in the fact that 
networks have many distinguishing properties, 
including:  
\ben
\item {\bf Sparsity.} Many real-world networks are sparse, 
in the sense that the expected number of edges in the network grows at a rate slower than $\binom{N}{2}$. 
The phenomena of network sparsity manifests in a variety of different applications, 
usually due to constraints, 
such as time or financial constraints, 
which can limit the number of connections any node can maintain at a given point in time
\citep[][]{KrHaMo11,butts:jms:2018}.
\item {\bf Node heterogeneity.} Different actors in a social network will have different properties,
called node covariates,
which can lead to different propensities to form edges. 
A key example is assortative and disassortative mixing patterns in networks 
\citep{McSmCo01,KrHaRaHo07},
as well as differences in structural patterns in the network \citep{Albert02, LiXu12}.  
\item {\bf Edge dependence.} In addition to node-based effects that give rise to 
heterogeneity 
in propensities for nodes to form edges, 
scientific and statistical evidence suggests edges are dependent in many applications 
\citep{HpLs72,Fo80,block2015reciprocity},
and modeling a single system of multiple binary random variables  without replication is a challenging statistical problem
inherent to many statistical network analysis applications. 
\een
Each of the above gives rise to distinct challenges for modeling network data and performing statistical inference
in statistical network analysis applications,  
and it is not straightforward to construct
models that due justice to each of these and more.  
To address these challenges, 
a plethora of statistical models have been proposed to model network data, 
which for single-layer networks have included  
exponential-families of random graph models 
\citep[e.g.,][]{ergm.book,ScKrBu17},
stochastic block models
\citep[e.g.,][]{HoLaLe83},
latent metric space models 
\citep[e.g.,][]{HpRaHm01},
random dot product graphs \citep[e.g.,][]{Athreya2018}, 
exchangeable random graph models \citep[e.g.,][]{CaFo17,CrDe16}, 
and more. 
In this work, 
we build upon the many classes of single-layer network data models by introducing a separable multilayer network modeling framework. This framework enables existing single-layer network models to be extended to the multilayer setting and simultaneously enables learning cross-layer dependence and interactions across different layers in the multilayer network.

\subsection{Separable multilayer network models}  
\label{sec:2.1}
Multilayer networks are subject to the same forces and phenomena as single layer networks, as multiple modes of relation or interaction do not remove constraints or properties of nodes which are fundamental to network data applications. 
The same set of nodes is defined across all layers in a multilayer network, and because all layers share the same set of nodes, the dyadic connections among these nodes fundamentally define the network formation process. By specifying a single-layer network as the foundational structure reference, we can separate the network formation process from the layer formation process. In doing so, the single-layer network serves as the baseline for establishing dyadic relationships that represent the relational structure across all layers of the multilayer network. 
As a result, the network formation process determines which dyads have the potential to form connections, i.e., which pair of nodes may exhibit at least one edge in any of the layers. 
In contrast,
The layer formation process dictates the particular layers in which these  connections appear. 
To learn the effects of cross-layer dependence in multilayer networks, we propose the class of separable multilayer network models, 
which extend the broad literature on single-layer network models into the multilayer realm.
These models can incorporate an arbitrary single-layer network structure as the foundational baseline and ensures that the underlying single-layer network can be recovered from observations of the multilayer network. We illustrate this approach and its advantages through our proposed modeling framework.
\input{general_model}

\input{prop_inference}

Proposition \ref{prop:inference} establishes a few key facts for the inference of cross-layer dependence structures 
in multilayer networks. 
First, 
we are able to observe $\bY$ through $\bX$,
as given any observation $\bx \in \mbX$ of the multilayer network $\bX$,
$\sepmodel(\bY = \by \,|\, \bX = \bx) = 1$ for one, and only one, $\by \in \mbY$. 
In other words, 
through the observation of $\bx$,
we can infer with probability $1$ the corresponding $\by$
due to the form of \eqref{general_model}.  
The significance of this result is that we do not need to treat the basis network $\bY$ as a latent network, 
which would require additional statistical and computational methodology to handle the latent missing network data.  
Second, 
we see that the inference for $\truth$ is unaffected by the choice of $g(\by)$; 
although, the statistical guarantees for estimators of $\truth$ will be indirectly influenced by the choice of $g(\by)$,
a point which we discuss in later sections.
Moreover, 
the above choice for $f(\bx, \nat)$
and the functional form of $\sepmodel(\bX = \bx \,|\, \bY = \by)$ derived in Proposition \ref{prop:inference}
establishes that $\log \, \sepmodel(\bX = \bx \,|\, \bY = \by)$
corresponds to the log-likelihood of a minimal exponential family, 
accessing a broad literature of statistical methodology and theory \citep[e.g.,][]{Su19}.
We note that other specifications for $f(\bx, \nat)$ are possible, 
but that Markov random field specifications provide a powerful class of models for dependent data 
\citep[e.g.,][]{WaJo08},
and in the case of the saturated model with maximal interaction term $H = K$, it
completely specifies all possible probabilities of outcomes $\bx_{i,j} \in \{0, 1\}^K$,
presenting a non-parametric model class for multilayer networks.

\subsection{Example of a multilayer network with pairwise interactions} 
\input{example_ml2}

\section{Estimation of cross-layer dependence structure} 
\label{sec3}
\input{estimation_setup}

\input{min_eig_lemma}

\s 

In classical settings with independent and identically distributed observations, 
the expected negative Hessian of the log-likelihood function 
is the Fisher information matrix and 
is expected to scale with the number of observations.  
In such cases, 
standard matrix theory indicates that the smallest eigenvalue of this expected negative Hessian matrix 
will scale with the sample size,
provided the smallest eigenvalue of the Fisher information matrix is bounded from below.
Lemma \ref{lem:min-eig} extends this notion by establishing similar scaling behavior concerning the expected number of activated dyads $\mbE \, \norm{\bY}_1$,
proxying as an effective sample size. 
Analogously,  $\mcI(\nat)$ can be seen as the Fisher information of the population distribution governing individual activated dyads in $\bY$, mirroring the role of Fisher information for population distributions in classical independent and identically distributed scenarios.  

Before we present our theoretical guarantees for maximum likelihood estimators in Theorem \ref{thm1}, we define some notations and outline some regularity assumptions for our theorem to follow.  
As we will show in Theorem \ref{thm1}, the choice of $g(\by)$ influences the estimation error through the expected 
number of edges in $\bY$ and through the covariances of edge variables in $\bY$. 
Define
\beno
D_{g}
&\coloneqq& \dsum_{\{i,j\} \prec \{v,w\} \subset \mN} \, \cov(Y_{i,j}, \, Y_{v,w}),
\ee
where $\{i,j\} \prec \{v,w\}$ implies the sum is taken with respect to the lexicographical ordering of pairs of nodes. 
Define $[D_{g}]^{+} \coloneqq \max\{0, \, D_{g}\}$ to be the positive part of $D_{g}$.
Let $\epsilon > 0$ be fixed independent of $N$ and $p$, and denote the $\epsilon$-ball of the data-generating parameter $\truth$ by $\mB_2(\truth, \epsilon) = \{\nat \in \mbR^p : \norm{\truth - \nat}_2 \leq \epsilon\}$. Define 
\beno
\widetilde{\lambda}_{\min}^{\epsilon}
\;\coloneqq\; \inf\limits_{\nat \in \mB_2(\truth, \epsilon)} \, \lambda_{\min}(\mcI(\nat)) 
&&\mbox{and}&&
\widetilde{\lambda}_{\max}^{\star}
\;\coloneqq\; \lambda_{\max}(\mcI(\truth)),  
\ee
where $\lambda_{\min}(\bA)$ and $\lambda_{\max}(\bA)$ are the smallest and the largest eigenvalue of matrix $\bA \in \mbR^{p \times p}$, respectively. 
\begin{assumption}
\label{assump1}
 Assume there exists a $C_0 > 0$ such that $\mbE\,\norm{\bY}_1 \ge 1$ and 
\beno
\dfrac{[D_{g}]^{+}}{\mbE\, \norm{\bY}_1} &\leq& C_0,
\ee
for all network sizes $N$. 
\end{assumption}

\begin{assumption}
\label{assump2}
Assume the parameter dimension $p$ satisfies 
\beno
p 
&\leq& \sqrt{\widetilde{\lambda}_{\max}^{\star} \, \mbE\, \norm{\bY}_1},
\ee
for all network sizes $N$. 
\end{assumption}

\begin{assumption}
\label{assump3}
Assume that 
$\widetilde{\lambda}_{\max}^{\star}$ and $\widetilde{\lambda}_{\min}^{\epsilon}$ satisfy,
as a function of the network size $N$,  
 \beno
    \dfrac{\sqrt{\widetilde{\lambda}_{\max}^{\star}} }{\widetilde{\lambda}_{\min}^{\epsilon}} \= o \;\left(\sqrt{\dfrac{\mbE \norm{\bY}_1}{p}} \right).
 \ee
\end{assumption}

Assumptions \ref{assump1}–\ref{assump3} provide a foundation for Theorem \ref{thm1} to establish the consistency result of the maximum likelihood estimator in large network settings. Assumption \ref{assump1} imposes a lower bound on the expected number of activated dyads relative to the covariance as the network size \(N\) grows. Assumption \ref{assump2} restricts the growth rate of \(p\) in relation to the network size and the largest eigenvalue of the Fisher information $\mcI(\truth)$. Finally, Assumption \ref{assump3} sets a constraint on the ratio between $\sqrt{\widetilde{\lambda}_{\max}^{\star}}$ and $\widetilde{\lambda}_{\min}^{\epsilon}$, balancing eigenvalue magnitudes in a way that preserves estimator consistency under increasing network size.

\input{thm1}

\section{Error of the normal approximation and model selection} 
\label{sec:normal}
In this section, 
we establish the asymptotic multivariate normality of the maximum likelihood estimator (MLE) for the data-generating parameter vector $\truth$ as its dimension grows.
Specifically,
we derive a non-asymptotic bound on the quality of the multivariate normal approximation and exhibit scaling conditions on both the model dimension $p$ 
and the expected number of activated dyads $\mbE \, \norm{\bY}_1$---under which the approximation error vanishes as the network size tends to infinity.  
Based on this result, 
we present a model selection method using multiple hypothesis testing procedures that control the false discovery rate. 
The main result is presented in Theorem \ref{thm2},
the proof of which 
is based on a Taylor expansion of the log-likelihood function  
and through the application of a Lyapunov type bound presented in \citet{Raic19}.

\input{thm2}


\subsection{Model selection via univariate testing with FDR control}

\input{global_test}

\section{Simulation studies}
\label{sec:sim}
Directly simulating maximum likelihood estimators for network data with dependent edges is challenging because the normalizing constants are often computationally intractable. Computing the normalizing constant requires enumerating all $2^{\binom{N}{2}}$ possible edge combinations for each layer to maximize the true likelihood function. Additionally, dependencies among network dyads prevent factorization of the likelihood, which further complicates direct maximization. As a result, direct maximization of likelihood functions is generally infeasible in these cases.
Two predominant methods of approximating the maximum likelihood estimator $\truth$ when the likelihood function is computationally intractable 
have emerged in the literature. 
Monte Carlo maximum likelihood estimation (MCMLE) \citep{GeTh92},
which constructs a simulation-based approximation to the likelihood function 
in order to approximate the maximum likelihood estimator,
is an established method for approximating maximum likelihood estimators 
in the statistical network analysis literature \citep{HuHa06}.
While able to provide accurate estimates of maximum likelihood estimators for complex models 
\citep[e.g.,][]{StScBoMo19,ScKrBu17}, 
a drawback of MCMLE,
and other simulation-based estimation methodology, 
is the computational burden which can scale with both the complexity of the model and the size of the network \citep{BaBrSl11}.  
In settings where the computation of the MCMLE is impractical, 
a computationally efficient alternative is provided via the maximum pseudolikelihood estimator (MPLE) \citep{Bj74},
whose application to social network analysis and to statistical network analysis dates back to \citet{StIk90}. 
Pseudolikelihood-based estimators have the following computational advantages:
\ben
\item Algorithms are generally deterministic and do not require simulation-based approximation schemes,
which aids in reproducibility of results; 
\item Algorithms are generally more scalable, 
relative to alternatives such as MCMLE and other simulation-based approximations, 
and are able to be parallelized to take advantage of larger  multicore computing infrastructures which are becoming increasingly common.  
\een

\hide{
Proposition \ref{prop:inference} establishes that $\bY$ is observable through $\bX$,
i.e., 
\beno
\mbP(Y_{i,j} = y_{i,j} \,|\, \bX = \bx, \bY_{-\{i,j\}} = \by_{-\{i,j\}})
\= 1,
\ee
when $y_{i,j} = \one(\norm{\bx_{i,j}}_1 \,>\, 0)$ and $\bY_{-\{i,j\}}$ is defined to be the 
($\binom{N}{2}$-$1$)-dimensional vector of edge variables in $\bY$ which excludes $Y_{i,j}$. 
As a result, 
if $(\bx, \by)$ is network concordant,
then  
\beno
\log \, \mbP(Y_{i,j} = y_{i,j} \,|\, \bX = \bx, \bY_{-\{i,j\}} = \by_{-\{i,j\}}) \= 0, 
&& \mbox{for all } \{i,j\} \subset \mN.
\ee 

The log-pseudolikelihood of \eqref{general model} can then be written down as 
\be
\label{eq:log-pseudo}
\pl(\nat; \bx, \by) 
\,\coloneqq \dsum_{\{i,j\} \subset \mN}  \dsum_{k=1}^K \, \log  
\sepmodel(X_{i,j}^{(k)} = x_{i,j}^{(k)} \, |\, \bX_{i,j}^{(-k)} = \bx_{i,j}^{(-k)}, \bY = \by), 
\ee
provided $(\bx, \by)$ is network concordant
and by exploiting the conditional independence properties implied by \eqref{general model}. 
We denote the set of maximum pseudolikelihood estimators of the data-generating parameter vector $\truth$ by 
\beno
\Mple &\coloneqq& \left\{ \nat \in \mbR^p \,: \, \pl(\nat; \bx, \by) = \sup\limits_{\nat^{\prime} \in \mbR^p} \, \pl(\nat^{\prime}; \bx, \by) \right\}.
\ee
Individual elements are referenced by $\mple \in \Mple$.
When it exists, 
the maximum pseudolikelihood estimator may or may not be unique.
However,
our theoretical results establish that all elements $\mple \in \Mple$ will all be within the same Euclidean distance to $\truth$.
The assumption that $(\bx, \by)$ is network concordant comes at no cost,
since $\bY$ is predictable through $\bX$,
as discussed above. 
The advantage of \eqref{eq:log-pseudo} is that 
the conditional probabilities 
$\sepmodel(X_{i,j}^{(k)} = x_{i,j}^{(k)} \, |\, \bX_{i,j}^{(-k)} = \bx_{i,j}^{(-k)}, \bY = \by)$ 
of edges in the multilayer network 
are often computationally tractable 
since the conditional distribution is a Bernoulli distribution when $Y_{i,j} = 1$,
and is a degenerate point mass at $0$ when $Y_{i,j} = 0$.

}

In this simulation, 
we consider the maximum pseudolikelihood estimator, denoted by $\mple$.
We conduct simulation studies to investigate the performance of the maximum pseudolikelihood estimator $\mple$ (MPLE),
supplementing the theoretical results established in Sections \ref{sec3} and \ref{sec:normal} for maximum likelihood estimators. 
As will be discussed later in Section \ref{sec:app}, we successfully reproduced the sufficient statistics using the MPLE in the application, suggesting that the MPLE for the multilayer network model solves a score equation similar to that of the MLE. This indicates that the MPLE serves as a close approximation and can be a good proxy for the MLE.
In section \ref{sec:sim_con}, 
we demonstrate the consistency results of Theorem \ref{thm1} 
in settings of different data-generating parameters and increasing model dimensions.
We conduct simulation studies of the multivariate normal approximation established by Theorem \ref{thm2} 
in Section \ref{sec:sim_norm}.  
Lastly, 
we discuss several testing procedures for selecting non-zero effects 
while controlling the false discovery rate (FDR) at a given family-wise significance level $\alpha$. 
 
In all simulation studies, 
we sample concordant multilayer networks $(\bX,\bY)$ from \eqref{general_model} 
with the maximum order of corss-layer interaction $H=2$:
\be
\label{eq:sim_model}
f(\bx, \nat) = \dprod_{\{i,j\} \subset \mN} \, \exp\left( \dsum_{k=1}^K\theta_{k}\,x_{i,j}^{(k)} + \dsum_{\substack{k < l}}^{K}\, \theta_{k,l} \, x_{i,j}^{(k)} \, x_{i,j}^{(l)} \right).
\ee
Unless otherwise specified, the basis network $\bY$ is generated from the Bernoulli random graph model.

\begin{figure}[t]
    \centering
	\includegraphics[width=3.5in, height=2in]{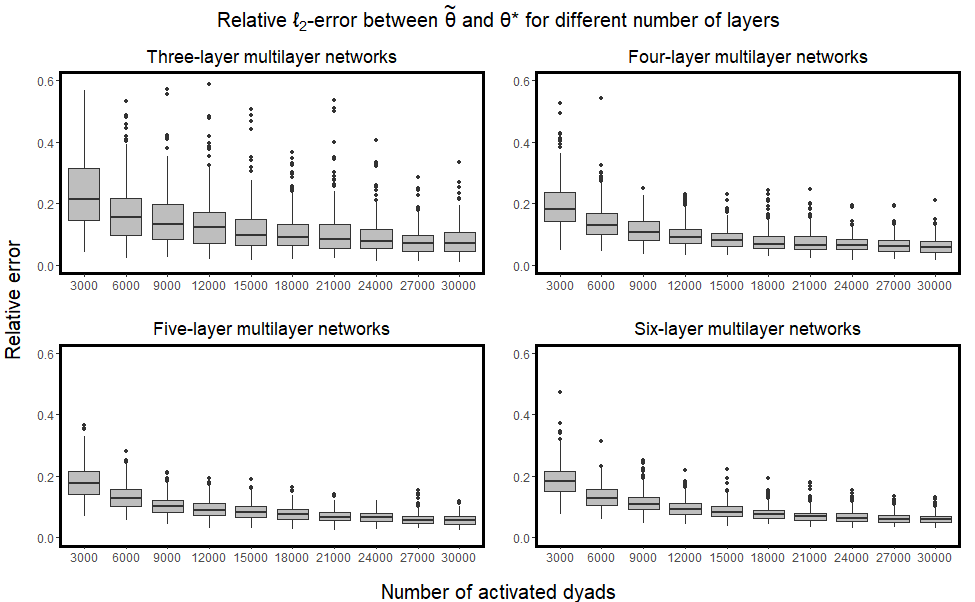}
	\caption{The relative $\ell_2$-errors between $\mple$ and $\truth$ decrease as the number of activated dyads increases. Each box is created by 250 replicates of multilayer networks.}
	\label{M_theta_error}
\end{figure}

\begin{figure}[t]
    \centering
	\includegraphics[width=0.9\textwidth]{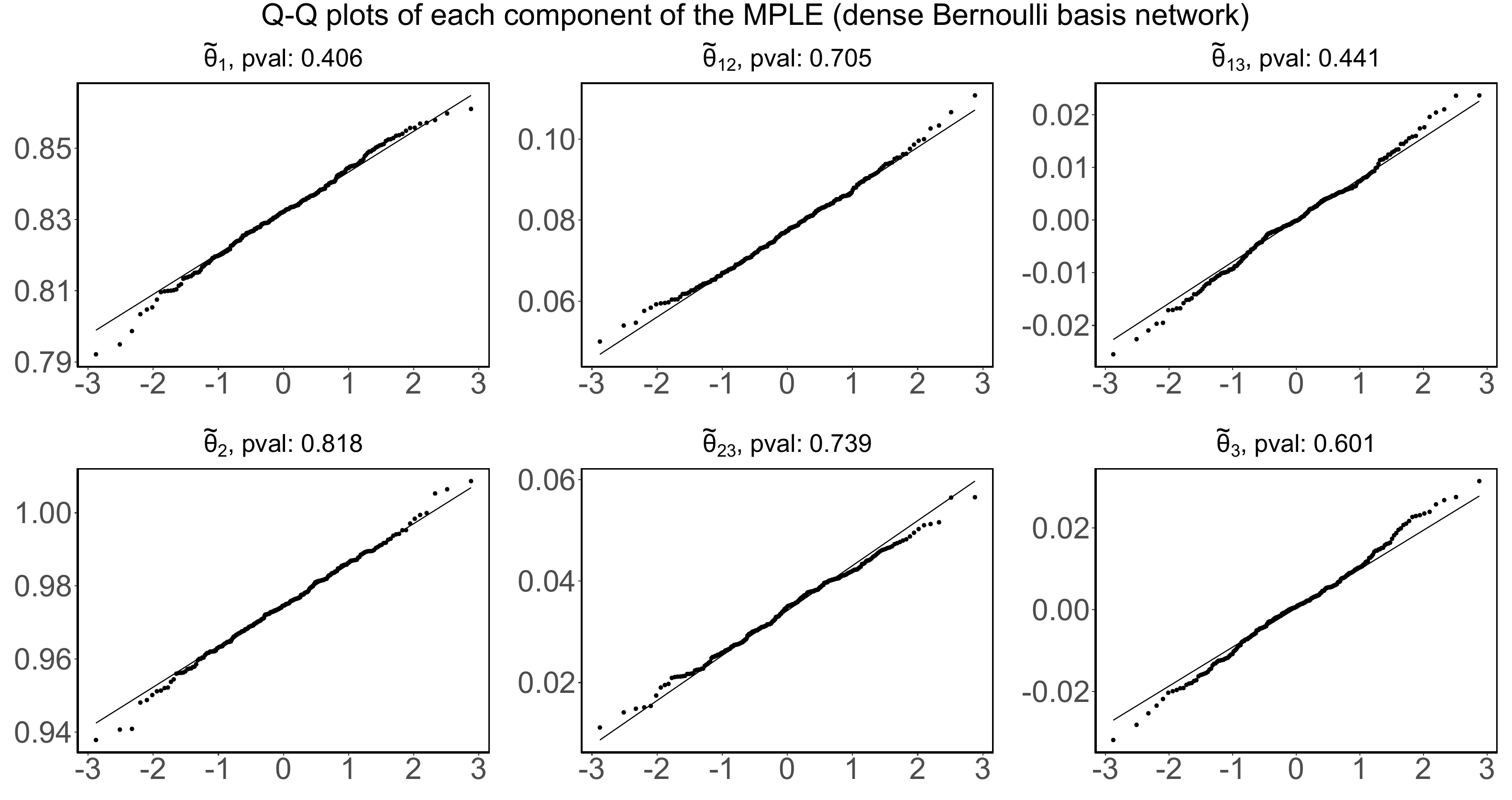}
	\caption{Q-Q plots and $p$-values of six components of $\mple$ estimated from 250 multilayer network samples at size 1000 on the dense Bernoulli basis network.}
	\label{qqplot_denseBer}
\end{figure}

\subsection{Consistency} 
\label{sec:sim_con}
The consistency is demonstrated through the decay of the relative $\ell_2$-errors between $\mple$ and the data-generating parameter $\truth$ as the expected number of activated dyads $\mbE\norm{\bY}_1$ increases. 
We generated $M = 250$ multilayer networks with $N=300$ nodes, using $M$ different data-generating parameters. We created these networks for each of ten evenly spaced numbers of activated dyads increasing from $3000$ to $30000$, and for four different numbers of layers increasing from $K=3$ to $6$. The model dimension increases from 6 to 21 as $K$ increases from 3 to 6.
For each number of activated dyads, 
number of layers $K$, 
and replicate, 
we sample a multilayer network $\bX$ from \eqref{general_model} using the specification in \eqref{eq:sim_model}
with the data-generating parameter vector $\truth$ populated by randomly selecting each component from the uniform distribution on $(-1,1)$.
We make the exception that components $\theta_{3}^\star$ and $\theta_{1,3}^\star$ are set to $0$.
In each replicate, 
we compute the maximum pseudolikelihood estimator. 
The results of this simulation study are given in Figure \ref{M_theta_error}, 
which shows the decay of the relative $\ell_2$-errors between $\mple$ and $\truth$ as the number of activated dyads increases in networks with different number of layers. 
The broad selection of data-generating parameter values on networks with increasing number of layers verifies that Theorem \ref{thm1} holds in many practical settings with increasing model dimensions.

\hide{
The decay rate of the relative $\ell_2$-error as the network size increases can be estimated by the slope of the OLS fitted line in the $\log$-$\log$ plot as shown in Figure \ref{Figure1}(b). The slope of the OLS line is $-1$, indicating a decay rate of order $1/N$, which follows the result in Theorem \ref{thm1}. We report the MPLE $\mple$ averaged from $M = 500$ samples with the selected data-generating parameter $\truth$ at network size 1000 in table \ref{t1}.

\begin{table}[t]
\begin{center}
\caption{\label{t1}Values of the data-generating parameter $\truth$ and mean of MPLE $\mple$ from 500 replications at network size 1000. Standard errors are in the parenthesis.}
\begin{tabular}{| c | c | c | c | c | c | c |} 
\hline
  & $\theta_{1}$&$\theta_{2}$&$\theta_{3}$& $\theta_{1,2}$ & $\theta_{1,3}$ & $\theta_{2,3}$ \\ 
\hline
$\truth$ & $-3$ & $-2$ & $-1$ & $.5$ & $0$ & $0$  \\
\hline
$\mple$ & $-2.999 \, (.02)$ & $-1.998 \, (.02)$ & $-.999 \, (.02)$ & $.499 \, (.02)$ & $-.001 \, (.02)$ & $-.002 \, (.02)$\\
\hline
\end{tabular}
\end{center}
\end{table}

}


\subsection{Multivariate normality and model selection} 
\label{sec:sim_norm}

\begin{table}[t]
\begin{center}
\caption{\label{fdr} False discovery rates of four procedures for detecting non-zero effects of 6 data-generating parameters 
($\truth_1$, $\truth_2$, $\truth_3$, $\truth_{4}$, $\truth_{5}$, $\truth_{6}$) estimated from 250 multilayer network samples at size 1000 on the dense Bernoulli basis network. All FDRs are smaller than $.05$.} 
\begin{tabular}{ l  r  r  r  r  r  r  } 
\hline
  Procedure & $\truth_1$ & $\truth_2$  & $\truth_3$ & $\truth_4$ & $\truth_5$ & $\truth_6$ \\ 
\hline
  Bonferroni    & .004  & .002 & .001 & .002 & .001 & .005  \\

  Benjamini-Hochberg   & .014 & .014 & .014 & .011 & .017 & .020  \\

 Hochberg &    .012 & .008 & .009 & .008 & .011 & .016 \\

 Holm  &   .010 & .008 & .006 & .008 & .007 &  .013 \\
\hline
\end{tabular}
\end{center}
\end{table}

As stated in Section \ref{sec:normal} and Theorem \ref{thm2}, 
the distribution of the maximum likelihood estimator $\mle$ 
converges in distribution to a multivariate normal distribution asymptotically. 
In order to study the quality of the normal approximation---especially for univariate testing 
which would be used for the false discovery rate control and model selection---we 
randomly select $6$ of the $250$ data-generating parameter vectors $\truth$
used to study the consistency results of Theorem \ref{thm1} 
in the simulation study conducted in Section \ref{sec:sim_con}.
We then generate $250$ replicates of multilayer network samples by each of these $6$ parameter vectors, using specification \eqref{eq:sim_model} on four basis network structures with the number of layers $K=3$:
the Bernoulli random graph model (dense and sparse), 
the stochastic block model, 
and the latent space model. 

The multivariate normality of $\mple$ passed Zhou-Shao's multivariate normal test \citep{Zhou13}, 
with $p$-values provided in the Appendix \ref{subsec:norm_sim} in the supplement to this paper.  
We visualize the marginal normality of individual component in $\mple$ with a dense Bernoulli basis network
in Figure \ref{qqplot_denseBer},
through Q-Q plots of the simulated maximum pseudolikelihood estimators. 
Univariate tests for normality failed to reject the null hypothesis that
each component of $\mple$ is marginally normal at a significance level of $.05$. 
Additional results studying the multivariate normality of $\mple$ on different basis network structures 
are provided in Appendix \ref{subsec:norm_sim} in the supplement to this paper.

\hide{
We first demonstrate through Q-Q plots in Figure \ref{Figure3} that each component of the MPLE $\mple$ follows a marginal normal distribution. 

\begin{table}
\begin{center}
\caption{\label{ZS-test} $p$-values of the Zhou-Shao's test for multivariate normality of $\mple$ for 6 data-generating parameters ($\truth_1$, $\truth_2$, $\truth_3$, $\truth_4$, $\truth_5$, $\truth_6$) estimated from 250 network samples at size 1000 on four basis network structures. All $p$-values are larger than .05. \s} 
\begin{tabular}{| c | c | c | c | c | c | c | } 
\hline
 Basis network model  & $\truth_1$ & $\truth_2$  & $\truth_3$ & $\truth_4$ & $\truth_5$ & $\truth_6$ \\ 
\hline
  Dense Bernoulli & .138 & .473 & .053 & .699 & .587 & .983  \\
\hline
 Sparse Bernoulli & .554 & .132 & .232 & .634 & .904 & .373  \\
\hline
 SBM & .65 & .891 & .982 & .975 & .871 & .674 \\
\hline
 LSM  & .859 & .831 & .5 & .227 & .613 & .409  \\
\hline
\end{tabular}
\end{center}
\end{table}
}
\hide{
We additionally performed marginal tests for normality for each component. In each case, the marginal test failed to reject the null hypotheses that $\mple_i$ is marginal normal at the same significance level, for each component $i = 1,\ldots, p$, the result of which is consistent with the $Q$-$Q$ plots. 
}

\hide{
\begin{figure}[t]
	\includegraphics[width=\textwidth]{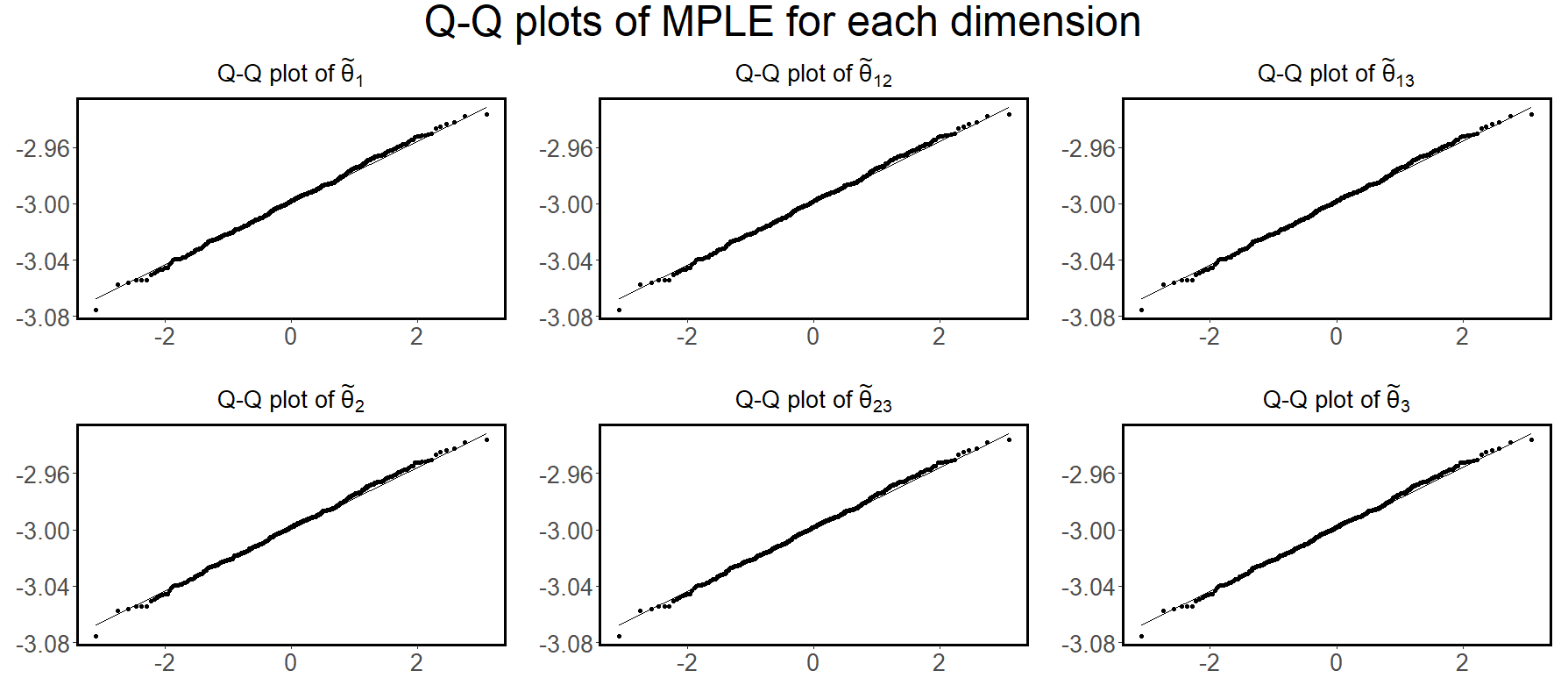}
	\caption{Q-Q plots of all dimensions of $\Tilde{\nat}$ show that the MPLEs fall on a straight line, indicating marginal normality of each dimension of $\mple$.}
	\label{Figure3}
\end{figure}

\begin{table}[t]
\begin{center}
\caption{\label{t3} Empirical FDR and power at significance level $.05$ for multiple testing of MPLE $\mple$ 
estimated from $500$ networks of size $1000$.  } 
\begin{tabular}{| c | c | c | c |} 
\hline
 Procedure  & Empirical FDR & FDR 95\% one-sided CI  & Empirical power \\ 
\hline
  Bonfferoni & .138 & 0.473 & .038 & .699 & .587 & .983  \\
\hline
  Benjamini-Hochberg's &  .018 & (0, .134) & 1  \\
\hline
 Hochberg's  & .016 & (0, .127) & 1 \\
\hline
 Holm's  & .013 &  (0, .114) & 1  \\
\hline
\end{tabular}
\end{center}
\end{table}
}

We then implement the multiple testing correction procedures of
Bonferroni, Benjamini-Hochberg, Hochberg, and Holm, for the 6 selected data-generating parameter vectors $\truth$ with 250 replicates to detect components that are significantly different from $0$ 
while controlling the false discovery rate (FDR) at a family-wise significance level of $\alpha = .05$---recall 
that $\theta_{1,3}^\star$ and $\theta_{3}^\star$ of $\truth$ are set to $0$ 
in each simulation replicate.  
We estimate the FDR of the four procedures by averaging the false discovery proportions from 250 replicates of each of the 6 randomly selected data-generating parameters $\truth$.
We provide the estimated FDRs for $\truth$ on a dense Bernoulli basis network in Table \ref{fdr}. 
In addition, we show the receiver operating characteristic (ROC) curves for $\mple$ estimating the 6 selected data-generating parameters in each of the subplot of Figure \ref{ROC}, on four basis network structures in Appendix \ref{more fdr} in the supplement to the paper. 
Simulation results suggest that the false discovery rate is controlled below the preset threshold $\alpha$. Different data-generating parameter values affect the trade-off between the sensitivity and the specificity of the model selection. In general, multilayer networks with a larger effective sample size lead to a larger area under the ROC curve which offers a tool to choose appropriate correction procedures and thresholds for model selection in different scenarios. 
Additional results on the false discovery rate with different basis network structures are provided in Appendix \ref{more fdr} in the supplement to the paper.

\section{Application} 
\label{sec:app}

\input{Lazega}

\section{Discussion} 
\label{sec:disc} 

In this work, 
we introduced a flexible class of statistical models for multilayer networks. 
Key to our approach lies in the integrative nature by which we establish our framework, 
extending arbitrary strictly positive probability distributions for single-layer networks 
to multilayer-network models through a novel separable framework with Markov random field specifications. 
We established the foundations for statistical inference through consistency and multivariate normality results,
the results of which have been demonstrated in simulation studies and in an application.  
The key assumption to our approach lies in the network separability assumption, 
which necessitates network dyads be conditionally independent given the basis network. 
This assumption may or may not be valid in practice,
which would necessitate the development of generalizations of the framework we established in this work 
through the relaxation of the conditional independence assumption.  
Such relaxations would result in more complex dependence structures,
requiring 
new and careful theoretical treatment in order to establish similar statistical foundations of models 
to the ones we have developed here,  
representing potential avenues for future research. 

\section*{Acknowledgements} 

Jonathan R. Stewart was supported by 
NSF award SES-2345043 and 
the Department of Defense Test Resource Management Center under contracts FA8075-18-D-0002 and FA8075-21-F-0074.

\bibliographystyle{agsm}
\bibliography{base}

\newpage

\begin{appendices}

\input{supplement.tex}
\end{appendices}

\label{last.page}

\end{document}

%% file: preamble.tex
\usepackage[mathscr]{eucal}
\usepackage{hyperref}
\usepackage{dlfltxbcodetips}
\usepackage{caption}
\usepackage{bbm}
\usepackage{mathtools}
\usepackage{subfigure}
\usepackage{amsfonts}
\usepackage{amsmath}
\usepackage{amssymb}
\usepackage{fancybox}
\usepackage{graphicx}
\usepackage{bm}
\usepackage{latexsym}

\newcommand{\Mle}{\widehat\bTheta}
\newcommand{\mle}{\widehat\btheta}
\newcommand{\mple}{\widetilde\btheta}
\newcommand{\Mple}{\widetilde\bTheta}
\newcommand{\truth}{\btheta^\star}

\newcommand{\given}{\,|\,}

\newcommand{\mA}{\mathscr{A}}
\newcommand{\mB}{\mathscr{B}}

\newcommand{\mD}{\mathscr{D}}
\newcommand{\mE}{\mathscr{E}}
\newcommand{\mF}{\mathscr{F}}

\newcommand{\mN}{\mathscr{N}}

\newcommand{\mP}{\mathcal{P}}

\newcommand{\mU}{\mathscr{U}}
\newcommand{\mV}{\mathscr{V}}

\newcommand{\mcF}{\mathcal{F}}

\newcommand{\mcI}{I}

\newcommand{\mcR}{\mathcal{R}}

\newcommand{\mbE}{\mathbb{E}}

\newcommand{\mbN}{\mathbb{N}}

\newcommand{\mbP}{\mathbb{P}}

\newcommand{\mbR}{\mathbb{R}}

\newcommand{\mbX}{\mathbb{X}}
\newcommand{\mbY}{\mathbb{Y}}

\newcommand{\bs}{\bm{s}}

\newcommand{\bu}{\bm{u}}
\newcommand{\bv}{\bm{v}}

\newcommand{\bx}{\bm{x}}
\newcommand{\by}{\bm{y}}

\newcommand{\bA}{\bm{A}}

\newcommand{\bH}{\bm{H}}
\newcommand{\bI}{\bm{I}}

\newcommand{\bR}{\bm{R}}
\newcommand{\bS}{\bm{S}}

\newcommand{\bW}{\bm{W}}
\newcommand{\bX}{\bm{X}}
\newcommand{\bY}{\bm{Y}}
\newcommand{\bZ}{\bm{Z}}

\newcommand{\bmu}{\mbox{\boldmath$\mu$}}

\newcommand{\bTheta}{\bm\Theta}

\newcommand{\btheta}{\boldsymbol{\theta}}

\newcommand{\one}{\mathbbm{1}}

\newcommand{\dsum}{\displaystyle\sum\limits}

\newcommand{\dprod}{\displaystyle\prod\limits}

\newcommand{\var}{\mathbb{V}}
\newcommand{\cov}{\mathbb{C}}
\newcommand{\norm}[1]{|\!|#1|\!|} 
\newcommand{\mnorm}[1]{|\!|\!|#1|\!|\!|}

\newcommand{\sepmodel}{\mbP_{\nat}}
\newcommand{\sepfam}{\mcF}
\newcommand{\gradL}{\nabla_{\nat}(\bx, \by)}
\newcommand{\GradL}{\nabla_{\nat}(\bX, \bY)}
\newcommand{\Lam}{\widetilde{\lambda}_{\min}^{\epsilon} \; \mbE \, \norm{\bY}_1}

\newcommand{\pl}{\widetilde\ell}

\newcommand{\btt}{\begin{box}}
\newcommand{\ett}{\end{box}}

\newcommand{\btheorem}{\begin{bclogo}[couleur={rgb:orange,0;yellow,0;white,1},arrondi=0.1,logo=\bcplume,ombre=true]{Theorem}}
\newcommand{\ettheorem}{\end{bclogo}}

\newcommand{\bst}{\begin{bclogo}[couleur={rgb:orange,1;yellow,1;white,0.5},arrondi=0.1,logo=\bcpanchant]}
\newcommand{\est}{\end{bclogo}}

\newcommand{\ba}{\begin{array}{llllllllll}}
\newcommand{\ea}{\end{array}}

\newcommand{\bea}{\begin{equation}\begin{array}{llllllllll}}
\newcommand{\eea}{\end{array}\end{equation}}

\newcommand{\be}{\begin{equation}\begin{array}{lllllllllllllllll}}
\newcommand{\beno}{\begin{equation}\begin{array}{lllllllllllll}\nonumber}
\newcommand{\ee}{\end{array}\end{equation}}

\newcommand{\bel}{\begin{equation}\begin{array}{lllllllllllll}\nonumber}
\newcommand{\eel}{\Box\end{array}\end{equation}}

\newcommand{\bi}{\begin{itemize}}
\newcommand{\ei}{\end{itemize}}

\newcommand{\ben}{\begin{enumerate}}
\newcommand{\een}{\end{enumerate}}

\newcommand{\bbq}{\begin{quote}\bf\em}
\newcommand{\ebq}{\end{quote}}

\renewcommand{\=}{&=&}

\newcommand{\hide}[1]{}
\newcommand{\ghost}[1]{}

\newcommand{\s}{\vspace{0.25cm}}

\newcounter{ex}

\newcounter{counterexample}

\newcounter{proposition}
\newenvironment{proposition}[1][]{\refstepcounter{proposition}\par\smallskip\indent%
\textbf{Proposition~\theproposition #1}\em \rmfamily}{}

\newcounter{ttproof}
\newcounter{pproof}
\newcommand{\pproof}{%
\addtocounter{pproof}{1}%
{\medskip}\textsc{Proof of Proposition}{}
}

\newcounter{ccproof}
\newcounter{ccsproof}
\newcounter{llproof}
\newcommand{\llproof}{%
\addtocounter{llproof}{1}%
{}\textsc{Proof of Lemma}{}
}

\newcounter{com}

\newcommand{\nat}{\btheta}

%% file: introduction_ml.tex
Multilayer networks have become a recent focal point of research in the field of statistical network analysis  
\citep[e.g.,][]{LeChLy20, CaGo20, arroyo_multi, KrKoMa20, Mac21, chen2022, sosa2022, huang2022},
arising in applications where a common set of elements in a population
interact through multiple modes or relationships with other elements in the population.
A prototypical example in the literature might be the Lazega law firm network \citep{Lazega2001}, 
in which attorneys are linked through various forms of interaction,
such as advice seeking, friendship, collaboration, etc.,
each of which would form a distinct layer in the multilayer network  \citep{KrKoMa20}.  
In essence, a multilayer network is a composite structure, where each layer captures a specific type of interaction or relationship between the same set of elements.

Edges in one layer of the multilayer network may depend on edges in other layers, 
creating what is known as cross-layer dependence. 
Understanding the drivers of edge formation in multilayer networks requires 
learning the dependence structures across these layers. 
A key challenge lies in the fact that the cross-layer dependence can be highly varied and complex,
and the development of statistical models with theoretical guarantees for network data with dependent edges 
is challenging.  
Current methodological frameworks for multilayer networks can be broadly categorized into two main groups:
\ben
\item Statistical models equipped with theoretical guarantees often rely on latent variable constructions \cite[e.g.,][]{Mac21, arroyo_multi, huang2022}. These models typically assume conditional independence of edges given the latent variables, following standard practices within the field.
\item Statistical models that do not provide formal theoretical guarantees \citep[e.g.,][]{CaGo20, KrKoMa20}. Instead, these methods extend existing approaches by explicitly allowing for edge dependencies, thereby relaxing the conditional independence assumptions present in the first class of models.
\een
In this work, 
we address a critical gap in the literature by introducing a separable multilayer network modeling framework for multilayer networks.
Our approach not only accommodates dependent edges but also provides theoretical guarantees for both estimation and inference without relying on any latent variables.
Specifically, we extend single-layer network models to the multilayer setting, with a central focus on identifying and understanding cross-layer dependence structures.
A key advantage of our proposed framework is that we are able to distinguish the network formation process  
from the layer formation process.
This allows us to create a wide range of novel multilayer network models derived from established single-layer network models, such as
exponential-family random graph models, stochastic block models, and latent space models.
By employing Markov random field specifications, we develop adaptable and comprehensive models to capture cross-layer dependencies in multilayer networks. 
As a result, 
our framework jointly models both network structures and cross-layer dependence, thus enabling any single-layer network model to be extended to the multilayer setting. 
Our main contributions in this work include: 
\ben
\item Introducing a novel framework for modeling 
cross-layer dependence in multilayer networks 
that synchronizes with current network models in the literature.
\item Deriving non-asymptotic theoretical guarantees in scenarios where the number of parameters tends to infinity, 
which establishes bounds on the:  
\ben
\item Statistical error of maximum likelihood estimators.  
\item Error of the multivariate normal approximation of estimators.
\een
\item Elaborating a model selection algorithm which controls the false discovery rate. 
\een

%% file: general_model.tex
We specify probability distributions 
on a double of networks $(\bX, \bY)$, 
where $\bY$ will represent the network formation process,
which we will call the {\it basis network},  
and $\bX$ will represent the realized multilayer network. 
We assume that $\bY \in \mbY \coloneqq \{0, 1\}^{\binom{N}{2}}$ 
is an undirected, single-layer network defined on the set of nodes $\mN$ 
where,
for all $ \{i,j\} \subset \mN$,  
\beno
Y_{i,j}
\= \begin{cases}
1 & \mbox{nodes } i \mbox{ and } j \mbox{ are connected in the basis network} \\ 
0 & \mbox{otherwise} 
\end{cases},
\ee 
making the usual conventions for undirected networks mentioned previously.
For $(\bX, \bY)$, we consider semi-parametric families of probability distributions $\sepfam \coloneqq \{\sepmodel : \nat \in \mbR^p\}$ which are absolutely continuous with respect to a $\sigma$-finite measure $\nu$ defined on $\mP(\mbX \times \mbY)$,
where $\mP(\mbX \times \mbY)$ is the power set of $\mbX \times \mbY$. 
Typically,
$\nu$ will be the counting measure.  
We say the probability mass function $\sepmodel \in \sepfam$ defines a {\it separable multilayer network model} if $\sepmodel$ admits the form: 
\be
\label{general_model}
\sepmodel(\{(\bx, \by)\}) 
\= f(\bx,\nat) \; g(\by) \; h(\bx, \by) \; \psi(\nat, \by),
&&  (\bx, \by) \in \mbX \times \mbY, 
\ee
where 
\bi
\item $f : \mbX \times \mbR^p \mapsto (0, 1)$ is given by 
\beno
f(\bx, \nat)
&= \, \dprod_{\{i,j\} \subset \mN} \,
\exp&\left(\dsum_{k=1}^K\,\theta_{k}\,x_{i,j}^{(k)} +  \dsum_{\substack{k < l}}^{K}\, \theta_{k,l} \, x_{i,j}^{(k)} \, x_{i,j}^{(l)} + \ldots \right.\\[20pt]
&& \left.+ \dsum_{k_1 < \,\ldots\,< k_H}^{K}\, \theta_{k_1,k_2,\ldots,k_H} \, x_{i,j}^{(k_1)} \cdots\, x_{i,j}^{(k_H)} \right),
\ee
where $H \le K$ is the highest order of cross-layer interactions included in the model.
We write $\theta_{k_1,k_2,\ldots,k_h}$ to reference the $h$-order interaction parameter 
for the interaction term among layers $\{k_1, \ldots, k_h\} \subseteq \{1, \ldots, K\}$. \s  
\item $g : \mbY \mapsto (0, 1)$ is the marginal probability mass function of $\bY$ 
and is assumed to be strictly positive on $\mbY$. \s
\item $h : \mbX \times \mbY \mapsto \{0, 1\}$ is given by 
\beno
h(\bx, \by) 
\= \dprod_{\{i,j\} \subset \mN} \, 
\one(\norm{\bx_{i,j}}_1 > 0)^{y_{i,j}} \; \one(\norm{\bx_{i,j}}_1 = 0)^{1 - y_{i,j}}, 
\ee 
where $\bx_{i,j} = (x_{i,j}^{(1)}, \ldots, x_{i,j}^{(K)}) \in \mbX_{i,j}$ ($\{i,j\} \subset \mN$). \s   
\item $\psi : \bTheta \times \mbY \mapsto (0, \infty)$ is defined by 
\beno 
\psi(\nat, \by) 
\= \left[ \, \dsum_{\bx \in \mbX} \, f(\bx, \nat) \, h(\bx, \by) \right]^{-1}, 
\ee
ensuring \eqref{general_model}
will be a valid probability mass function for $(\bX, \bY)$.
\ei

The notation $\sepmodel(\{(\bx, \by)\})$ is well-defined for each pair $(\bx, \by) \in \mbX \times \mbY$, 
as $\sepmodel$ is a probability measure defined on $\mP(\mbX \times \mbY)$. 
Frequently, 
we will  write the probability expressions
$\sepmodel(\bX = \bx, \bY = \by)$ for the joint probability of $\{(\bx, \by)\}$, 
and $\sepmodel(\bX = \bx \,|\, \bY = \by)$ for the conditional probability of the event $\bX = \bx$ 
conditional on the event $\bY = \by$. 
We denote the data-generating parameter vector by $\truth \in \mbR^p$,
and the corresponding probability measure and expectation operator  
by $\mbP \equiv \mbP_{\truth}$ and $\mbE \equiv \mbE_{\truth}$, respectively.

The specification in equation \eqref{general_model} 
separates the network formation process $\bY$,
specified by $g(\by)$, 
from the layer formation process, 
specified by $f(\bx, \nat)$. 
The two are joined by the function $h(\bx, \by)$,
which ensures $\norm{\bx_{i,j}}_1 = 0$ 
whenever $Y_{i,j} = 0$
and $\norm{\bx_{i,j}}_1 > 0$ whenever $Y_{i,j} = 1$, 
as we allow edges between nodes $i \in \mN$ and $j \in \mN$ in $\bX$ if and only if $Y_{i,j} = 1$.
We call dyads $\{i, j\} \subset \mN$ with $Y_{i,j} = 1$ {\it activated dyads}, and a pair $(\bx, \by) \in \mbX \times \mbY$ that satisfies $h(\bx, \by) =1$ is said to be a {\it concordant} pair.
We will only focus on concordant pairs of multilayer networks since our modeling framework guarantees the recovery of the basis network $\bY$ given an observation of $\bX$, a point that will be made clear shortly in Proposition \ref{prop:inference}.
The function $\psi(\nat, \by)$ ensures the resulting product of functions will be a valid probability mass function, 
and it has less of a direct role in modeling the cross-layer dependence,
essentially fulfilling the role of a normalizing constant for the conditional probability distribution of $\bX$ given $\bY$,
as derived in Proposition \ref{prop:inference}. 
Such specifications have the advantage of being able to specify the network formation process separately 
from the process that populates the layers of activated dyads,
thus modeling the cross-layer dependence conditional on $\bY$. 
\begin{figure}[t]
\centering 
  \includegraphics[width=4in, height = 2.5in]{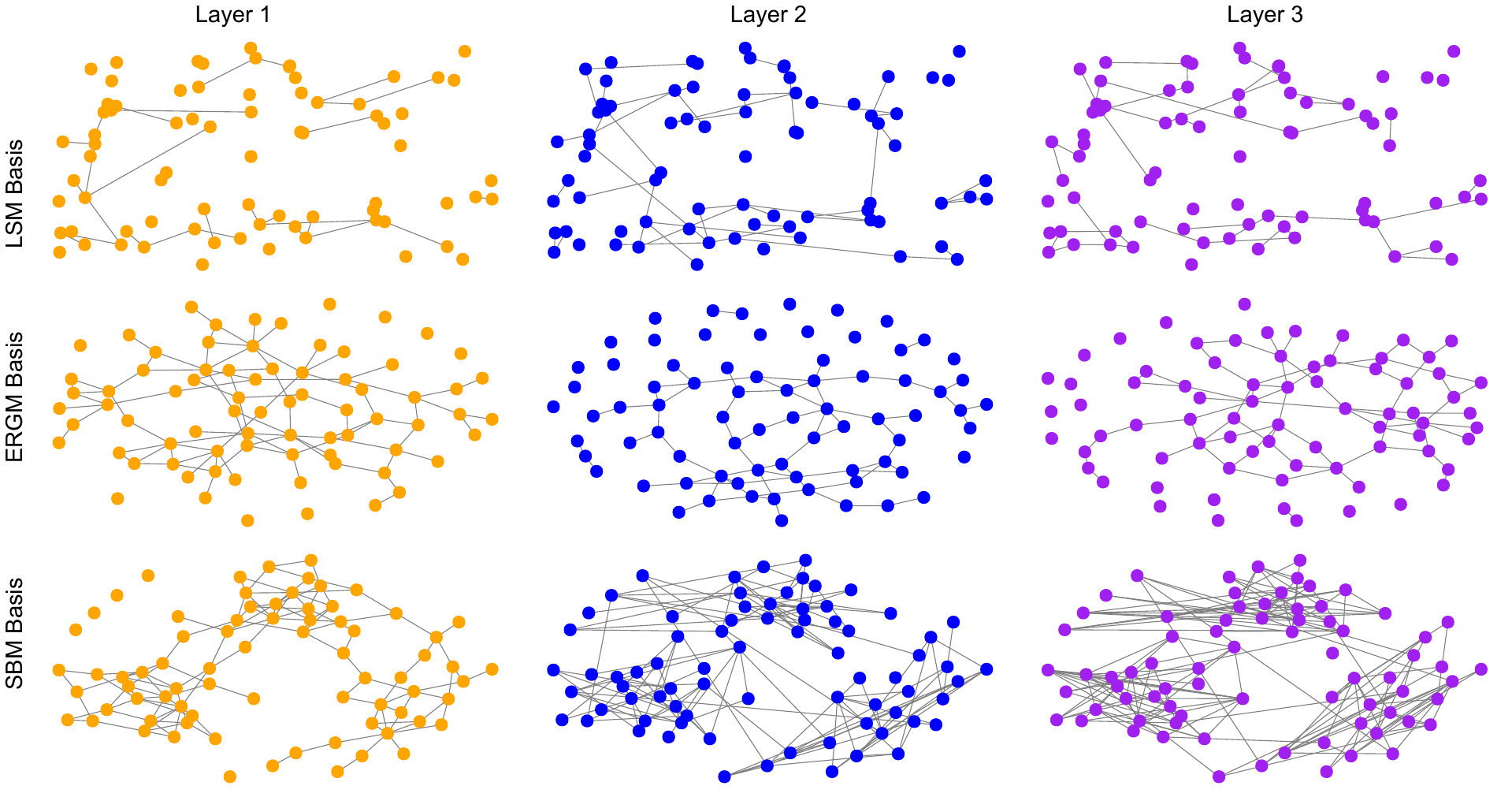}
  \caption{Multilayer networks specified by three different basis network structures: the latent space model (LSM), the exponential random graph model (ERGM), and the stochastic block model (SBM).}
  \label{fig:nets}
\end{figure}
To illustrate the flexibility and generality of \eqref{general_model},
observe that $g(\by)$ is allowed to  be any probability mass function for a single-layer network $\bY$ 
(e.g., exponential-family random graph model, stochastic block model, latent space model),
provided $g(\by) > 0$ for all $\by \in \mbY$. 
To illustrate this point, 
Figure \ref{fig:nets} displays various multilayer networks with $K = 3$ layers where the basis network 
is specified via three different models,
demonstrating that our modeling framework is capable of ensuring that the multilayer network 
respects structural properties of the underlying basis network. 
\hide{
The basis network can be thought of as a network describing the fundamental connection of two nodes,
which then propagates into specific realizations of connections in the multilayer network.
}
We view our framework as semi-parametric as $g(\by)$ need not assume a specific parametric form. 
Moreover, 
our framework can be viewed as non-parametric 
when the maximal possible order of interaction terms are included in \eqref{general_model},
a point on which we further elaborate later. 
An important feature of our framework lies in the fact that 
the choice of the probability distribution for the network formation process does not directly 
influence the estimation for the cross-layer dependence structure,
i.e.,
the choice of $g(\by)$ does not directly influence estimation for $\truth$. 
Proposition \ref{prop:inference} demonstrates this point in the case of likelihood-based inference.

%% file: prop_inference.tex
\begin{proposition}
\label{prop:inference}
Let $\{\sepmodel : \nat \in \mbR^p\}$ satisfy \eqref{general_model}. 
Then the following hold:   
\ben
\item 
For each $\bx \in \mbX$, 
$\bY = \by$ ($\sepmodel$-a.s.) for one and only one $\by \in \mbY$. \s  
\item $\bY$ is predictable via $\bX$, 
i.e.,
for each $\bx \in \mbX$,  
$\mbP_{\nat}(\bY = \by \,|\, \bX = \bx) = 1$ where 
\beno
y_{i,j} \= \one(\norm{\bx_{i,j}}_1 \,>\, 0), 
&& \{i,j\} \subset \mN. 
\ee  
\item For all $(\bx, \by) \in \mbX \times \mbY$ with $h(\bx, \by) = 1$, 
\beno
\log \, \sepmodel(\bX = \bx, \bY = \by)
\;=\; \log \, \sepmodel(\bX = \bx \mid \bY = \by) + \log \, g(\by),
\ee
where $\sepmodel(\bX = \bx \,|\, \bY = \by)$ belongs to a minimal exponential family 
with natural parameter vector $\nat \in \mbR^p$ and is given by 
\beno
\sepmodel(\bX = \bx \mid \bY = \by)
\;=\; \exp(\log f(\bx, \nat) + \log \psi(\nat, \by)).  
\ee
\een 
\end{proposition}

%% file: example_ml2.tex
We illustrate cross-layer dependence among layers in our modeling framework
by considering a separable multilayer network model using the Markov random field specification 
for $f(\bx, \nat)$ given in the previous section and maximal interaction term $H = 2$:
\be
\label{pairwise}
f(\bx, \nat)
\= \dprod_{\{i,j\} \subset \mN} \, 
\exp\left( \, \dsum_{k=1}^{K}\,\theta_{k} \,  x_{i,j}^{(k)} +  
\dsum_{k<l}^{K} \, \theta_{k,l} \,  x_{i,j}^{(k)} \, x_{i,j}^{(l)} \right).
\ee
The dimension of the parameter vector $\nat$ is $\dim(\nat) = K + \binom{K}{2}$,
with $K$ parameters governing the single-layer effects for the $K$ layers and $\binom{K}{2}$ combinations 
of layers to form the pairwise interactions for the cross-layer dependence effects.  

Define the ($K$-$1$)-dimensional vector 
$X_{i,j}^{(-k)} \coloneqq (X_{i,j}^{(l)} : l \in \{1, \ldots, K\} \setminus \{k\})$ 
to be the 
vector of edge variables in $\bX_{i,j}$ which excludes the edge variable $X_{i,j}^{(k)}$,
i.e.,
excluding the edge variable between nodes $i$ and $j$ in layer $k$.  
The conditional log-odds of edge $X_{i,j}^{(k)}$ takes the form:    
\beno
\scalebox{0.95}{$\log \, \dfrac{\mbP(X_{i,j}^{(k)} = 1 \,|\, \bX_{i,j}^{(-k)} = \bx_{i,j}^{(-k)}, Y_{i,j} = 1)}
{\mbP(X_{i,j}^{(k)} = 0 \,|\, \bX_{i,j}^{(-k)} = \bx_{i,j}^{(-k)}, Y_{i,j} = 1)}$} 
= \begin{cases}
\theta_{k} +  \dsum_{l \neq k}^{K}\, \theta_{k,l} \, x_{i,j}^{(l)}, & \norm{\bx_{i,j}^{(-k)}}_1 > 0 \\ 
+\infty, & \norm{\bx_{i,j}^{(-k)}}_1 = 0  
\end{cases}.
\ee
A primary advantage and motivation
of using a parametric Markov random field specification for $f(\bx, \nat)$ lies in the interpretability of the model. 
An effective approach to analyzing and understanding marginal network effects in such specifications is to study 
conditional log-odds of edges under different conditioning statements \citep[e.g.,][]{StScBoMo19}.
By the form of $h(\bx, \by)$,
when $Y_{i,j} = 1$,
we require $\norm{\bx_{i,j}}_1 > 0$,
meaning nodes $i$ and $j$ must have at least one connection in $\bX$. 
This is seen through the log-odds formula above,
where the log-odds of edge $X_{i,j}^{(k)}$ is equal to $+\infty$ when $\norm{\bx_{i,j}^{(-k)}}_1 = 0$. 
In contrast, 
when $\norm{\bx_{i,j}^{(-k)}}_1 > 0$,
the constraint $\norm{\bx_{i,j}}_1 > 0$ is already satisfied,
and the log-odds of edge $X_{i,j}^{(k)}$ depends on the layer specific parameter $\theta_k$, 
as well as the pairwise interaction effects where edges present in other layers $l \in \{1, \ldots, K\} \setminus \{k\}$ 
can influence the likelihood of the edge $X_{i,j}^{(k)}$ depending on the signs and magnitudes of the pairwise interaction 
parameters $\theta_{k,l}$ ($\{k,l\} \subseteq \{1, \ldots K\})$.

%% file: estimation_setup.tex
For separable multilayer network models satisfying \eqref{general_model}, 
Proposition \ref{prop:inference} establishes that the log-likelihood function takes the form 
\be
\label{loglikelihood}
\ell(\nat; \bx, \by) 
&\coloneqq& \log \,\sepmodel(\bX = \bx, \bY = \by) \s \\
\= \log \, \mbP_{\nat}(\bX = \bx \,|\, \bY = \by) + \log \, g(\by).
\ee
Given an observation $\bx \in \mbX$ of the multilayer network $\bX$,  
and therefore an observation $\by \in \mbY$ of $\bY$ by Proposition \ref{prop:inference},
we denote the set of maximum likelihood estimators by 
\beno
\Mle &\coloneqq& \left\{ \nat \in \mbR^p 
\,:\, \ell(\nat; \bx, \by) = \sup\limits_{\nat^{\prime} \in \mbR^p} \, \ell(\nat^{\prime}; \bx, \by) \right\},
\ee 
and reference individual elements of the set by $\mle \in \Mle$.
As Proposition \ref{prop:inference} establishes 
$\log \, \mbP_{\nat}(\bX = \bx \,|\, \bY = \by)$ to be a minimal, 
and by construction regular, 
exponential family, 
$|\Mle| \in \{0, 1\}$, 
i.e., 
when the maximum likelihood estimator exists, 
the set $\Mle$ will contain a unique element when non-empty
\citep[Proposition 3.11, pp. 32--33,][]{Su19}.
As seen from the forms of $\ell(\nat; \bx, \by)$ given above, 
the gradients and Hessians of the log-likelihood equations 
do not directly depend on $g(\by)$. 
However, 
the following lemma shows how theoretical guarantees for estimators of $\truth$ will be indirectly influenced by the choice of $g(\by)$.

%% file: min_eig_lemma.tex
\begin{lemma}
\label{lem:min-eig} 
Consider a family $\{\sepmodel : \nat \in \mbR^p\}$ of separable multilayer network models satisfying \eqref{general_model} 
and an observation $\bx \in \mbX$ of $\bX$.
Let $(\bx, \by)$ be the concordant pair where $\by$ is given by Proposition \ref{prop:inference}. 
Define,
for each pair of nodes $\{i,j\} \subset \mN$,
\beno
L_{i,j}(\nat, \bx_{i,j}, \by)
&\coloneqq& \log \, \sepmodel(\bX_{i,j} = \bx_{i,j} \,|\, \bY = \by).
\ee
Then there exists a $p \times p$ matrix 
$\mcI(\nat)$ such that 
\beno
\mbE\left[ - \nabla_{\nat}^2 L_{i,j}(\nat, \bX_{i,j}, \bY) \,|\, \bY = \by \right] 
\=
\begin{cases}
    \mcI(\nat) & Y_{i,j} = 1 \\
    \bm{0}_{p,p} & Y_{i,j} = 0,
\end{cases}
\ee
for all $\{i,j\} \subset \mN$,  
where $\bm{0}_{p,p}$ is the $p\times p$ matrix with all $0$ entries, 
and 
\beno 
&& \lambda_{\min}(-\mbE \, \nabla_{\nat}^2 \, \ell(\nat; \bX, \bY)) 
\= \lambda_{\min}(\mcI(\nat)) \, \mbE \, \norm{\bY}_1 \s \\
&& \lambda_{\max}(-\mbE \, \nabla_{\nat}^2 \, \ell(\nat; \bX, \bY)) 
\= \lambda_{\max}(\mcI(\nat)) \, \mbE \, \norm{\bY}_1,  
\ee
where $\lambda_{\min}(\bA)$ and $\lambda_{\max}(\bA)$ are the smallest and the largest eigenvalue of matrix $\bA \in \mbR^{p \times p}$, respectively.
\end{lemma}

%% file: thm1.tex
\begin{theorem}
\label{thm1}
Consider a multilayer network model following the form of equation \eqref{general_model} and is 
defined on a set of $N \geq 3$ nodes. If Assumptions \ref{assump1}, \ref{assump2}, and \ref{assump3} are satisfied,
there exists $N_0 \geq 3$ such that, 
for all $N \geq N_0$,
with probability at least $1 - \exp\,(-2\,p) \, - \, (\mbE \, \norm{\bY}_1)^{-1}$,
the set $\Mle$ is non-empty and the unique element $\mle \in \Mle$ satisfies  
\be
\label{consistency}
\norm{\mle - \truth}_2 &\leq&
C \; \dfrac{\sqrt{\widetilde{\lambda}_{\max}^{\star}} }{\widetilde{\lambda}_{\min}^{\epsilon}} \; \sqrt{\dfrac{p}{\mbE \norm{\bY}_1}},
\ee
where $C > 0$ is a constant independent of $N$ and $p$. 

\end{theorem}

The results of Theorem \ref{thm1} establish a few key facts concerning statistical estimation 
of the data-generating parameter vector $\truth$. 
First, 
we can view the quantity $\widetilde{\lambda}_{\min}^{\epsilon}\, \sqrt{\mbE \, \norm{\bY}_1 \, / \, \widetilde{\lambda}_{\max}^{\star}}$ 
as the effective sample size in order to compare our results to 
classical settings with independent and identically distributed data. 
The effective sample size, 
together with the dimension of the model $p$, 
helps to determine the rate of convergence (with respect to the Euclidean distance) 
of maximum likelihood estimators. 
As previously mentioned, 
the quantity $\mbE \norm{\bY}_1$ 
is determined by properties of $g(\by)$,
the marginal probability mass function of $\bY$. 
While the specification of $g(\by)$ does not directly influence the estimation algorithm, 
the statistical guarantees of estimators will depend on $g(\by)$ producing enough activated dyads and not possessing overly strong 
dependence among edges in the single-layer basis network $\bY$ (Assumption \ref{assump1}).
The requirement (Assumption \ref{assump3}) that the right-hand side of the bounds in Theorem \ref{thm1} tends to $0$ as $N \to \infty$ 
ensures that all regularity assumptions remain valid. 
Namely, 
key to our approach lies in the ability to control minimum eigenvalues of matrices 
$\mcI(\nat)$
in a neighborhood of the data-generating parameter vector $\truth$. 
The condition that the bounds tend to $0$ ensures that it is sufficient to control the smallest eigenvalue
in a bounded set,
i.e.,
we may let $\epsilon$ be fixed independent of $N$ and $p$, 
and moreover, 
to ensure consistency in the sense that $\norm{\mle - \truth}_2 \to 0$ with probability approaching $1$ as $N, p \to \infty$.

\begin{corollary}
\label{corollary}
Under the assumptions of Theorem \ref{thm1}, and in the case that the parameter dimension $p$ is fixed, there exists $N_0 \geq 3$ such that, 
for all $N \geq N_0$ and $\alpha_N \in (2\,(\mbE\norm{\bY}_1)^{-1} \, , \, 1/2)$,
with probability at least $1 \, - \, \alpha_N$,
the set $\Mle$ is non-empty and the unique element $\mle \in \Mle$ satisfies  
\beno
\norm{\mle - \truth}_2 &\leq&
C \; \left| \, \log\, \left(\dfrac{\alpha_N}{2} \right) \, \right|  \; \dfrac{\sqrt{\widetilde{\lambda}_{\max}^{\star}} }{\widetilde{\lambda}_{\min}^{\epsilon}} \; \sqrt{\dfrac{1}{\mbE \norm{\bY}_1}},
\ee
where $C>0$ is a constant independent of $N$.
\end{corollary}

The corollary builds on the consistency result established in Theorem \ref{thm1} by providing a similar bound in the situation that the parameter dimension $p$ remains fixed. This simplifies the convergence rate by removing the dependence of $p$ in the error term, which can yield sharper asymptotic guarantees. The introduction of the probability bound $\alpha_N$ offers an explicit control over the confidence level for the estimate's accuracy, which improves interpretability and practical applicability in finite samples. The bound on $\norm{\mle - \truth}_2$ in Corollary \ref{corollary} now depends logarithmically on $\alpha_N$, introducing a trade-off between the confidence level and the convergence rate. While the key dependencies remain on the effective sample size $\widetilde{\lambda}_{\min}^{\epsilon}\, \sqrt{\mbE \, \norm{\bY}_1 \, / \, \widetilde{\lambda}_{\max}^{\star}}$ as in Theorem \ref{thm1}, Corollary \ref{corollary} provides a useful refinement of the consistency result when the model's dimensionality is constrained.

We next present that the upper bound in Theorem \ref{thm1} is minimax optimal up to a constant. Define the minimax risk to be
\beno
\mcR_N &\coloneqq& \inf\limits_{\mle}\;\sup\limits_{\nat \in \mbR^p} \; \mbE_{\nat}\,\norm{\mle - \nat}_2,
\ee
and
\beno
\widetilde{\lambda}_{\max}^{\epsilon}
&\coloneqq & \sup\limits_{\nat \in \mB_2(\truth, \epsilon)} \, \lambda_{\max} \, (\mcI(\nat)),
\ee
where $\lambda_{\max}(\bA)$ is the largest eigenvalue of matrix $\bA \in \mbR^{p \times p}$.

\begin{theorem}
\label{thm:minimax}
Consider a separable multilayer network model following the form of equation \eqref{general_model} and is 
defined on a set of $N \geq 3$ nodes. If Assumptions \ref{assump1}, \ref{assump2} and \ref{assump3} are satisfied,
there exists a constant $C > 0$ independent of $N$ and $p$, such that, 
the lower bound of the minimax risk $\mcR_N$ satisfies
\beno
\mcR_N
& \geq &
C \, \dfrac{\widetilde{\lambda}_{\min}^{\epsilon}}{\widetilde\lambda_{\max}^{\epsilon}} \, \dfrac{\sqrt{\widetilde{\lambda}_{\max}^{\star}}}{\widetilde{\lambda}_{\min}^{\epsilon}} \, \sqrt{\dfrac{p}{\mbE\,\norm{\bY}_1} }.
\ee
\end{theorem}

Theorem \ref{thm:minimax} establishes a lower bound for the minimax risk $\mcR_N$, differing from the upper bound of the $\ell_2$-error for the maximum likelihood estimator in Theorem \ref{thm1} by a factor of $\widetilde{\lambda}_{\min}^{\epsilon} \, / \, \widetilde\lambda_{\max}^{\epsilon}$. Building on this result, we establish conditions for the minimax optimality of the maximum likelihood estimators in Corollary \ref{cor:minimax}.

\begin{corollary}
    \label{cor:minimax}
    Under the assumptions of Theorem \ref{thm:minimax} and the assumption that
    \be
    \label{assump:minimax}
    \widetilde\lambda_{\max}^{\epsilon}= O\, \left(\widetilde{\lambda}_{\min}^{\epsilon}\right),
    \ee
    the maximum likelihood estimator $\mle$ achieves the minimax rate of convergence, in the sense that the upper bound on the $\ell_2$-error of $\mle$ given in Theorem \ref{thm1} matches the lower bound of the minimax risk $\mcR_N$ in Theorem \ref{thm:minimax}, up to a constant.
\end{corollary}

The condition in \eqref{assump:minimax} ensures that the rate of convergence for the maximum likelihood estimator achieves the minimax optimality by imposing a more direct and stringent relationship between the minimum and maximum eigenvalues than that required by Assumption \ref{assump3}. The control on the minimum and maximum eigenvalues for high-dimensional graphical models are common \citep[e.g.,][]{witten14, Zhao06, RaWaLa10}, ensuring that the minimum and maximum eigenvalues of the information matrices within a neighborhood of the data-generating parameter are bounded away from 0 and bounded from above, respectively, and do not diverge relative to one another.

\hide{
\begin{corollary}[version 2]
Under the assumptions of Theorem \ref{thm1}, and in the case that the parameter dimension $p$ is fixed, there exists $N_0 \geq 3$ such that, 
for all $N \geq N_0$ and $\gamma > 0$,
with probability at least $1 \, - \, \exp \, (-\gamma \, p ) \, - \, (\mbE \, \norm{\bY}_1)^{-1}$,
the set $\Mle$ is non-empty and the unique element $\mle \in \Mle$ satisfies  
\beno
\norm{\mle - \truth}_2 &\leq&
\beta \, (\gamma) \; \dfrac{\sqrt{\widetilde{\lambda}_{\max}^{\star}} }{\widetilde{\lambda}_{\min}^{\epsilon^\star}} \; \sqrt{\dfrac{p}{\mbE \norm{\bY}_1}},
\ee
where $\beta : (0 \, , \, \infty) \mapsto (C_2 \; + \; \sqrt{C_2^2 \; + \; 72 \, C_2} \; / \; 2 \, , \,  \infty)$ is an unbounded increasing function of $\gamma$ given by
\beno
\beta\, (\gamma) \= \dfrac{C_2 \; + \;  C_2 \, \gamma \; + \; \sqrt{(C_2\; + \;  C_2 \, \gamma )^2 \; + \; 72\,(C_2\; + \; C_2 \, \gamma )}}{2},
\ee
where $C_2 > 0$ is a constant independent of $N$ and $\gamma$.
\end{corollary}
}

%% file: thm2.tex
In the following, 
$\bZ$ will denote a standard multivariate normal random vector,
i.e.,
with mean vector equal to the zero vector and covariance matrix equal to the identity matrix
(each of appropriate dimension), 
and $\Phi$ will denote the corresponding probability measure.

\begin{theorem}
\label{thm2}
Consider a separable multilayer network model following the form of equation \eqref{general_model} and is 
defined on a set of $N \geq 3$ nodes.
There exists $N_0 \geq 3$ such that,
for all $N \geq N_0$ and any measurable convex set $\mA \subseteq \mbR^p$, 
the error of the multivariate normal approximation
\beno
\left|\mbP((I(\truth) \, \norm{\bY}_1)^{1/2} \, (\mle - \truth) - \Delta \in \mA) - \Phi(\bZ \in \mA)\right|
\ee
is bounded above by
\beno
\dfrac{83}{(\widetilde{\lambda}_{\min}^{\epsilon})^{3/2}} \,
\sqrt{\dfrac{p^{7/2}}{\mbE \, \norm{\bY}_1}}
+  \dfrac{4}{\mbE \, \norm{\bY}_1} + \dfrac{8 \, \left[D_{g}\right]^{+}}{\left(\mbE \, \norm{\bY}_1 \right)^2}
\ee
and $\Delta$ satisfies
\beno
\mbP\left(\norm{\Delta}_2 \leq 
\dfrac{\sqrt{2} \, C^2 \, p^{5/2}}{\sqrt{\mbE \norm{\bY}_1}} \, 
\dfrac{\widetilde{\lambda}_{\max}^{\star} }{(\widetilde{\lambda}_{\min}^{\epsilon})^{5/2}} \right)
&\geq& 1 - \exp\,(-2\,p) \, - \, \dfrac{5+ 8 \, C_0}{\mbE \, \norm{\bY}_1},
\ee
where $C>0$ is the constant given in Theorem \ref{thm1} and $C_0>0$ is the constant given in Assumption \ref{assump1}, both independent of $N$ and $p$.
\end{theorem}

\s

Theorem \ref{thm2} serves as a foundation for establishing the asymptotic normality of the maximum likelihood estimator $\mle$. 
If 
\beno
\lim\limits_{N \to \infty} \; 
\left[ \dfrac{83}{(\widetilde{\lambda}_{\min}^{\epsilon})^{3/2}} \,
\sqrt{\dfrac{p^{7/2}}{\mbE \, \norm{\bY}_1}}
+  \dfrac{4}{\mbE \, \norm{\bY}_1} + \dfrac{8 \, \left[D_{g}\right]^{+}}{\left(\mbE \, \norm{\bY}_1 \right)^2} \right] 
\= 0,  
\ee
Theorem \ref{thm2} implies 
$(I(\truth) \, \norm{\bY}_1)^{1/2} \, (\mle - \truth) - \Delta$ will converge
in distribution
to a standard multivariate normal random vector,
as the error bound on the multivariate normal approximation will vanish in this case.  
The term $\Delta$ can be viewed as an error term,
resulting from the fact that the normal approximation in Theorem \ref{thm2} is obtained via a multivariate Taylor approximation 
in order to bridge the distributional gap between key statistics which admit forms amenable to existing 
theorems for the normal approximation 
and the parameter vectors of interest,
thus introducing an additional source of error in the normal approximation. 

While involved, 
the above condition for asymptotic multivariate normality essentially places restrictions 
on the dependence induced through the single-layer basis network $\bY$ 
measured by $[D_{g}]^{+}$,
as well as the smallest eigenvalue of the dyad-based information matrix $\mcI(\nat)$ 
in a neighborhood of the data-generating parameter vector $\truth$ as measured by $\widetilde{\lambda}_{\min}^{\epsilon}$, 
and the model dimension $p$. 
As a result, 
if the information matrix $\mcI(\nat)$ is nearly singular at $\truth$, 
in which case $\widetilde{\lambda}_{\min}^{\epsilon}$ will be small,
the error of the normal approximation will be uniformly larger (all else equal).
Likewise, 
if the edge dependence in $\bY$ is large as measured by $[D_{g}]^{+}$,
we may not have sufficient activated dyads to ensure the error bound is small,
as $\norm{\bY}_1$ may not be tightly concentrated around $\mbE \, \norm{\bY}_1$. 
The dependence of the error approximation on the dimension of the random vector is a known challenge 
in establishing multivariate normality \citep[e.g.,][]{Raic19}.  
All quantities which are not explicit constants can increase or decrease with $N$,
with the rates of these increases or decreases having implications for the rate of convergence in distribution. 
Theorem \ref{thm2} demonstrates 
that the allowable scaling for most of quantities is with respect to the expected number of activated dyads 
$\mbE \, \norm{\bY}_1$.

We further examine Theorem \ref{thm2} through an example where  
$\bY$ is a Bernoulli random graph model,
which assumes edge variables are independent Bernoulli random variables with probability $\pi \in (0, 1)$. 
Under this model, 
$[D_{g}]^{+} = 0$ owing to the independence of edge variables
and $\mbE \norm{\bY}_1 = \pi \, \binom{N}{2}$. 
Under this scenario, 
we can show that 
\beno
\left|\mbP((I(\truth) \, \norm{\bY}_1)^{1/2} (\mle - \truth) - \Delta \in \mA) - \Phi(\bZ \in \mA)\right|
\ee
is bounded above by
\beno
\dfrac{166}{\sqrt{\pi \, (\widetilde{\lambda}_{\min}^{\epsilon})^{3}}} \; \dfrac{p^{1.75}}{N} + \dfrac{16}{\pi  N^2}, 
\ee
with the additional bound   
\beno
\mbP\left(\norm{\Delta}_2 \leq 
\dfrac{\sqrt{2} \, C^2 \, p^{2.5}}{\sqrt{\pi\,\binom{N}{2}}} \, 
\dfrac{\widetilde{\lambda}_{\max}^{\star} }{(\widetilde{\lambda}_{\min}^{\epsilon})^{2.5}} \right)
&\geq& 1 - \exp\,(-2\,p) \, - \, \dfrac{5+ 8 \, C_0}{\pi\,\binom{N}{2}},
\ee
where $C>0$ is the constant given in Theorem \ref{thm1} and $C_0>0$ is the constant given in Assumption \ref{assump1}, both independent of $N$ and $p$.
If $\widetilde{\lambda}_{\min}^{\epsilon}$ and $\pi$ are both bounded away from $0$,
then the error of the normal approximation will convergence to $0$ 
provided $(p^{2.5} \, \widetilde{\lambda}_{\max}^{\star}) \,/\, N \to 0$ as $N \to \infty$,
which is sufficient to ensure $\norm{\Delta}_2$ converges in probability to $0$. 
Under the fully saturated model specification for \eqref{general_model} ($H = K$), 
the Binomial theorem shows that $p = 2^{K} - 1 \leq 2^K$. 
Hence, 
the dimension restriction on $p$ in turn implies a restriction on the allowable rate of growth of the number of layers $K$ with $N$,
where a sufficient condition for $(p^{2.5} \, \widetilde{\lambda}_{\max}^{\star}) \,/\, N \to 0$ 
is for $K \leq .5 \, \log N$. 
In other words, 
the number of layers $K$ can grow at most logarithmically with $N$ in the fully saturated model. 
In cases when the number of interaction terms included in the cross-layer dependence probability model 
is fixed, 
$K$ may admit a sublinear scaling  with $N$.

%% file: global_test.tex
We outline a procedure for model selection that controls the false discovery rate,
leveraging the results of Theorems \ref{thm1} and \ref{thm2}. 
Hotelling's $T$-squared statistic can be used to conduct a global test for 
$H_0: \truth = \bmu$ versus $H_1:\truth \neq \bmu$, where $\bmu \in \mbR^p$ is the value of $\nat$ we want to test.
We will mostly be interested in the case when $\bmu = \bm{0}_p$,
i.e.,
the zero vector of dimension $p$. 
If the global test is rejected, 
or is not of interest,  
we can perform model selection by leveraging the multivariate normal approximation 
to obtain univariate normal approximation results for the components of $\mle$
and proceed to test each component: $ H_{i,0} : \theta^\star_i=\mu_i$ 
versus $H_{i,1}: \theta^\star_i\neq\mu_i$,
for $i = 1, \ldots p$ and $\mu_i \in \mbR$. 
In general, 
$\mu_i = 0$ will allow us to test whether the estimated effect $\widehat\theta_i$ is present in the model
(i.e., whether $\theta^\star_i \neq 0$).  
One challenge in this approach lies in the fact that the model selection procedure 
is sensitive to multiple testing error. 

To ensure a more reliable procedure for identifying cross-layer dependence effects in multilayer networks while mitigating the risk of spurious discoveries, we elaborate a model selection algorithm that employs suitable multiple testing adjustments to control the false discovery rate. 
We provide simulation examples of four different univariate testing procedures including Bonferroni, Benjamini-Hochberg, 
Hochberg, 
and Holm procedures in Section \ref{sec:sim_norm}. 
In simulation studies, 
all four univariate testing procedures exhibit strong statistical power 
for detecting non-zero parameters,
while controlling the false discovery rate at a preset family-wise significance level.

%% file: Lazega.tex
We present a case study using a dataset on corporate law partnership among a Northeastern US corporate law firm in New England collected by \citet{Lazega2001}. The dataset collected information about three types of cooperation among 71 lawyers in the corporate law firm, resulting in three networks including the strong-coworker network, the advice network, and the friendship network. Since the cooperation relationship collected are not symmetric, we only consider a connection to be present when both sides acknowledged their cooperation. We treat these three types of networks as a three-layer multilayer network embedded among the 71 lawyers. A summary of this multilayer network is provided in Table \ref{t_datasum}. We apply the separable multilayer network model in \eqref{general_model} with the specification in \eqref{eq:sim_model} to Lazega's lawyer network, i.e., $\nat = (\theta_{1},\,\theta_{2},\,\theta_{3},\,\theta_{1,2},\,\theta_{1,3}, \,\theta_{2,3})$. The maximum pseudolikelihood estimator $\mple$ is computed from the observed network,
the results of which are provided in Table \ref{t_Lazeg_mple}.
 \begin{figure}[t]
\centering 
  \includegraphics[width=4in, height = 1.5in]{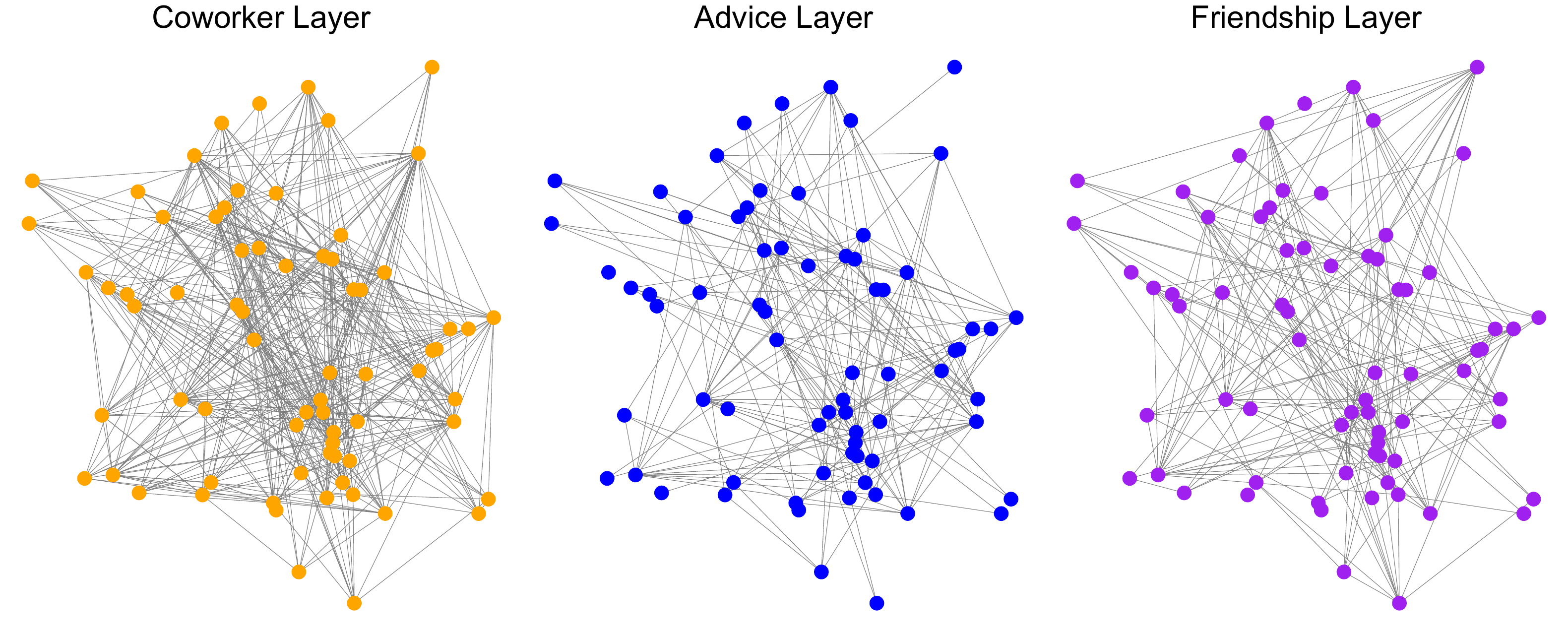}
  \caption{The coworker layer, the advice layer and the friendship layer of Lazega's corporate law partnership network.}
  \label{Figure4}
\end{figure}

\begin{table}[t]
\begin{center}
\caption{\label{t_datasum} Summary of Lazega's corporate law partnership data with 71 lawyers (nodes).  }
\begin{tabular}{ l  r  r }
\hline
 & Average Node Degree & Number of Edges \\
\hline
 Co-Worker Layer  & 11 & 378  \\

Advice Layer & 5 & 175  \\

Friendship Layer & 5 & 176  \\
\hline
\end{tabular}
\end{center}
\end{table}

\begin{table}[t]
\begin{center}
\caption{\label{t_Lazeg_mple}MPLEs (and standard errors) of the separable multilayer network model for the Lazega's lawyer network.}
\resizebox{\columnwidth}{!}{
\begin{tabular}{ c  c  c  c  c  c } 
\hline
$\widetilde\theta_{1}$ &$\widetilde\theta_{2}$ &$\widetilde\theta_{3}$ & $\widetilde\theta_{1,2}$ & $\widetilde\theta_{1,3}$ & $\widetilde\theta_{2,3}$ \\ 
\hline
$-1.450 \; (.263)$ & $-3.334 \; (.244)$ & $-2.695\; (.256)$ & $1.801 \; (.244)$ & $0.218 \; (.247)$ & $2.458 \; (.231)$ \\

Coworker (C) & Advice (A) & Friendship (F) & C $\times$ A & C $\times$ F & A $\times$ F \\
\hline
\end{tabular}
}
\end{center}
\end{table}

As shown in Table \ref{t_Lazeg_mple}, 
the maximum pseudolikelihood estimates 
$\widetilde\theta_{1}, \, \widetilde\theta_{2}$, and $\widetilde\theta_{3}$ correspond to the estimated single-layer effects of the coworker layer, the advice layer, and the friendship layer, respectively,
whereas $\theta_{1,2}, \, \theta_{1,3}$, and $\theta_{2,3}$ correspond to the layer interaction effects.  
We can calculate the conditional log-odds of each edge being present in the multilayer network given the rest of the network. For example,  if lawyer $i$ and lawyer $j$ are observed to be coworkers and are friends at the same time, the odds of these two lawyers to have an advice relationship is given by
\beno
\dfrac{\mbP(X_{i,j}^{(A)} = 1 \,|\, \bX_{i,j}^{(C)} = 1, \, \bX_{i,j}^{(F)} = 1)}
{\mbP(X_{i,j}^{(A)} = 0 \,|\, \bX_{i,j}^{(C)} = 1, \, \bX_{i,j}^{(F)} = 1)}
\quad = \quad \exp\left(
\widetilde\theta_{2} +  \widetilde\theta_{1,2} \, x_{i,j}^{(C)} + \widetilde\theta_{2,3} \, x_{i,j}^{(F)}\right) \s \\
 = \quad \exp\left(
-3.334 +  1.801 \, x_{i,j}^{(C)} + 2.458 \, x_{i,j}^{(F)}\right) 
 \quad = \quad  2.522,
\ee
providing interpretation of the interaction and influence among the different layers.

\begin{figure}[t]
 \centering
  \includegraphics[width=0.8\textwidth]{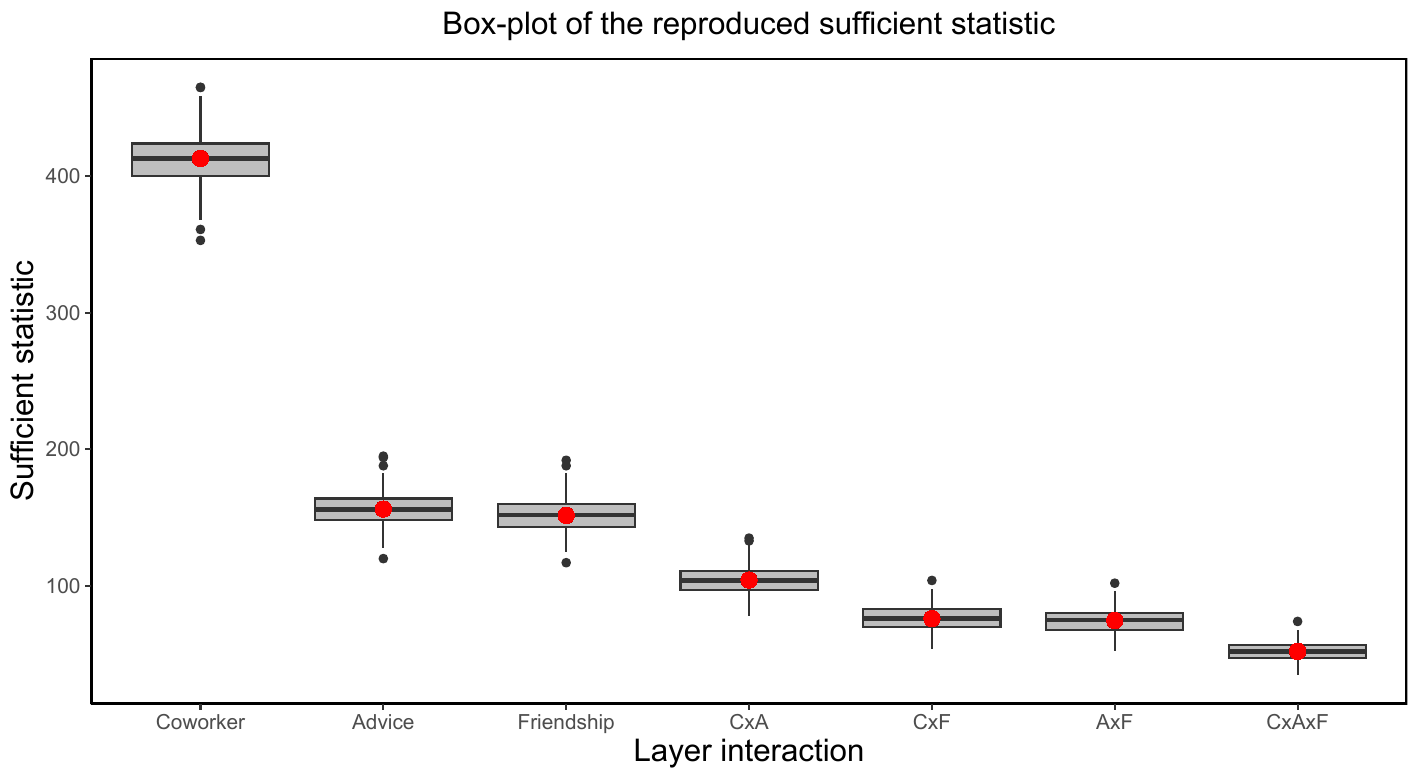} 
  \caption{Box-plot of reproduced statistics from 10 simulated samples using the MPLE obtained from the Lazega's lawyer network.
Red dots are values of the observed sufficient statistics of Lazega's lawyer network.} 
  \label{Figure_Lazega}
\end{figure}

Next, 
we use the MPLE to reproduce multilayer networks of the same size and compare the sufficient statistics 
of the simulated networks and the Lazega's lawyer network.
We recover the basis network according to Proposition \ref{prop:inference}, i.e.,
a dyad is activated if and only if at least one of its layers has a present edge in the Lazega's lawyer multilayer network.
We then populate layers of all activated dyads according to equation \eqref{eq:sim_model} by the MPLEs obtained in Table \ref{t_Lazeg_mple}.
Comparisons of the sufficient statistics between the observed Lazega's lawyer network and the simulated networks with 10 replications are provided in 
Figure \ref{Figure_Lazega}. 
Such comparisons serve two key purposes. 
First, 
such comparisons are an established method of diagnosing model fit in the statistical network analysis literature 
\citep{HuGoHa08}, 
and second, 
provide a check on the approximate solution to the score equation.
Note that MPLEs are not guaranteed to reproduce (on average) observed values of sufficient statistics in exponential families---in contrast to MLEs.
The relative $\ell_2$-error of the sufficient statistics between the observed and the average of the 10 simulated networks is $0.09$, suggesting a successful re-construction of the observed network statistics.

\hide{
The R package we developed for the simulation analysis in Section \ref{sec:sim} and the application analysis of Lazega's corporate law partnership data in Section \ref{sec:app} is available on GitHub: https://github.com/jiaheng-li/cross-layer-dependence-mlyrnetwork-simulation.
}

%% file: supplement.tex
\pagebreak

\appendix

\makeatletter

\setcounter{page}{1}

\setcounter{section}{0}

\setcounter{com}{0}

\begin{center}
\Large\bf\textsc\bf
Supplement:\\
Learning cross-layer dependence structure in multilayer networks \s\s

\normalsize
{\normalfont\textsc{By Jiaheng Li and Jonathan R. Stewart}} \s 
\\
{\normalfont\em Department of Statistics, Florida State University} \s\s
\end{center}

\s\s

\noindent
{Appendix \ref{sec:prop_proof}: Proof of Proposition \ref{prop:inference}}\dotfill\pageref{sec:prop_proof}\s\\ 
{Appendix \ref{sec:lem1_proof}: Proof of Lemma \ref{lem:min-eig}}\dotfill\pageref{sec:lem1_proof}\s\\
{Appendix \ref{sec:concentration}: Concentration inequalities for multilayer networks}\dotfill\pageref{sec:concentration}\s\\
{Appendix \ref{sec:pf_thm1}: Proof of Theorem \ref{thm1} and Corollary \ref{corollary}}\dotfill\pageref{sec:pf_thm1}\s\\
{Appendix \ref{sec:pf_minimax}: Proof of Theorem \ref{thm:minimax} and Corollary \ref{cor:minimax}}\dotfill\pageref{sec:pf_thm1}\s\\
{Appendix \ref{sec:pf_prop2}: Proposition \ref{prop:suff_norm} and proof}\dotfill\pageref{sec:pf_prop2}\s\\
{Appendix \ref{sec:pf_thm2}: Proof of Theorem \ref{thm2}}\dotfill\pageref{sec:pf_thm2}\s\\
{Appendix \ref{sec:add_sim}: Additional simulation results}\dotfill\pageref{sec:add_sim}\s\\

\section{Proof of Proposition \ref{prop:inference}} 
\label{sec:prop_proof} 

We prove Proposition \ref{prop:inference} from Section \ref{sec2}. 
\input{proof_prop_inference}
\s\s

\section{Proof of Lemma \ref{lem:min-eig}}
\label{sec:lem1_proof}

We prove Lemma \ref{lem:min-eig} from Section \ref{sec2}.
\input{proof_min_eig_lemma}
\s

\section{Concentration inequalities for multilayer networks} 
\label{sec:concentration}

We establish the concentration inequality of gradients of log-likelihood functions 
of multilayer networks in Lemma \ref{lem:concentration_likelihood}.
Recall the definition $[D_{g}]^{+} \coloneqq \max\{0, \, D_{g}\}$, 
where  
\beno 
D_{g}
&\coloneqq& \dsum_{\{i,j\} \prec \{v,w\} \subset \mN} \, \cov(Y_{i,j}, \, Y_{v,w}), 
\ee 
with $\{i,j\} \prec \{v,w\}$ implying the sum is taken with respect to the lexicographical ordering 
of pairs of nodes, 
and where $g : \mbY \mapsto (0, 1)$ is the marginal probability mass function of $\bY$.  

\s\s

\input{concentration_likelihood}



\s\s
\subsection{Auxiliary results}
\input{suff_exp}

\section{Proof of Theorem \ref{thm1}}
\label{sec:pf_thm1}
We prove Theorem \ref{thm1} from Section \ref{sec3}.
\input{proof_thm1}

\s\s
\section{Proposition \ref{prop:suff_norm} and proof} 
\label{sec:pf_prop2}

In order to establish a bound on the error of the multivariate normal approximation 
for estimators of data-generating parameters, 
we first establish an error bound on the multivariate normal approximation 
of a standardization of the 
sufficient statistic vector $\bs(\bX)$ of the exponential family distribution of $\bX$ given $\bY$,
derived in Lemma \ref{lem:s_hetero},
in Proposition \ref{prop:suff_norm} 
using a Lyapunov type bound presented in \citet{Raic19}.
Proposition \ref{prop:suff_norm} provides the basis for our normality proof for estimators
which we present in Theorem \ref{thm2}.
\input{prop2_normal}

\input{proof_prop2}
\s\s
\section{Proof of Theorem \ref{thm2}}
\label{sec:pf_thm2}

\input{proof_thm2}

\section{Additional simulation results} 
\label{sec:add_sim} 

Additional simulation results that enhance those contained in Section \ref{sec:sim} are provided in this section.

\subsection{Normal approximation with different basis networks}
\label{subsec:norm_sim}

The multivariate normality of $\mple$ is tested by Zhou-Shao's multivariate normal test \citep{Zhou13}, and the p-values are provided in tabel \ref{ZS-test}. Q-Q plots of $\mple$ estimated from 6 different model-generating parameters with a dense Bernoulli basis network, a sparse Bernoulli basis network, a stochastic block model (SBM) generated basis network, and a latent space model (LSM) generated basis network are shown in Figure \ref{qqplot_dense}, \ref{qqplot_sparse}, \ref{qqplot_SBM} and \ref{qqplot_LSM}, respectively.

\begin{table}[t]
\begin{center}
\caption{\label{ZS-test} P-values of the Zhou-Shao's test for multivariate normality of $\mple$ for 6 model-generating parameters ($\truth_1$, $\truth_2$, $\truth_3$, $\truth_4$, $\truth_5$, $\truth_6$) estimated from 250 network samples at size 1000 on four basis network structures. All p-values are larger than .1. \s} 
\begin{tabular}{ l  r r r  r r  r  } 
\hline
 Basis network model  & $\truth_1$ & $\truth_2$  & $\truth_3$ & $\truth_4$ & $\truth_5$ & $\truth_6$ \\ 
\hline
  Dense Bernoulli & .138 & .473 & .053 & .699 & .587 & .983  \\

 Sparse Bernoulli & .554 & .132 & .232 & .634 & .904 & .373  \\

 SBM & .650 & .891 & .982 & .975 & .871 & .674 \\

 LSM  & .859 & .831 & .500 & .227 & .613 & .409  \\
\hline
\end{tabular}
\end{center}
\end{table}

\vspace{.25in}

\begin{center}
\includegraphics[width=0.9\textwidth]{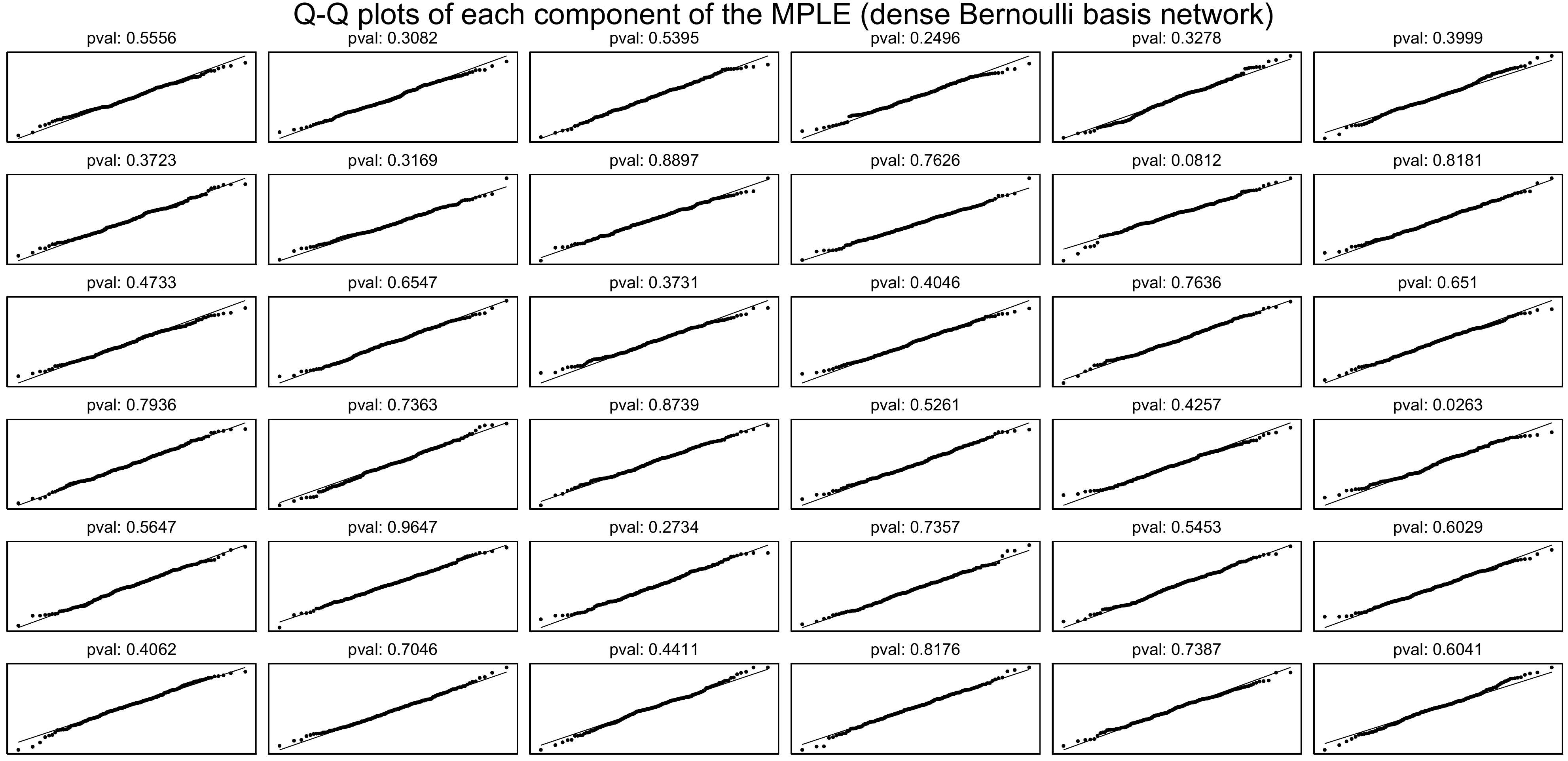} 
\captionof{figure}{Q-Q plots and p-values of six components of $\mple$ estimated from 250 multilayer network samples at size 1000 on the dense Bernoulli basis network for 6 model-generating parameters on each row.}\label{qqplot_dense}
\end{center}

\begin{center}
\includegraphics[width=0.9\textwidth]{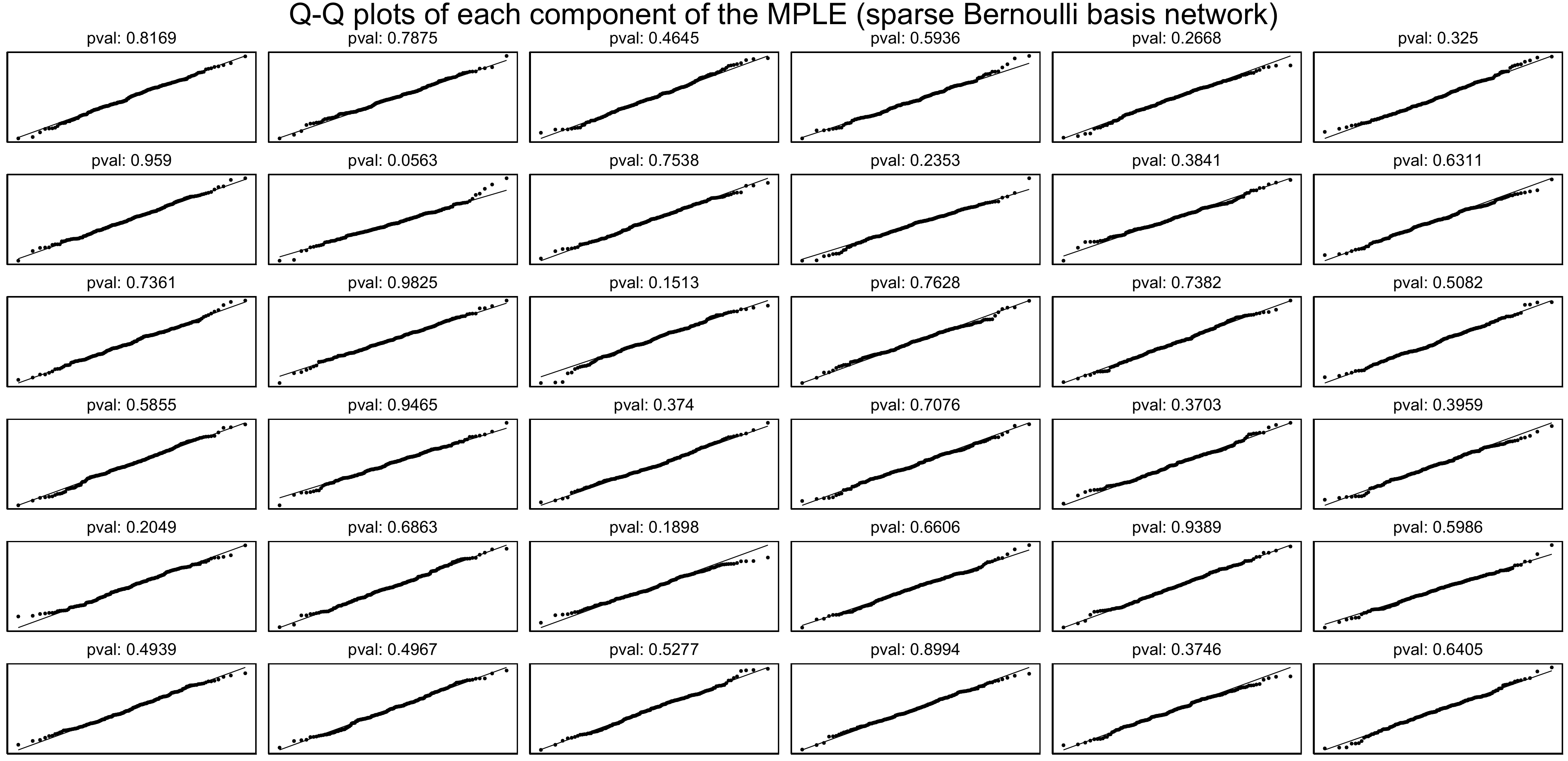} 
\captionof{figure}{Q-Q plots and p-values of six components of $\mple$ estimated from 250 multilayer network samples at size 1000 on the sparse Bernoulli basis network for 6 model-generating parameters on each row.}\label{qqplot_sparse}
\end{center}

\begin{center}
\includegraphics[width=0.9\textwidth]{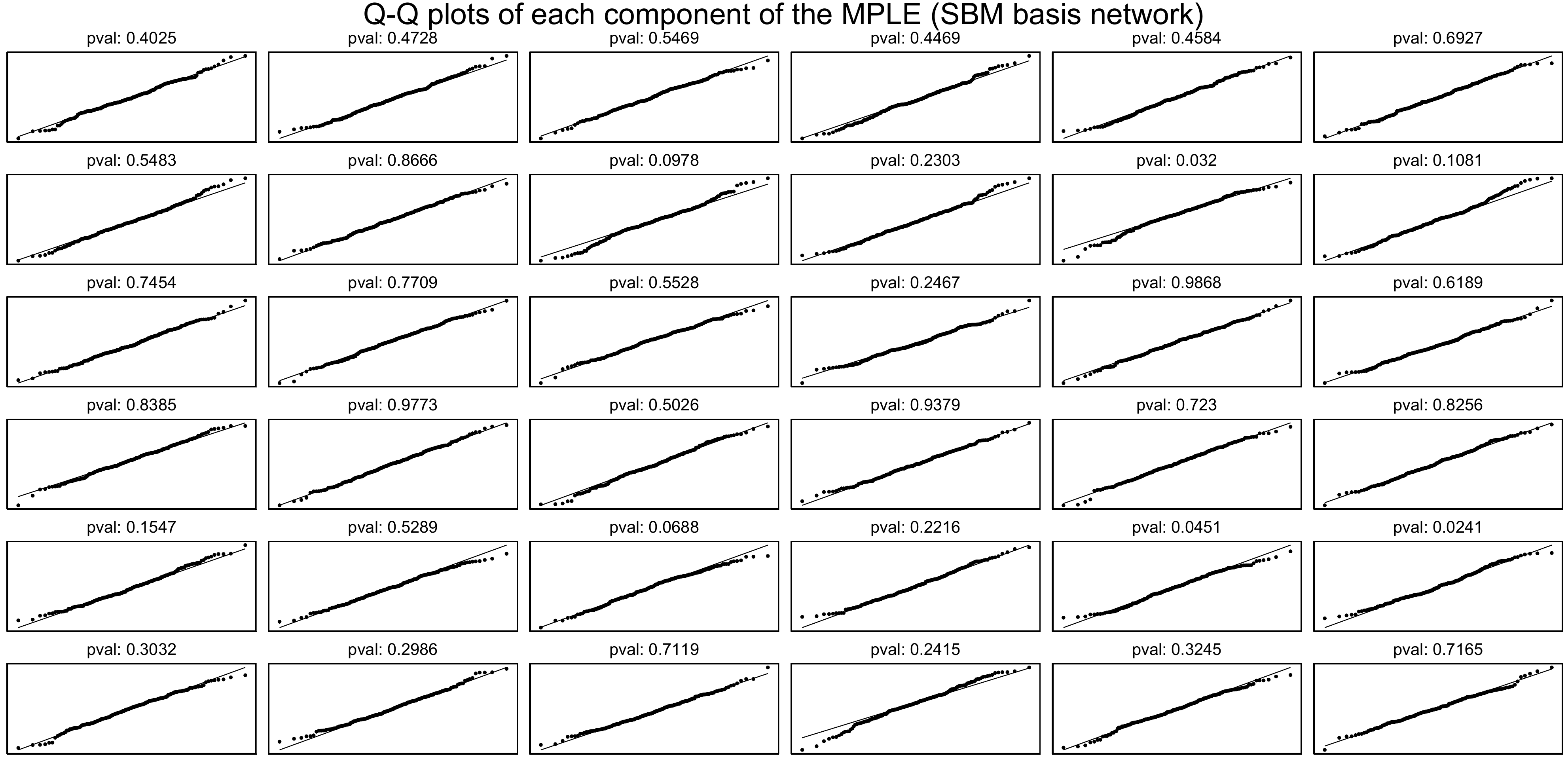}
\captionof{figure}{Q-Q plots and p-values of six components of $\mple$ estimated from 250 multilayer network samples at size 1000 on the SBM generated basis network for 6 model-generating parameters on each row.}\label{qqplot_SBM}
\end{center}

\begin{center}
\includegraphics[width=0.9\textwidth]{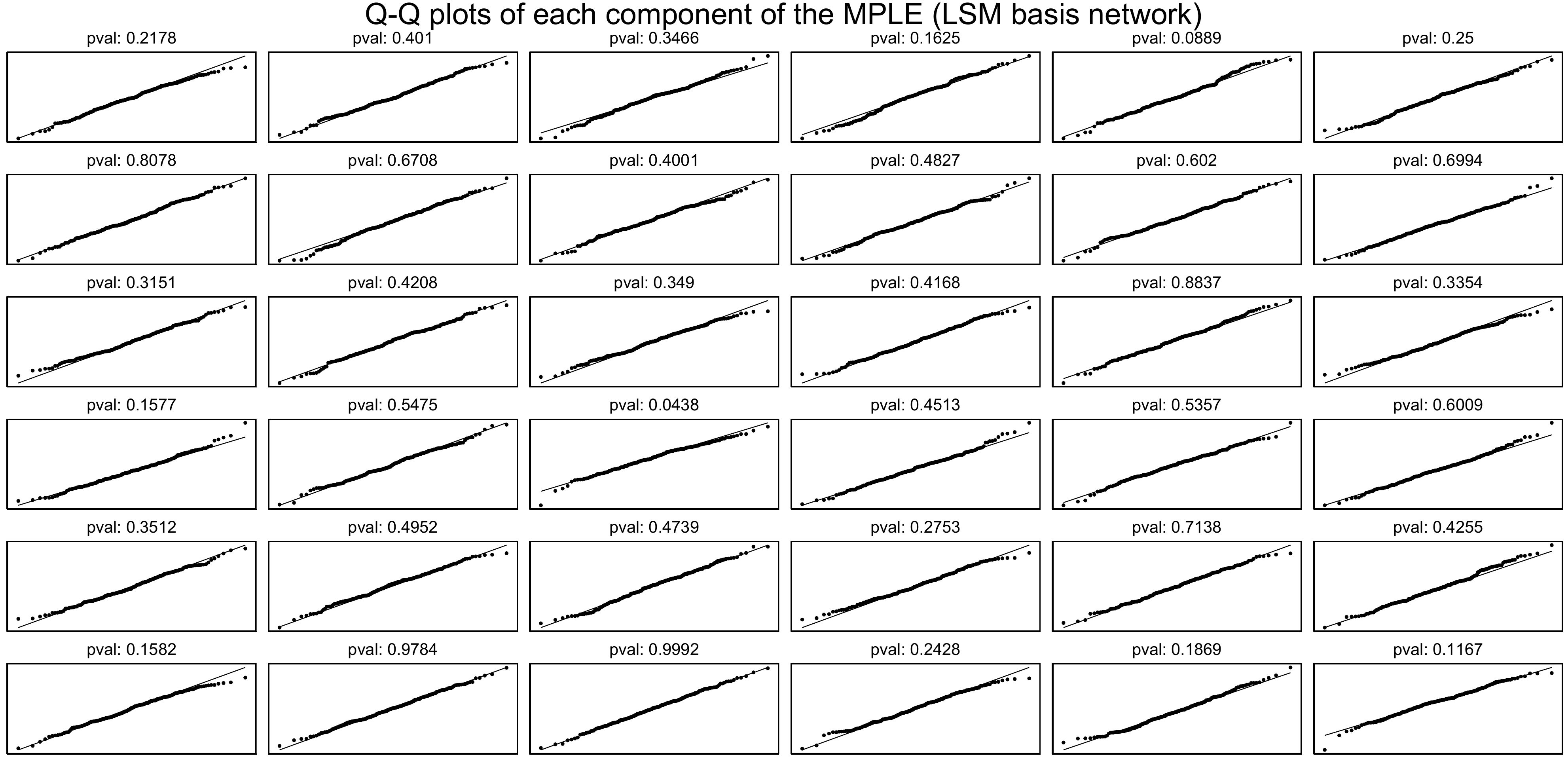}
\captionof{figure}{Q-Q plots and p-values of six components of $\mple$ estimated from 250 multilayer network samples at size 1000 on the LSM generated basis network for 6 model-generating parameters on each row.}\label{qqplot_LSM}
\end{center}

\subsection{Additional results on the false discovery rate}
\label{more fdr}
The false discovery rate (FDR) of the multiple testing correction procedures of
Bonferroni, Benjamini-Hochberg, Hochberg, and Holm to detect non-zero components of $\truth$ at a family-wise significance level of $\alpha = 0.05$ with a sparse Bernoulli basis network, an SBM generated basis network and an LSM generated basis network are provided in Table \ref{fdr_sparse}, \ref{fdr_SBM} and \ref{fdr_LSM}, respectively (recall that components $\theta_{1,3}^\star$ and $\theta_{3}^\star$ of $\truth$ are set to 0). The receiver operating characteristic (ROC) curves for $\mple$ of 6 selected model-generating parameters on four basis network structures are provided in each of the subplot of Figure \ref{ROC}.

\begin{table}[t]
\begin{center}
\caption{\label{fdr_sparse} False discovery rates of four procedures for detecting non-zero effects of six model-generating parameters ($\truth_1$, $\truth_2$, $\truth_3$, $\truth_4$, $\truth_5$, $\truth_6$) estimated from 250 multilayer network samples at size 1000 on the sparse Bernoulli basis network. All FDRs are smaller than 0.05.} 
\begin{tabular}{l  r r r  r r  r } 
\hline
  Procedure & $\truth_1$ & $\truth_2$  & $\truth_3$ & $\truth_4$ & $\truth_5$ & $\truth_6$ \\ 
\hline
  Bonferroni    & .002  & .003 & .003 & .003 & .003 & .011  \\

  Benjamini-Hochberg   & .020 & .011 & .022 & .022 & .014 & .017  \\

 Hochberg's &    .009 & .008 & .012 & .010 & .010 & .014 \\

 Holm's  &   .007 & .008 & .011 & .009 & .006 &  .014 \\
\hline
\end{tabular}
\end{center}
\end{table}

\begin{table}
\begin{center}
\caption{\label{fdr_SBM} False discovery rates of four procedures for detecting non-zero effects of six model-generating parameters ($\truth_1$, $\truth_2$, $\truth_3$, $\truth_4$, $\truth_5$, $\truth_6$) estimated from 250 multilayer network samples at size 1000 on the SBM generated basis network. All FDRs are smaller than 0.05.} 
\begin{tabular}{l  r r r  r r  r } 
\hline
  Procedure & $\truth_1$ & $\truth_2$  & $\truth_3$ & $\truth_4$ & $\truth_5$ & $\truth_6$ \\ 
\hline
  Bonferroni    & .002  & .002 & .003 & .001 & .001 & .004  \\

  Benjamini-Hochberg   & .022 & .013 & .014 & .015 & .015 & .018  \\

 Hochberg's &    .009 & .014 & .01 & .008 & .011 & .014 \\

 Holm's  &   .009 & .013 & .005 & .009 & .009 &  .011 \\
\hline
\end{tabular}
\end{center}
\end{table}

\begin{table}
\begin{center}
\caption{\label{fdr_LSM} False discovery rates of four procedures for detecting non-zero effects of six model-generating parameters ($\truth_1$, $\truth_2$, $\truth_3$, $\truth_4$, $\truth_5$, $\truth_6$) estimated from 250 multilayer network samples at size 1000 on the LSM generated basis network. All FDRs are smaller than 0.05.} 
\begin{tabular}{l  r r r  r r  r } 
\hline
  Procedure & $\truth_1$ & $\truth_2$  & $\truth_3$ & $\truth_4$ & $\truth_5$ & $\truth_6$ \\ 
\hline
  Bonferroni    & .004  & .006 & .000 & .005 & .003 & .004  \\

  Benjamini-Hochberg   & .016 & .013 & .011 & .015 & .016 & .017  \\

 Hochberg's &    .009 & .014 & .009 & .011 & .010 & .011 \\

 Holm's  &   .008 & .014 & .009 & .011 & .007 &  .010 \\
\hline
\end{tabular}
\end{center}
\end{table}

\begin{figure}
\centering
\includegraphics[width=0.9\textwidth]{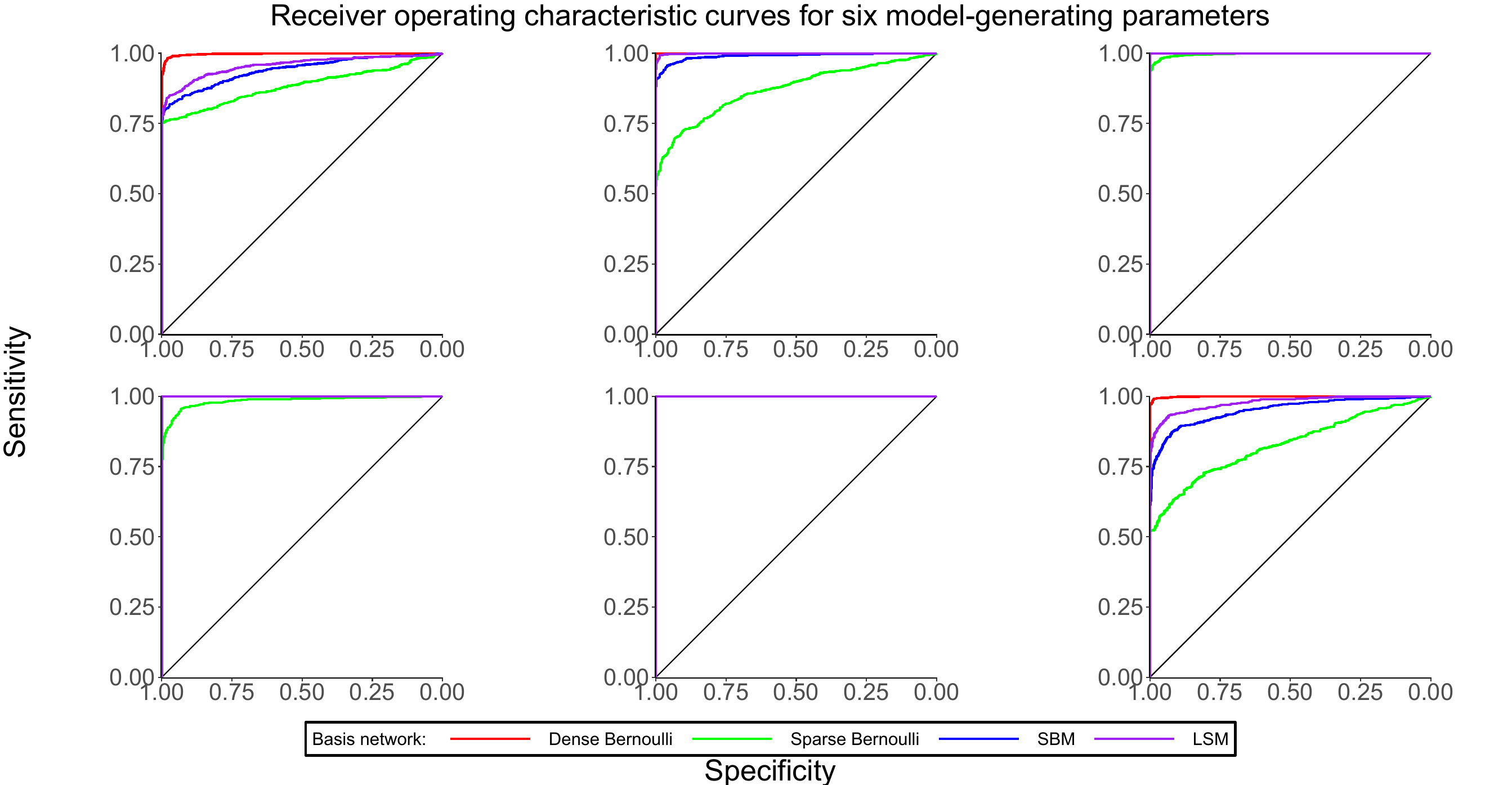} 
\label{ROC}
\captionof{figure}{ROC curves for $\mple$ estimated from 250 multilayer network samples at size $1000$ of six model-generating parameters on four different basis networks.}\label{ROC}
\end{figure}

\pagebreak

\input{new-submission.bbl}

%% file: proof_prop_inference.tex
For the first and second results,
define the set
\beno
\mA_{+}
&\coloneqq& \left\{ (\bx, \by) \in \mbX \times \mbY \;:\; h(\bx, \by) = 1 \right\},
\ee
and the vector-valued map $\bm{\varphi} : \mbX \mapsto \mbY$ by defining its components to be
\beno
\varphi_{i,j}(\bx)
\= \one(\norm{\bx_{i,j}}_1 > 0),
&& \{i,j\} \subset \mN, 
\ee
populating the vector $\bm{\varphi}(\bx)$ in the lexicographical ordering of the dyad indices $\{i,j\} \subset \mN$.  
By the definition of $h : \mbX \times \mbY \mapsto \{0, 1\}$ and $\bm{\varphi} : \mbX \mapsto \mbY$,
$\bm{\varphi}(\bx)  = \by$ for each pair $(\bx, \by) \in \mA_{+}$.  
Furthermore, 
the element $\by$ is unique for a given $\bx \in \mbX$,
because if there would exists some $\by^\prime \in \mbY$ 
such that $\by \neq \by^\prime$ 
with the property that 
$\{(\bx, \by), (\bx, \by^\prime)\} \subseteq \mA_{+}$, 
then there would exist a  pair $\{i,j\} \subset \mN$
such that $y_{i,j} = 1 - y^\prime_{i,j}$ which implies
$\one(\norm{\bx_{i,j}}_1 > 0) \neq y^\prime_{i,j}$.
In this case, $h(\bx, \by^\prime) = 0$ is
contradicting the assumption that $\{(\bx, \by)\} \in \mA_{+}$. 
By \eqref{general_model},
the functions $f$ and $g$ are assumed to be strictly positive in their respective domains.
Hence, 
$(\mbX \times \mbY) \setminus \mA_{+}$ is the largest null set of $\mbX \times \mbY$,
i.e.,
$\sepmodel(\mA) = 0$ if and only if $\mA \subseteq (\mbX \times \mbY) \setminus \mA_{+}$.
Thus, 
the first and second results are established.

For the third result,
note that $g$ is assumed to be strictly positive on its domain $\mbY$. 
Hence, 
$g(\by) = \mbP_{\nat}(\bY = \by) > 0$ for all $\by \in \mbY$
and  
$\sepmodel(\bX = \bx \,|\, \bY = \by)$
is therefore well-defined.
By definition,
\beno
\sepmodel(\bX = \bx \,|\, \bY = \by)
\= \dfrac{\sepmodel(\bX = \bx, \, \bY = \by)}{\sepmodel(\bY = \by)},
\ee
where $\sepmodel(\bY = \by)$ is the marginal probability of event $\bY = \by$ and is assumed to be equal to $g(\by)$.
The model form for $\sepmodel$ given in \eqref{general_model} implies 
\beno
\dfrac{\sepmodel(\bX = \bx, \, \bY = \by)}{\sepmodel(\bY = \by)}
\= \dfrac{f(\bx, \nat) \, g(\by) \, \psi(\nat, \by)}{g(\by)} \s \\
\= \exp(\log f(\bx, \nat) + \log \psi(\nat, \by)),
\ee
under the assumption that $h(\bx, \by) = 1$.
Hence,
\beno
\sepmodel(\bX = \bx, \bY = \by)
\= \sepmodel(\bX = \bX \;|\; \bY = \by) \; \sepmodel(\bY = \by)
\ee
so that
\beno
\log \, \sepmodel(\bX = \bx, \bY = \by)
\= \log \, \sepmodel(\bX = \bX \;|\; \bY = \by) + \log \, g(\by),
\ee
as $g(\by)$ is the marginal probability mass function of $\bY$,
i.e.,
$\sepmodel(\bY = \by) = g(\by)$. 
Lemma \ref{lem:s_hetero} establishes that $\sepmodel(\bX = \bX \;|\; \bY = \by)$ belongs to a minimal exponential family,
completing the proof of the third and last result of the proposition.  

\qed 

%% file: proof_min_eig_lemma.tex
Using \eqref{general_model},  
\beno
-\mbE \, \nabla_{\nat}^2 \, \ell(\nat; \bX, \bY)
\= \dsum_{\by \in \mbY} \, \dsum_{\bx \in \mbX} \, - \nabla_{\nat}^2 \, \ell(\nat; \bx, \by) \, 
\sepmodel(\bX = \bx \,|\, \bY = \by) \, g(\by) \s \\ 
\= \dsum_{\by \in \mbY} \, g(\by)  \dsum_{\bx \in \mbX} \, - \nabla_{\nat}^2 \, \ell(\nat; \bx, \by) \, 
\sepmodel(\bX = \bx \,|\, \bY = \by) \s \\ 
\= \dsum_{\by \in \mbY} \, g(\by) \, \dsum_{\{i,j\} \subset \mN \,:\, y_{i,j} = 1} \, \mcI(\nat) \s \\ 
\= \mcI(\nat) \, \dsum_{\by \in \mbY} \, g(\by)  \, \norm{\by}_1 \s\\
\= \mcI(\nat) \, \mbE \norm{\bY}_1.  
\ee
The above follows by 
exploiting the conditional independence of vectors $\bx_{i,j}$ ($\{i,j\} \subset \mN$) 
given $\bY = \by$ under \eqref{general_model}, 
which implies 
\beno
\ell(\nat; \bx, \by)
\= \dsum_{\{i,j\} \subset \mN} \log \sepmodel(\bX_{i,j} = \bx_{i,j} \,|\, \bY = \by), 
\ee
and from the fact that the conditional probability distribution of $\bX_{i,j}$ 
given $\bY$ 
is a degenerate point mass at $\bm{0}$ 
when $Y_{i,j} = 0$ so that $- \nabla_{\nat}^2 \, \ell(\nat; \bx, \by)$ is a sum of $\norm{\by}_1$ 
matrices each equal to $\mcI(\nat)$,
i.e.,
given $\by \in \mbY$,
we have 
\beno
&& \dsum_{\bx \in \mbX} - \nabla_{\nat}^2 \, \ell(\nat; \bx, \by) \, \sepmodel(\bX = \bx \,|\, \bY = \by) \s\\ 
\= \dsum_{\{i,j\} \subset \mN} \, \mbE\left[- \nabla_{\nat}^2 \, L_{i,j}(\nat, \bX_{i,j}, \bY) \,|\, \bY = \by \right]  
\;\;= \dsum_{\{i,j\} \subset \mN \,:\, y_{i,j} = 1} \, \mcI(\nat). 
\ee
The fact that $\mcI(\nat)$ is constant for all pairs $\{i,j\} \subset \mN$ satisfying $Y_{i,j} = 1$ 
follows from the form of \eqref{general_model},
which assumes each vector $\bX_{i,j}$ ($\{i,j\} \subset \mN$)
is conditionally independent and identically distributed, 
conditional on $\bY$.  
Hence, 
\beno
\mbE \left[ -\nabla_{\nat}^2 \, \ell(\nat; \bX, \bY)  \right]
\= \mcI(\nat) \, \mbE \, \norm{\bY}_1,
\ee 
which in turn implies 
\beno
\lambda_{\min}(-\mbE \nabla_{\nat}^2 \,\ell(\nat; \bX, \bY)) = \lambda_{\min}(\mcI(\nat)) \, \mbE \, \norm{\bY}_1\s\\
\lambda_{\max}(-\mbE \nabla_{\nat}^2 \,\ell(\nat; \bX, \bY)) = \lambda_{\max}(\mcI(\nat)) \, \mbE \, \norm{\bY}_1.
\ee

\qed 

%% file: concentration_likelihood.tex
\begin{lemma}
\label{lem:concentration_likelihood}
Consider a multilayer network model following the form of equation \eqref{general_model} and is 
defined on a set of $N \geq 3$ nodes and $K \geq 1$ layers. 
Define 
$\gradL \coloneqq - \nabla_{\nat} \, \ell(\nat; \bx, \by)$,  
where $\ell(\nat; \bx, \by)$ is the log-likelihood function. 
Then,
for all $t > 0$ and $\nat \in \mbR^p$, the probability
\beno
\mbP\left(\norm{\GradL - \mbE \, \GradL}_2 \geq t \right) 
\ee
is bounded above by
\beno
 \exp\left( -\dfrac{t^2}{36 \, \widetilde{\lambda}_{\max}^{\star} \, (\mbE \, \norm{\bY}_1 + [D_{g}]^{+}) \,  + 2 \, \sqrt{p} \, t} + \log \, p \right) + \dfrac{1}{\mbE \norm{\bY}_1}.
\ee
\end{lemma}

\llproof \ref{lem:concentration_likelihood}. 
By Proposition \ref{prop:inference},
\beno 
\ell(\nat;\bx,\by)
\= \log \, \mbP_{\nat}(\bX = \bx \,|\, \bY = \by)
+ \log \, g(\by).  
\ee
Thus, 
\be
\label{eq:887}
-\gradL
\= \nabla_{\nat} \, \log \, \mbP_{\nat}(\bX = \bx \,|\, \bY = \by)
 + \nabla_{\nat} \, \log \, g(\by) \s \\
\= s(\bx) - \mbE_{\nat} \, s(\bX),  
\ee
as $g(\by) = \mbP_{\nat}(\bY = \by)$ is assumed to not be a function of $\nat$. 
The last equation in \eqref{eq:887} follows from 
Lemma \ref{lem:s_hetero},
which showed that 
$\mbP_{\nat}(\bX = \bx \,|\, \bY = \by)$ is a minimal exponential family
with the natural parameter vector $\nat \in \mbR^p$ and the sufficient statistic vector $\bs(\bx)$ defined in Lemma \ref{lem:s_hetero},
inserting the familiar form of the score equation of an exponential family with respect to the natural parameter vector 
\citep[e.g., Proposition 3.10, p. 32,][]{Su19}.
Thus, 
\beno
- (\GradL   - \mbE \, \GradL ) 
\= \bs(\bX) - \mbE_{\nat} \, \bs(\bX) 
- \mbE \left[ \bs(\bX) - \mbE_{\nat} \, \bs(\bX)\right] \s \\
\= \bs(\bX) - \mbE \, \bs(\bX). 
\ee
Let $t > 0$ and $\nat \in \mbR^p$ be arbitrary and fixed 
and define  
$\mD_{2}(\nat,t)$ to be the event that 
$\norm{\GradL - \mbE \, \GradL}_{2}  \geq t$, i.e., 
\beno
\mD_{2}(\nat,t) \= \left\{ \bx \in \mbX \,:\, \norm{\bs(\bx)-\mbE  \, \bs(\bX)}_{2} \geq t \right\}.
\ee
Let $\epsilon > 0$ and define  
$\mE(\epsilon)$ to be the event that $|\norm{\bY}_1 - \mbE \norm{\bY}_1| \le \epsilon$,
i.e.,
\beno
\mE(\epsilon)
\= \left\{ \by \in \mbY \,:\, \left|\norm{\by}_1 - \mbE \norm{\bY}_1\right| \leq \epsilon \right\}. 
\ee 
We assume that $\epsilon > 0$ is chosen so that $\mE(\epsilon)$ is not empty, 
which implies $\mbP(\mE(\epsilon)) > 0$ 
as $g(\by)$ is assumed to be strictly positive on $\mbY$.  
By the law of total probability, 
\be
\label{divide and conquer}
\mbP\left(\mD_{2}(\nat,t) \right) 
&=& \mbP\left( \mD_{2}(\nat, t) \,|\, \mE(\epsilon) \right)\,
\mbP\left(\mE(\epsilon) \right) + 
\mbP\left( \mD_{2}(\nat,t) \,|\, \mE(\epsilon)^c \right) \, \mbP\left( \mE(\epsilon)^c \right) \s \\
&\leq& \mbP\left( \mD_{2}(\nat,t) \,|\, \mE(\epsilon) \right) + \mbP\left( \mE(\epsilon)^c \right).
\ee
Note that we have not necessarily guaranteed that $\mbP\left( \mE(\epsilon)^c \right) > 0$.  
However, 
if $\mbP\left( \mE(\epsilon)^c \right) = 0$ the non-conditional form of the law of total probability would yield the bound   
\beno
\mbP\left(\mD_{2}(\nat,t) \right)
&\leq&  \mbP\left( \mD_{2}(\nat,t) \,|\, \mE(\epsilon) \right),
\ee
which is strictly sharper than the bound we give in \eqref{divide and conquer}. 
We will use a divide-and-conquer strategy to bound each probability in \eqref{divide and conquer} in turn. 

To bound the first term in \eqref{divide and conquer}, let $\mU \coloneqq \{\bu \in \mbR^p \, : \, \norm{\bu}_2 \leq 1\}$ be a closed unit ball in $\mbR^p$. Define an $\epsilon$-net $\mV_\epsilon$ of $\mU \subset \mbR^p$. By Corollary 4.2.13 of \citep{Vershynin18}, there exists an $\epsilon$-net $\mV_\epsilon \subset \mU$ such that its cardinality satisfies $\log\, |\mV_\epsilon| \leq p \, \log \, (2\,\epsilon^{-1} + 1)$. Taking $\epsilon = 1/2$, for each $\bu \in \mU$, there exists a $\bv \in \mV_{1/2}$ such that $\norm{\bu - \bv}_2 \leq 1/2$, and by the Cauchy-Schwarz inequality, 
\be
\label{inner product of score eq}
\langle \, \bu \, , \, \nabla_{\nat} \, \ell(\nat;\bx,\by) \, \rangle \= \langle \, \bv \, ,\, \nabla_{\nat} \, \ell(\nat;\bx,\by) \,\rangle \, + \, \langle \, \bu - \bv \, ,\, \nabla_{\nat} \, \ell(\nat;\bx,\by)\,\rangle \s \\
& \leq & \langle \, \bu \, , \, \nabla_{\nat} \, \ell(\nat;\bx,\by) \, \rangle \, + \, \norm{\bu - \bv}_2 \, \norm{\nabla_{\nat} \, \ell(\nat;\bx,\by)}_2 \s \\
& \leq & \langle \, \bu \, , \, \nabla_{\nat} \, \ell(\nat;\bx,\by) \, \rangle \, + \, \dfrac{1}{2} \, \norm{\nabla_{\nat} \, \ell(\nat;\bx,\by) }_2.
\ee
If $\norm{\nabla_{\nat} \, \ell(\nat;\bx,\by) }_2 \neq 0$, we can choose 
\beno
u_i \= \dfrac{\nabla_{\nat} \, \ell(\nat;\bx,\by)_i}{\norm{\nabla_{\nat} \, \ell(\nat;\bx,\by) }_2},
\ee
so that $\norm{\bu}_2 \leq 1$ and $\bu \in \mU$. Next, re-write 
\beno
\langle \, \bu \, , \, \nabla_{\nat} \, \ell(\nat;\bx,\by) \, \rangle \= \dfrac{1}{\norm{\langle \, \bu \, , \, \nabla_{\nat} \, \ell(\nat;\bx,\by) \, \rangle }_2} \, \dsum_{i=1}^p \, (\nabla_{\nat} \, \ell(\nat;\bx,\by)_i)^2 \s \\
\= \norm{\nabla_{\nat} \, \ell(\nat;\bx,\by)}_2,
\ee
and together with \eqref{inner product of score eq}, we have
\be
\label{score ineq}
\norm{\nabla_{\nat} \, \ell(\nat;\bx,\by)}_2 \, \leq \, 2 \, \max\limits_{\bv \in \mV_{1/2}} \, \langle \, \bv \, , \, \nabla_{\nat} \, \ell(\nat;\bx,\by) \, \rangle. 
\ee
If $\norm{\nabla_{\nat} \, \ell(\nat;\bx,\by) }_2 = 0$, the inequality \eqref{score ineq} holds trivially. As a result of \eqref{score ineq}, for any $t > 0$,
\beno
\mbP\left(\mD_{2}(\nat,t) \,|\, \mE(\epsilon) \right) & \leq & \mbP \left(  2 \, \max\limits_{\bv \in \mV_{1/2}}  \, \langle \, \bv \, , \, \nabla_{\nat} \, \ell(\nat;\bx,\by) \, \rangle \, \ge \, t \right) \s \\
& \leq & \dsum_{\bv \in \mV_{1/2}} \, \mbP \left( \,\langle \, \bv \, , \, \nabla_{\nat} \, \ell(\nat;\bx,\by) \, \rangle \, \ge \,\dfrac{t}{2} \,  \right) \s \\
& \leq & \exp \, \left( p\, \log \, 5\right) \max\limits_{\bv \in \mV_{1/2}} \, \mbP \left( \,\langle \, \bv \, , \, \nabla_{\nat} \, \ell(\nat;\bx,\by) \, \rangle \, \ge \,\dfrac{t}{2} \,  \right). 
\ee
The last inequality is true because $\log\, |\mV_{1/2}| \leq p \, \log \, 5$. Note that
\be
\label{sum_over_dim}
\langle \, \bv \, , \, \nabla_{\nat} \, \ell(\nat;\bx,\by) \, \rangle  \= \dsum_{l=1}^p \, v_l \, [\nabla_{\nat} \, \ell(\nat;\bx,\by)]_l \s \\
\= \dsum_{l=1}^p \, v_l \, [s_l(\bx) - \mbE \, s_l(\bX)].
\ee
The form of \eqref{general_model} implies, 
through factorization principles,
that the dyad-based vectors $\bX_{i,j}$ ($\{i,j\} \subset \mN$) are conditionally independent given $\bY$ 
\citep[e.g.,][p. 11--13]{graphical_model_handbook}. 
Hence,
using Lemma \ref{lem:s_hetero},  
the components of the sufficient statistic vector decompose into the sum  
\beno
s_l(\bX) 
\= \dsum_{\{i,j\} \subset \mN} \, s_{l,i,j}(\bX_{i,j}),
&& l \in \{1, \ldots, p\},
\ee
so that the components of $\bs(\bX)$ are sums of bounded 
conditionally independent random variables given $\bY$.
As a result, equation \eqref{sum_over_dim} can be further decomposed into sums of independent random variables:
\beno
\langle \, \bv \, , \, \nabla_{\nat} \, \ell(\nat;\bx,\by) \, \rangle  \= \dsum_{\{i,j\} \subset \mN} \, \dsum_{l=1}^p \, v_l \, [s_{l,i,j}(\bx_{i,j}) \, - \, \mbE \, s_{l,i,j}(\bX_{i,j})].
\ee
Using the forms for $ s_{l,i,j}(\bX_{i,j})$ given in Lemma \ref{lem:s_hetero},
we have $0 \le s_{l,i,j}(\bX_{i,j}) \le Y_{i,j}$ $\mbP$-almost surely,
because $s_{l,i,j}(\bX_{i,j}) \in \{0, 1\}$ and $s_{l,i,j}(\bX_{i,j}) = 0$ if $Y_{i,j} = 0$ $\mbP$-almost surely. Then for each $\{i,j\} \subset \mN$, we have
\beno
\mbE \, \dsum_{l=1}^p \, v_l \, [s_{l,i,j}(\bx_{i,j}) \, - \, \mbE \, s_{l,i,j}(\bX_{i,j})] \= 0,
\ee
and by the Cauchy-Schwarz inequality, we obtain
\beno
\left|\,\dsum_{l=1}^p \, v_l \, [s_{l,i,j}(\bx_{i,j}) \, - \, \mbE \, s_{l,i,j}(\bX_{i,j})] \, \right| &\leq & \scalebox{0.97}{$\norm{\bv}_2 \, \sqrt{p} \, \norm{\bs_{l,i,j}(\bx_{i,j}) - \mbE\,\bs_{l,i,j}(\bX_{i,j})}_{\infty}$} \s \\
& \leq & \dfrac{3}{2} \, \sqrt{p}.
\ee
The last inequality follows from
\beno
\norm{\bv}_2 &\leq & \norm{\bu}_2 \, + \, \norm{\bu - \bv}_2 &\leq &1 \, + \, \dfrac{1}{2} &\leq & \dfrac{3}{2}.
\ee
The inequality is true because the construction of the $\epsilon$-net $\mV_{1/2} \subset \mU$ with $\epsilon=1/2$ ensures that such a $\bu \in \mU$ exists. We next bound the variance by 
\beno
\var \, \dsum_{\{i,j\}\subset \mN} \, \dsum_{l=1}^p \, v_l \, [s_{l,i,j}(\bx_{i,j}) \, - \, \mbE \, s_{l,i,j}(\bX_{i,j})] \s \\
= \dsum_{\{i,j\}\subset \mN} \, \dsum_{m=1}^p \, \dsum_{n=1}^p \, \cov \, (v_m \, s_{m,i,j} (\bX_{i,j}) \, ,\, v_n \, s_{n,i,j(\bX_{i,j})}) \s \\ 
= \dsum_{\{i,j\}\subset \mN}\, \dsum_{m=1}^p \, \dsum_{n=1}^p \, v_m \, v_n \, \cov(s_{m,i,j} (\bX_{i,j})\, s_{n,i,j} (\bX_{i,j})) \s \\
= \langle \, \bv \, , \, \norm{\by}_1 \, \mcI(\truth) \, \bv \,  \rangle \s \\
\leq  \norm{\bv}_2^2 \, \norm{\by}_1 \, \widetilde{\lambda}_{\max}^{\star} \s \\
 \leq  \dfrac{9}{4} \, \norm{\by}_1 \, \widetilde{\lambda}_{\max}^{\star},
\ee
where $\widetilde{\lambda}_{\max}^{\star}$ is the largest eigenvalue of the Fisher information of individual activated dyad defined in Lemma \ref{lem:min-eig} evaluated at the data-generating parameter $\truth$.
We then apply the one-sided Bernstein's inequality to obtain the upper bound for the conditional probability of $\mD_{2}(\nat,t)$ as follows: \cite[e.g.,][Theorem 2.8.4]{Vershynin18}
\be
\label{eq:hoef_lik}
\scalebox{0.9}{$\mbP\left( \mD_{2}(\nat,t) \given \bY = \by \right)$}
&\leq& \scalebox{0.9}{$\exp \, \left( p\, \log \, 5\right) \max\limits_{\bv \in \mV_{1/2}} \, \mbP \left( \,\langle \, \bv \, , \, \nabla_{\nat} \, \ell(\nat;\bx,\by) \, \rangle \, \ge \,\dfrac{t}{2} \,  \right)$} \s \\
& \leq & \scalebox{0.9}{$\exp\left(\dfrac{- \dfrac{(t/2)^2}{2}}{\dfrac{9}{4} \, \norm{\by}_1 \,\widetilde{\lambda}_{\max}^{\star} \, + \, \dfrac{1}{3}\, \dfrac{3}{2} \, \sqrt{p} \, \dfrac{t}{2} } \, + \, p \, \log \, 5 \right)$} \s \\
\= \scalebox{0.9}{$\exp\left(\dfrac{- t^2}{18 \, \norm{\by}_1 \,\widetilde{\lambda}_{\max}^{\star} \, + \, 2\, \sqrt{p} \, t } \, + \, p \, \log \, 5 \right)$}.
\ee
Using the law of total probability, 
we bound $\mbP\left( \mD_{2}(\nat, t) \given \mE(\epsilon) \right)$ as follows: 
\be
\label{divide_conquer}
\scalebox{0.9}{$\mbP\left( \mD_{2}(\nat, t) \,|\, \mE(\epsilon) \right)$}
&=& \scalebox{0.9}{$\dsum_{\by \in \mbY} \, \mbP\left( \mD_{2}(\nat,t) \cap [\bY = \by] \,|\, \mE(\epsilon) \right)$}   \s \\
\= \scalebox{0.9}{$\dsum_{\by \in \mE(\epsilon)} \, \mbP\left( \mD_{2}(\nat,t) \cap [\bY = \by] \,|\, \mE(\epsilon) \right)$}   \s \\
\= \scalebox{0.9}{$\dsum_{\by \in \mE(\epsilon)} \, \mbP\left( \mD_{2}(\nat,t) \,|\, [\bY = \by] \cap \mE(\epsilon)\right) \, 
\mbP(\bY = \by \,|\, \mE(\epsilon))$} \s \\ 
\= \scalebox{0.9}{$\dsum_{\by \in \mE(\epsilon)} \, \mbP(\mD_{2}(\nat,t) \,|\, \bY = \by) \, 
\dfrac{\mbP(\bY = \by)}{\mbP(\mE(\epsilon))}$},
\ee
noting that $[\bY = \by] \cap \mE(\epsilon) = [\bY = \by]$ whenever $\by \in \mE(\epsilon)$
and in the case when $\by \not\in \mE(\epsilon)$,
the intersection is empty,
implying 
\beno
\mbP(\bY = \by \,|\, \mE(\epsilon))
\= \dfrac{\mbP([\bY = \by] \cap \mE(\epsilon))}{\mbP(\mE(\epsilon))}
\= \begin{cases} 
\dfrac{\mbP(\bY = \by)}{\mbP(\mE(\epsilon))} & \by \in \mE(\epsilon) \\ 
0 & \by \not\in \mE(\epsilon)
\end{cases}. 
\ee
We now bound \eqref{divide_conquer} using the bound in \eqref{eq:hoef_lik}:  
\beno
& \dsum_{\by \in \mE(\epsilon)} \, \mbP(\mD_{2}(\nat,t) \,|\, \bY = \by) \,
\dfrac{\mbP(\bY = \by)}{\mbP(\mE(\epsilon))}  
\s \\ 
& \leq \,  \dsum_{\by \in \mE(\epsilon)} \, \exp\left(\dfrac{- t^2}{18 \, \norm{\by}_1 \,\widetilde{\lambda}_{\max}^{\star} \, + \, 2\, \sqrt{p} \, t } \, + \, p \, \log \, 5 \right) \, 
\dfrac{\mbP(\bY = \by)}{\mbP(\mE(\epsilon))} \s\\
&\leq  \, \exp\left(\dfrac{- t^2}{18 \, (\mbE\,\norm{\bY}_1 \, + \, \epsilon )\,\widetilde{\lambda}_{\max}^{\star} \, + \, 2\, \sqrt{p} \, t }  \, + \, p \, \log \, 5\right)  \, 
\dsum_{\by \in \mE(\epsilon)} \, \dfrac{\mbP(\bY = \by)}{\mbP(\mE(\epsilon))} \s\\ 
& = \exp\left(\dfrac{- t^2}{18 \, (\mbE\,\norm{\bY}_1 \, + \, \epsilon )\,\widetilde{\lambda}_{\max}^{\star} \, + \, 2\, \sqrt{p} \, t }  \, + \, p \, \log \, 5\right),
\ee
showing 
\beno
\mbP\left( \mD_{2}(\nat, t) \,|\, \mE(\epsilon) \right)
&\leq& \exp\left(\dfrac{- t^2}{18 \, (\mbE\,\norm{\bY}_1 \, + \, \epsilon )\,\widetilde{\lambda}_{\max}^{\star} \, + \, 2\, \sqrt{p} \, t }  \, + \, p \, \log \, 5\right). 
\ee
The replacement of $\norm{\by}_1$ by $\mbE\norm{\bY}_1 + \epsilon$ follows because 
$\norm{\by}_1 \le \mbE\norm{\bY}_1 + \epsilon$ for $\by \in \mE(\epsilon)$, 
resulting in the upper bound above.  
We bound the second term in the inequality \eqref{divide and conquer} 
using Chebyshev's inequality:
\beno
\mbP(\mE(\epsilon)^c) 
\= \mbP(\left| \norm{\bY}_1 - \mbE \, \norm{\bY}_1 \right| > \epsilon) \s \\
&\leq& \mbP(\left| \norm{\bY}_1 - \mbE \, \norm{\bY}_1 \right| \geq \epsilon) \s \\
&\leq&\dfrac{\var(\norm{\bY}_1)}{\epsilon^2}.
\ee
We bound the variance $\var(\norm{\bY}_1)$ as follows:
\beno
\var(\norm{\bY}_1)
\= \dsum_{\{i,j\} \subset \mN} \, \var \, Y_{i,j} 
+ 2 \, \dsum_{\{i,j\} \prec \{v,w\} \subset \mN} \, \cov(Y_{i,j}, \, Y_{v,w})\s\\
&\leq& \mbE \, \norm{\bY}_1 + 2 \, \dsum_{\{i,j\} \prec \{v,w\} \subset \mN} \, \cov(Y_{i,j}, \, Y_{v,w}),  
\ee
noting $Y_{i,j} \in \{0,1\}$ so that 
$\var \, Y_{i,j} = \mbP(Y_{i,j} = 1) \, \mbP(Y_{i,j} = 0) \leq \mbE \, Y_{i,j}$. 
Hence, 
\be
\label{ineq:var}
\mbP(\mE(\epsilon)^c)  
&\leq& \dfrac{\mbE \, \norm{\bY}_1 + 2 \, \sum_{\{i,j\} \prec \{v,w\} \subset \mN} \, \cov(Y_{i,j}, \, Y_{v,w})}{\epsilon^2}  \s \\
\= \dfrac{\mbE \, \norm{\bY}_1 + 2 \, \left[ D_{g} \right]^{+}}{\epsilon^2}.  
\ee 
Taking $\epsilon = \mbE \, \norm{\bY}_1 + 2 \, \left[ D_{g} \right]^{+} > 0$
shows that 
$\mbP(\mE(\epsilon)^c) \leq (\mbE \, \norm{\bY}_1)^{-1}$ and 
\beno
\mbP\left( \mD_{2}(\nat, t) \,|\, \mE(\epsilon) \right)
&\leq& \exp\left(\dfrac{-t^2}{36\,\widetilde{\lambda}_{\max}^{\star} \, (\mbE \, \norm{\bY}_1 + [D_{g}]^{+}) \, + \, 2\, \sqrt{p} \, t } \, + \, p \, \log \, 5 \right).
\ee
Combining all results shows that 
\beno
\mbP\left(\norm{\GradL - \mbE \, \GradL}_{2} \geq t \right)
\ee
is bounded above by 
\beno
\exp\left(\dfrac{-t^2}{36\,\widetilde{\lambda}_{\max}^{\star} \, (\mbE \, \norm{\bY}_1 + [D_{g}]^{+}) \, + \, 2\, \sqrt{p} \, t } \, + \, p \, \log \, 5 \right) \, + \, \dfrac{1}{\mbE \norm{\bY}_1}. 
\ee
As a final matter, 
note that this choice of $\epsilon > 0$ ensures $\mE(\epsilon)$ contains all $\by \in \mbY$ 
with $\norm{\by}_1 \in [0, \, 2 (\mbE \, \norm{\bY} + [D_g]^{+})]$
as the empty graph is an element of $\mbY$ with $0$ edges. 

\qed

\s 

\hide{
We next prove a related result for gradients of log-pseudolikelihood functions 
of multilayer networks in Lemma \ref{lem:concentration_pl}.  
The proof of Lemma \ref{lem:concentration_pl} essentially follows the same proof of Lemma \ref{lem:concentration_likelihood},
and as a result we do not repeat key arguments, instead opting to only outline the changes in the proof.  
}

%% file: suff_exp.tex
\begin{lemma}
\label{lem:s_hetero}
Consider a multilayer network model following the form of equation \eqref{general_model} 
with maximum interaction term $H \leq K$
and 
is defined on a set of $N \geq 3$ nodes and $K \geq 1$ layers.  
Then the following hold: 
\ben
\item The conditional probability mass function of $\bX$ given $\bY$ is an exponential family:  
\beno
\sepmodel(\bX = \bx \mid \bY = \by)
&\propto& h(\bx, \, \by) \; \exp\left( \langle \nat, \, \bs(\bx) \rangle \right),
\ee
where 
\beno
h(\bx, \, \by)
\= \dprod_{\{i,j\} \subset \mN} \, \one(\norm{\bx_{i,j}}_1 > 0)^{y_{i,j}} \; 
\one(\norm{\bx_{i,j}}_1 = 0)^{1 - y_{i,j}},
\ee
the sufficient statistic vector $s : \mbX \mapsto \mbR^p$ and the natural parameter vector $\nat \in \mbR^p$. 

\item For each $l \in \{1, \ldots, p\}$,
there exists $h \in \{1, \ldots, H\}$ and $\{k_1, \ldots, k_h\} \subseteq \{1, \ldots, K\}$ such that 
the $l^{\text{th}}$ component of the sufficient statistic vector $s(\bx)$ can be written as 
\be
\label{eq:suff}
s_l(\bx)
\= \dsum_{\{i,j\} \subset \mN} \, s_{l,i,j}(\bx)
\= \dprod_{r=1}^{h} \, x_{i,j}^{(k_r)}. 
\ee
\item The exponential family outlined above is minimal, full, and regular. 
\een
\end{lemma}

\s 

\llproof \ref{lem:s_hetero}.
First, 
the form of the conditional probability distribution of $\bX$ given $\bY$
derived in Proposition \ref{prop:inference} is given by 
\be
\label{eq:1999}
\mbP_{\nat}(\bX = \bx \,|\, \bY = \by)
\= \exp\left( \log \, f(\bx, \nat) + \log \, \psi(\nat, \by) \right), 
\ee
provided $h(\bx, \by) = 1$.
The form of \eqref{general_model} 
suggests that \eqref{eq:1999} will be a minimal exponential family in canonical form
due to the form of the Markov random field specification for $f(\nat, \bx)$
and the definition of $\psi(\nat, \by)$. 
From the form of $f(\bx, \nat)$ in \eqref{general_model},
\beno
\scalebox{0.9}{$\log f(\bx, \nat)$}
\,=\, \scalebox{0.8}{$\dsum_{\{i,j\} \subset \mN} \,
\left(\dsum_{k=1}^K \theta_{k} x_{i,j}^{(k)} 
+  \dsum_{\substack{k < l}}^{K}  \theta_{k,l} x_{i,j}^{(k)} x_{i,j}^{(l)} + \ldots  
+ \dsum_{k_1 <\ldots< k_H}^{K} \theta_{k_1,k_2,\ldots,k_H}  x_{i,j}^{(k_1)} \cdots x_{i,j}^{(k_H)} \right)$},
\ee
where $H \le K$ is the highest order of cross-layer interactions included in the model.
We write $\theta_{k_1,k_2,\ldots,k_h}$ to reference the $h$-order interaction parameter
for the interaction term among layers $\{k_1, \ldots, k_h\} \subseteq \{1, \ldots, K\}$.
As specified, 
$\psi(\nat, \by)$ is the normalizing constant for the exponential family. 
As such, 
the natural parameter space of the exponential family is $\mbR^p$ 
as the support of $\mbX$ is finite, 
which implies $\psi(\nat, \by) < \infty$ for all $\nat \in \mbR^p$ and $\by \in \mbY$.
We establish minimality by noting that the components of the parameter vector $\nat$ 
satisfy no linear or affine constraints.  
Attached to each parameter $\theta_{k_1,\ldots,k_h}$ 
($\{k_1,\ldots, k_h\} \subseteq \{1, \ldots K\}$, 
$h \in \{1, \ldots, H\}$) 
is the sufficient statistic 
\beno
s_{k_1, \ldots, k_h}(\bx) 
\= \dsum_{\{i,j\} \subset \mN} \, x_{i,j}^{(k_1)} \,\cdots\, x_{i,j}^{(k_h)}.
\ee
Each statistic $s_{k_1, \ldots, k_h}$ is a function of distinct, non-degenerate random variables,
provided $\norm{\by}_1 > 0$,
and so none of the statistics $s_{k_1, \ldots, k_h}$ satisfy any linear or affine constraints. 
Hence, 
\eqref{general_model}
specifies a minimal and full exponential family with natural parameter space $\mbR^p$
of dimension $p = \sum_{h=1}^{H} \, \binom{K}{h}$
and sufficient statistic vector $s(\bx)$ 
with components 
$s_{k_1, \ldots, k_h}(\bx)$ 
($\{k_1, \ldots, k_h\} \subseteq \{1, \ldots, K\}, h = 1, \ldots, H$). 
Regularity follows trivially 
\citep[e.g., Proposition 3.7, pp. 28,][]{Su19}. 
The form of \eqref{eq:suff} outlines this for a linear indexing of the components of the sufficient statistic vector. 

\qed

\s\s
\hide{
\begin{lemma}
\label{lem:exp_pseudo}
Consider multilayer networks satisfying \eqref{general_model} 
with maximum interaction term $H \leq K$ 
and defined on a set of $N \geq 3$ nodes and $K \geq 1$ layers. 
Then the conditional probability mass function of $X_{i,j}^{(k)}$ given 
$\bY = \by$ and $\bX_{i,j}^{(-k)} = \bx_{i,j}^{(-k)}$ is an exponential family 
\beno
\sepmodel(X_{i,j}^{(k)} = x_{i,j}^{(k)} \,|\, \bX_{i,j}^{(-k)} = \bx_{i,j}^{(-k)}, \bY = \by)
&\propto& h(\bx, \, \by) \; \exp\left( \langle \nat, \, \bs(\bx) \rangle \right),
\ee
with sufficient statistic vector $\bs : \mbX \mapsto \mbR^p$ defined in Lemma \ref{lem:s_hetero}, 
natural parameter vector $\nat \in \mbR^p$,  
and 
\beno
h(\bx, \, \by)
\= \dprod_{\{i,j\} \subset \mN} \, \one(\norm{\bx_{i,j}}_1 > 0)^{y_{i,j}} \; 
\one(\norm{\bx_{i,j}}_1 = 0)^{1 - y_{i,j}}.  
\ee
\end{lemma}

\s 

\llproof \ref{lem:exp_pseudo}. First,  
note that the form of \eqref{general_model} and Proposition \ref{prop:inference} suggests that 
\be
 \sepmodel(X_{i,j}^{(k)} = x_{i,j}^{(k)} \,|\, \bX_{i,j}^{(-k)} = \bx_{i,j}^{(-k)}, \bY = \by) \s \\
 \quad\quad\quad = \dfrac{h(x_{i,j}^{(k)}, \, \bx_{i,j}^{(-k)}, \, \by) \; \exp\left( \langle \nat, \, \bs(x_{i,j}^{(k)}, \, \bx_{i,j}^{(-k)}) \rangle \right)}{\dsum_{x_{i,j}^{(k)} \in \{0,1\} }\, h(x_{i,j}^{(k)}, \, \bx_{i,j}^{(-k)}, \, \by) \; \exp\left( \langle \nat, \, \bs(x_{i,j}^{(k)}, \, \bx_{i,j}^{(-k)}) \rangle \right)}
\ee
is an exponential family in canonical form
using the Markov random field specification for $f(\nat, \bx)$ 
and the form of the conditional probability distribution of $X_{i,j}^{(k)}$ given $\bY$ 
and $\bX_{i,j}^{(-k)}$.  
However,
this exponential family may not be full rank due to possible 0 values of components of the given ($K$-$1$)-dimensional vector $\bx_{i,j}^{(-k)}$ and thus may not be minimal. 

\qed
}

%% file: proof_thm1.tex
By Proposition \ref{prop:inference}, 
observing $\bX = \bx$ implies we observe $\bY = \by$,
as for each given $\bx \in \mbX$,
$\bY = \by$ ($\mbP$-a.s.) for one and only one $\by \in \mbY$ given by 
\beno
y_{i,j} 
\= \one\left( \norm{\bx_{i,j}}_1 \,>\, 0 \right),
&& \{i,j\} \subset \mN. 
\ee 
Denote the gradient of $-\ell(\nat; \bx, \by)$ by 
\beno
\gradL 
&\coloneqq& -\nabla_{\nat} \, \ell(\nat; \bx, \by)
\ee
and the expected Hessian matrix of the negative log-likelihood by 
\beno
\bH(\nat) 
&\coloneqq& -\mbE \, \nabla_{\nat}^2 \, \ell(\nat; \bX, \bY). 
\ee 
Theorem 6.3.4 of \citet{OrRh20} states that if 
\beno 
\label{eq:ex_cond} 
(\nat - \truth)^{\top} \, \gradL
&\geq& 0 
\ee
for all $\nat \in \partial \, \mB_2(\truth, \epsilon)$, 
where $\partial \, \mB_2(\truth, \epsilon)$ is the boundary of the set 
\beno
\mB_2(\truth, \epsilon) 
\= \{ \nat \in \mbR^p \,:\, \norm{\nat - \truth}_2 < \epsilon\},
\ee 
then $\gradL$ 
has a root in $\mB_2(\truth, \epsilon) \cup \partial \mB_2(\truth, \epsilon)$,
i.e.,
$\mle$ exists and satisfies $\norm{\mle - \truth}_2 \leq \epsilon$.
Note that a root of $\gradL$ is also a root of $-\gradL$;
in what follows, 
we consider finding a maximizer of $\ell(\nat;\bx,\by)$ by finding a minimizer of $-\ell(\nat;\bx,\by)$. 
The classification of roots as maximizers/minimizers 
is justified from the fact that that $\ell(\nat; \bx, \by)$ is concave in $\nat$,
a fact which follows from Proposition \ref{prop:inference},
as $g(\by)$ is constant in $\nat$ and $\log \, \mbP_{\nat}(\bX = \bx \,|\, \bY = \by)$ 
is the log-likelihood of a minimal, full, and regular exponential family with natural parameter vector $\nat$
and thus is strictly concave in $\nat$ \citep[Proposition 3.10, p. 32,][]{Su19}. 
By the multivariate mean-value theorem \citep[][Theorem 5]{FuMa91}, 
\beno
(\nat - \truth)^{\top}  \mbE \, \GradL 
\= (\nat - \truth)^{\top}  \mbE \nabla_{\truth}(\bX, \bY) 
+ (\nat - \truth)^{\top}  \bH(\dot\nat)  (\nat - \truth) \s \\ 
\=(\nat - \truth)^{\top}  \bH(\dot\nat)  (\nat - \truth), 
\ee
where $\dot\nat = t \, \nat + (1 - t) \, \truth$ (some $t \in [0, 1]$) 
and by invoking Lemma 2 of \citet{StSc21},
which shows that both the expected log-likelihood and log-pseudolikelihood 
of a minimal exponential family is uniquely maximized at the data-generating parameter vector $\truth$, 
implying $\mbE \, \nabla_{\truth}(\bX, \bY) = 0$. 
Let $\gamma \in (0, \epsilon)$
and arbitrarily take $\nat \in \partial \mB_2(\truth, \gamma)$. 
Then 
\beno
(\nat - \truth)^{\top} \bH(\dot\nat) (\nat - \truth)
\= \dfrac{(\nat - \truth)^{\top}  \bH(\dot\nat)  (\nat - \truth)}{(\nat - \truth)^{\top} (\nat- \truth)}  
\norm{\nat - \truth}_2^2\s \\
&\geq& \gamma^2 \, \lambda_{\min}(\bH(\dot\nat)),  
\ee
since $\norm{\nat - \truth}_2 = \gamma$ as $\nat \in \partial \mB_2(\truth, \gamma)$ 
and because the Rayleigh quotient of a matrix is bounded below by the smallest eigenvalue of that matrix so that
\beno
\dfrac{(\nat - \truth)^{\top}  \bH(\dot\nat)  (\nat - \truth)}{(\nat - \truth)^{\top} (\nat- \truth)}
\;\geq\; \lambda_{\min}(\bH(\dot\nat))
\;\geq\; \inf\limits_{\nat \in \mB_2(\truth, \epsilon)} \, \lambda_{\min}(\bH(\nat)),
\ee
where $\lambda_{\min}(\bH(\dot\nat))$ is the smallest eigenvalue of $\bH(\dot\nat)$,
noting that 
\beno
\norm{\dot\nat - \truth}_2 
\= \norm{t \, \nat + (1 - t) \, \truth - \truth}_2
\= t \, \norm{\nat - \truth}_2 
&\leq& \epsilon,  
\ee
since $t \in [0, 1]$. 
Lemma \ref{lem:min-eig} showed that 
\beno
\lambda_{\min}(\bH(\nat)) 
\= \lambda_{\min}(\mcI(\nat)) \; \mbE \, \norm{\bY}_1, 
\ee
which in turn implies 
\beno
\inf\limits_{\nat \in \mB_2(\truth, \epsilon)} \, \lambda_{\min}(\bH(\nat))
\= \widetilde{\lambda}_{\min}^{\epsilon} \; \mbE \, \norm{\bY}_1,
\ee
where 
\beno
\widetilde{\lambda}_{\min}^{\epsilon}
\;\coloneqq\; \inf\limits_{\nat \in \mB_2(\truth, \epsilon)} \, \lambda_{\min}(\mcI(\nat)) ,
\ee
with $\mcI(\nat)$ defined in Lemma \ref{lem:min-eig}. 
Hence,
for $\nat \in \partial \mB_2(\truth, \gamma)$ ($\gamma \in (0, \epsilon)$),  
\be
(\nat - \truth)^{\top}  \mbE \, \GradL
&\geq& \gamma^2 \;\widetilde{\lambda}_{\min}^{\epsilon} \; \mbE \, \norm{\bY}_1 \; \ge \; 0. 
\ee
We next turn to showing  
\beno
\mbP\left( \inf\limits_{\nat \in \mB_2(\truth,\gamma)} \, 
(\nat - \truth)^{\top} \, \GradL \,\geq\, 0 \right) 
&\geq& 1 - \, \exp \, (-2\,p) \, - \, (\mbE \norm{\bY}_1)^{-1},
\ee
by showing that the event  
\beno
\sup\limits_{\nat \in \mB_2(\nat^\star, \gamma)} \, 
|(\nat - \truth)^{\top} \, ( \mbE \, \GradL - \GradL)|
\;<\; \gamma^2 \; \Lam
\ee
occurs with probability at least $1 - \, \exp \, (-2\,p) \, - \, (\mbE \norm{\bY}_1)^{-1}$. This will 
in turn imply that the event that   
$\norm{\mle - \truth}_2 \; \leq \; \gamma$ will
happen with probability at least $1 - \, \exp \, (-2\,p) \, - \, (\mbE \norm{\bY}_1)^{-1}$. 
Applying the Cauchy-Schwarz inequality and utilizing standard vector norm inequalities, for $\nat \in \partial \mB_2(\truth, \gamma)$, we have
\beno
|(\nat - \truth)^{\top} \, ( \mbE \, \GradL - \GradL)| \s \\
\leq \norm{\nat - \truth}_2 \, \norm{\GradL  - \mbE \, \GradL}_2 \s \\ 
=  \gamma \; \norm{\GradL - \mbE \, \GradL}_{2}.
\ee
Therefore, it suffices to demonstrate,
for all $\nat \in \partial \mB_2(\truth, \gamma)$,
\beno
\mbP\left(\norm{\GradL - \mbE \, \GradL}_{2} \,<\, \gamma \, \Lam \right)
\ee
is bounded below by
\beno
1 - \, \exp \, (-2\,p) \, - \, (\mbE \norm{\bY}_1)^{-1}. 
\ee
For ease of presentation,
we define $\mD_{N,\gamma,p}$ to be the event 
\beno
\norm{\GradL - \mbE \, \GradL}_{2} 
&\geq& \gamma \, \Lam. 
\ee
Applying Lemma \ref{lem:concentration_likelihood}, the probability $\mbP\left(\mD_{N,\gamma,p} \right)$ is bounded above by
\be
\label{lem2_bound}
\scalebox{0.9}{$\exp\left( -\dfrac{(\gamma \, \Lam)^2}{36 \, \widetilde{\lambda}_{\max}^{\star} \, (\mbE \, \norm{\bY}_1 + [D_{g}]^{+}) \,  + 2 \, \sqrt{p} \, \gamma \, \Lam} \, + \, p \, \log \, 5 \right) + \dfrac{1}{\mbE \norm{\bY}_1}$},
\ee
recalling $\widetilde{\lambda}_{\max}^{\star} = \lambda_{\max}(\mcI(\truth))$, 
$[D_{g}]^{+} \coloneqq \max\{0, \, D_{g}\}$,
and 
\beno
D_{g}
&\coloneqq& \dsum_{\{i,j\} \prec \{v,w\} \subset \mN} \, \cov(Y_{i,j}, \, Y_{v,w}),
\ee
where $\{i,j\} \prec \{v,w\}$ implies the sum is taken with respect to the lexicographical ordering
of pairs of nodes.
Choose
\beno
\gamma 
\= \beta \, \sqrt{\dfrac{ p \, \widetilde{\lambda}_{\max}^{\star}}{\mbE \norm{\bY}_1}} \; \dfrac{1}{\widetilde{\lambda}_{\min}^{\epsilon}},
\ee
where $\beta > 0$ is a positive constant independent of $N$ and $p$ whose value will be determined later.
If 
\beno
\lim\limits_{N \to \infty} \, 
\beta \, \sqrt{\dfrac{ p \, \widetilde{\lambda}_{\max}^{\star}}{\mbE \norm{\bY}_1}} \; \dfrac{1}{\widetilde{\lambda}_{\min}^{\epsilon}}
\= 0,
\ee
then for $N$ sufficiently large,
we will have $\gamma < \epsilon$,
which ensures $\epsilon$ may be chosen independent of $N$ and $p$.  
While $\epsilon$ can be chosen independent of $N$ and $p$, 
note that $p$ is expected to be a function of $N$ and thus $\widetilde{\lambda}_{\min}^{\epsilon}$ 
will not (in general) be independent of $N$,
possibly holding implications for how fast $p$ may grow with $N$ for certain $\truth$ and $\epsilon$. 
This choice of $\gamma$ in turn implies that the first term of the exponent in \eqref{lem2_bound} becomes
\be
\label{reduced_exp}
 \exp\left( -\dfrac{\beta^2\, \widetilde{\lambda}_{\max}^{\star} \, \mbE\,\norm{\bY}_1 \, p }{36 \,\widetilde{\lambda}_{\max}^{\star}\, (\mbE\,\norm{\bY}_1+[D_g]^+) \, + \, 2 \,  \beta \, p \, \sqrt{\mbE\,\norm{\bY}_1\,\widetilde{\lambda}_{\max}^{\star}}} \right).
\ee
Canceling $\mbE\,\norm{\bY}_1$ in \eqref{reduced_exp} gives
\beno
\exp\left( -\dfrac{\beta^2\, \widetilde{\lambda}_{\max}^{\star}  \, p }{36 \,\widetilde{\lambda}_{\max}^{\star} \, \left(1+\dfrac{[D_g]^+}{\mbE\,\norm{\bY}_1} \right) \, + \, 2 \,  \beta \,p \, \sqrt{\dfrac{\widetilde{\lambda}_{\max}^{\star}}{\mbE\,\norm{\bY}_1} }} \right).
\ee
By Assumption \ref{assump2},
\beno
p \leq \sqrt{\widetilde{\lambda}_{\max}^{\star} \, \mbE\, \norm{\bY}_1} \, \left(1+\dfrac{[D_g]^+}{\mbE\,\norm{\bY}_1} \right),
\ee
and by Assumption \ref{assump1},
\beno
\dfrac{[D_g]^+}{\mbE\,\norm{\bY}_1}  \; \leq \; C_0,
\ee
where $C_0 > 0$ is a positive constant independent of $N$ and $p$,
we have 
\beno
p \,\sqrt{\dfrac{\widetilde{\lambda}_{\max}^{\star}}{\mbE\,\norm{\bY}_1} }
\; \leq \; 
\widetilde{\lambda}_{\max}^{\star} \, \left(1+\dfrac{[D_g]^+}{\mbE\,\norm{\bY}_1} \right)
\; \leq \; \widetilde{\lambda}_{\max}^{\star} \, \left(1+C_0 \right).
\ee
As a result, the upper bound in \eqref{lem2_bound} reduces to 
\be
\label{reduced_upper_bound}
\mbP\left(\mD_{N,\gamma,p} \right)
\;\leq\; \exp \left( \left(\dfrac{-\beta^2}{36(1+C_0)+2(1+C_0)\beta}  + \log 5 \right) p\right) + \dfrac{1}{\mbE\,\norm{\bY}_1}.
\ee
To obtain the desired convergence rate, require
\be
\label{quad_eq}
\dfrac{-\beta^2}{18\,C+C \, \beta}  + \log 5 = -2,
\ee
where $C = 2(1+C_0)$. 
The constant $\beta$ can be solved by using the quadratic formula and the positive root is given by 
\beno
\beta \= \dfrac{C\,\log\,5 \, -2\, C \; + \; \sqrt{(C\,\log\,5 \, -2 \, C )^2 \; + \; 72\,(C\,\log\,5 - 2\, C )}}{2},
\ee
which ensures $\mbP\left(\mD_{N,\gamma,p} \right) \leq \exp \, (-2\,p) + (\mbE\norm{\bY}_1)^{-1}$.  
We have thus shown,
for all $\nat \in  \partial \mB_2(\truth, \gamma)$, 
that 
\beno
\mbP\left(\norm{\GradL - \mbE \, \GradL}_{2} \,\leq\, \gamma \, \Lam \right)
\ee
is bounded below by
\beno
1 - \exp \, (-2\,p)  \, - \, (\mbE \, \norm{\bY}_1)^{-1}, 
\ee
under the above conditions. 
As a result, 
there exists $N_0 \geq 3$ such that,
for all $N \geq N_0$  
and with probability at least $1 - \exp \, (-2\,p)  \, - \, (\mbE \, \norm{\bY}_1)^{-1}$, 
the set $\Mle$ is non-empty and the unique element of the set $\mle \in \Mle$ 
satisfies (uniqueness following from minimality, as discussed in Section \ref{sec3}) 
\beno
\norm{\mle - \truth}_2 &\leq& 
C \; \dfrac{\sqrt{\widetilde{\lambda}_{\max}^{\star}} }{\widetilde{\lambda}_{\min}^{\epsilon}} \; \sqrt{\dfrac{p}{\mbE \norm{\bY}_1}}.
\ee

\s\s
\subsection{Proof of Corollary \ref{corollary}}
We prove Corollary \ref{corollary} from Section \ref{sec3}. Under the same assumptions as Theorem \ref{thm1} and in the case that the parameter dimension $p$ is fixed, the proof of Corollary \ref{corollary} remains unchanged from that of Theorem \ref{thm1} except that the exponent in equation \eqref{reduced_upper_bound} scales with $N$ as opposed to $p$ in Theorem \ref{thm1}. Following the same notations in the proof of Theorem \ref{thm1}, we rewrite equation \eqref{quad_eq} as
\be
\label{beta_gamma}
 \dfrac{-\beta^2}{18\,C+C \, \beta}  + \log 5 \= -\eta \, (N),
\ee
where $\eta : \mbN^+ \mapsto \mbR^+ $ is an increasing function of $N$. For the ease of notation, we write $\eta_N$ instead of $\eta(N)$ in the rest of the proof.
For any $\alpha_N \in (2\,(\mbE\,\norm{\bY}_1)^{-1}, 1/2)$, let
\beno
\exp \, \left( -\eta_N \, p\right) \= \dfrac{\alpha_N}{2}.
\ee
Note that as $N$ goes to infinity, $\alpha_N$ is allowed to approach 0 through the increasing of $\eta_N$.
Then the upper bound of $\mbP\left(\mD_{N,\gamma,p} \right)$ given in \eqref{reduced_upper_bound} becomes
\beno
\mbP\left(\mD_{N,\gamma,p} \right)
&\leq&  \dfrac{\alpha_N}{2} \; + \;\dfrac{1}{\mbE\,\norm{\bY}_1} &\leq& \alpha_N,
\ee
where the last inequality follows from $ \alpha_N \ge 2\,(\mbE\,\norm{\bY}_1)^{-1}$.
To obtain the desired result,
solve the positive root of $\beta$ in terms of $\eta_N$ from \eqref{beta_gamma}:
\beno
\beta \= \dfrac{C\,\log\,5 \; + \; \eta_N\, C \; + \; \sqrt{(C\,\log\,5 \; + \;\eta_N \, C )^2 \; + \; 72\,(C\,\log\,5 \;+\; \eta_N\, C )}}{2},
\ee
and write $\eta_N$ in terms of $\alpha_N$, where $\alpha_N \to 0$ and $\eta_N \to \infty$ as $N \to \infty$:
\beno
\eta_N \= -\,\dfrac{\log\,(\alpha_N \, / \,2)}{p}.
\ee
Let 
\beno
A_1 \= C\, \log\,5, && A_2 \= \dfrac{C}{p}.
\ee
Then for $\alpha_N \in (2\,(\mbE\,\norm{\bY}_1)^{-1}, 1/2)$,
\beno
\beta \= \scalebox{0.82}{$\dfrac{A_1 \; - \; A_2 \, \log \, (\alpha_N \, / \,2) \; + \; \sqrt{(A_1 \; - \; A_2 \, \log \, (\alpha_N \, / \,2))^2 \; + \; 72\,(A_1 \; - \; A_2\,\log\,(\alpha_N \, / \,2))}}{2}$} \s \\
\= \scalebox{0.66}{$\dfrac{\log\,\left(\dfrac{\alpha_N}{2}\right) \; \left(\dfrac{A_1}{\log\,(\alpha_N \, / \,2)} \; - \; A_2 \;+ \; \sqrt{\left(\dfrac{A_1}{\log\,(\alpha_N \, / \,2)} \; - \; A_2 \right)^2 \; + \; 72\,\left(\dfrac{A_1}{(\log\,(\alpha_N \, / \,2))^2} \; - \; \dfrac{A_2}{\log\,(\alpha_N \, / \,2)} \right)}\right) }{2}$} \s \\
& \leq & 
\scalebox{0.73}{$\dfrac{\log\,\left(\dfrac{\alpha_N}{2}\right)\; \left( \dfrac{A_1}{\log\,0.25} \; - \; A_2 \;+ \; \sqrt{\left(\dfrac{A_1}{\log\,0.25} \; - \; A_2 \right)^2 \; + \; 72\,\left(\dfrac{A_1}{(\log\,(0.25))^2} \; - \; \dfrac{A_2}{\log\,(0.25)} \right)}\right) }{2}$} \s \\
\= A \, \left| \log \, \left(\dfrac{\alpha_N}{2}\right) \right|,
\ee
where 
\beno
A \=\scalebox{0.8}{$\dfrac{ \left| \;\dfrac{A_1}{\log\,0.25} \; - \; A_2 \;+ \; \sqrt{\left(\dfrac{A_1}{\log\,0.25} \; - \; A_2 \right)^2 \; + \; 72\,\left(\dfrac{A_1}{(\log\,(0.25))^2} \; - \; \dfrac{A_2}{\log\,(0.25)} \right)}\;\right| }{2}$}.
\ee
As a result, when $p$ is fixed, we showed that for $\alpha_N \in (2\,(\mbE\,\norm{\bY}_1)^{-1}, 1/2)$, with probability at least $1-\alpha_N$,
\beno
\norm{\mle - \truth}_2 &\leq&
A^\prime \; | \, \log\, (\alpha_N \, / \,2) \, |  \; \dfrac{\sqrt{\widetilde{\lambda}_{\max}^{\star}} }{\widetilde{\lambda}_{\min}^{\epsilon}} \; \sqrt{\dfrac{1}{\mbE \norm{\bY}_1}},
\ee
where $A^\prime = A \, \sqrt{p}$ is a positive constant independent of $N$.
\s \s

\section{Proof of Theorem \ref{thm:minimax} and Corollary \ref{cor:minimax}}
\label{sec:pf_minimax}
We prove Theorem \ref{thm:minimax} and Corollary \ref{cor:minimax} from Section \ref{sec3} in one chapter. We first use Fano's method outlined in Chapter 15.3 of \citet{high-dim-stat} and the Kullback-Leibler divergence to derive the lower bound of the minimax risk for multilayer network models specified in \eqref{general_model}.
Let $\epsilon > 0$ be fixed and consider $\gamma \in (0,\epsilon)$. For $M\ge2$ and some $\delta > 0$, let $\{\nat_1, \ldots, \nat_M\} \subset \mB_2(\truth,\gamma)$ be a $2\delta$-separated set. We then have $\norm{\nat_i-\nat_j}_2 \ge 2\delta$ for any pair $\{i,j\} \subseteq \{1,\ldots,M\}$. Define the Kullback-Leibler divergence of $\nat_i$ and $\nat_j$ by
\beno
\text{KL}(\nat_i,\nat_j) &\coloneqq& \dsum_{\bx \in \mbX}\, \varphi_{\nat_i}(\bx) \, \log \, \dfrac{\varphi_{\nat_i}(\bx)}{\varphi_{\nat_j}(\bx)}, && \{i,j\} \subseteq \{1,\ldots,M\},
\ee
where $\varphi_{\nat}(\bx)$ belongs to a minimal exponential family defined in Proposition \ref{prop:inference}:
\beno
\varphi_{\nat}(\bx) & \coloneqq &\mbP_{\nat}(\bX = \bx \mid \bY = \by)
\= \exp\,(\log f(\bx, \nat) + \log \psi(\nat, \by)),
\ee
recalling $ f(\bx, \nat)$ and $\psi(\nat, \by)$ follow the same form of \eqref{general_model}.
For $\nat \in \mbR^p$, denote by $\bs(\bX) \in \mbR^p$ the sufficient statistic vector of the exponential family $\varphi_{\nat}(\bx)$. Then the Kullback-Leibler divergence can be written as
\be
\label{KL}
\scalebox{0.9}{$\text{KL}(\nat_i,\nat_j)$} \= \scalebox{0.9}{$\dsum_{\bx \in \mbX} \, \varphi_{\nat_i}(\bx) \, \left[ \langle\,\nat_i - \nat_j, \; \bs(\bx) \, \rangle + \log\,\psi(\nat_i,\by) - \log\,\psi(\nat_j,\by) \right]$} \s \\
\= \scalebox{0.9}{$\mbE_{\nat_i}\,\langle\,\nat_i - \nat_j, \; \bs(\bX) \, \rangle  + \log\,\psi(\nat_i,\by) - \log\,\psi(\nat_j,\by)$} \s \\
\= \scalebox{0.9}{$\langle\,\nat_i - \nat_j, \; \bmu(\nat_i) \, \rangle  + \log\,\psi(\nat_i,\by) - \log\,\psi(\nat_j,\by)$},
\ee
where $\bmu(\nat) \coloneqq \mbE_{\nat}\,\bs(\bX)$ is the mean-value parameter map of the exponential family. By Corollary 2.3 of \citet{Br86},
\be
\label{log_norm_const}
\scalebox{0.95}{$\log\,\psi(\nat_j)$} \= \scalebox{0.85}{$\log\,\psi(\nat_i) + \langle\,\nat_j - \nat_i, \; -\bmu(\nat_i) \, \rangle - \dfrac{1}{2}\, \langle\,\nat_j - \nat_i, \; \mcI_{\bX}(\dot\nat) \, (\nat_j - \nat_i) \, \rangle$} \s \\
\= \scalebox{0.85}{$\log\,\psi(\nat_i) + \langle\,\nat_i - \nat_j, \; \bmu(\nat_i) \, \rangle - \dfrac{1}{2}\, \langle\,\nat_i - \nat_j, \; \mcI_{\bX}(\dot\nat)\, (\nat_i - \nat_j) \, \rangle$}, 
\ee
where $\dot\nat = t\nat_i + (1-t)\nat_j$ for some $t \in (0,1)$, and $\mcI_{\bX}(\dot\nat)$ is the Fisher information matrix at $\dot\nat$ for $\bX \in \mbX$. For a fixed $\epsilon > 0$ such that $\gamma \in (0,\epsilon)$ and $\{\nat_1, \ldots, \nat_M\} \subset \mB_2(\truth,\gamma)$, define
\beno
\widetilde{\lambda}_{\max}^{\epsilon}
&\coloneqq & 
\sup\limits_{\nat \in \mB_2(\truth, \epsilon)} \,\dfrac{\lambda_{\max}(\mcI_{\bX}(\nat)) }{\mbE\,\norm{\bY}_1} \= \sup\limits_{\nat \in \mB_2(\truth, \epsilon)} \, \lambda_{\max} \, (\mcI(\nat)),
\ee
where $\lambda_{\max}(\bA)$ is the maximum eigenvalue of matrix $\bA$ and $\mcI(\nat)$ is the Fisher information matrix for an activated dyad defined in Lemma $\ref{lem:min-eig}$.  
Combining \eqref{KL} and \eqref{log_norm_const} and using the standard matrix norm inequality and the triangle inequality, we have
\beno
\text{KL}(\nat_i,\nat_j) \= \dfrac{1}{2} \, \langle\,\nat_i - \nat_j, \; \mcI_{\bX}(\dot\nat) \, (\nat_i - \nat_j) \, \rangle \s \\
&\leq& \dfrac{1}{2} \, \mbE\,\norm{\bY}_1\,\widetilde{\lambda}_{\max}^{\epsilon} \, \norm{\nat_i - \nat_j}_2^2 \s \\
&\leq& \dfrac{1}{2} \, \mbE\,\norm{\bY}_1\,\widetilde{\lambda}_{\max}^{\epsilon} \, (\norm{\nat_i - \truth}_2 + \norm{\nat_j - \truth}_2)^2  \s \\ 
&\leq & 2 \, \epsilon^2\,  \mbE\,\norm{\bY}_1\,\widetilde{\lambda}_{\max}^{\epsilon}.
\ee
Note that the size $M$ of the largest possible $2\,\delta$-separated set $\{\nat_1,\ldots,\nat_M\} \subset \mB_2(\truth,\gamma) \subset \mbR^p$ is the packing number of $\mB_2(\truth,\gamma)$. By Lemma 4.2.8 and Corollary 4.2.13 of \citet{Vershynin18}, we have
\beno
M &\geq& \left(\dfrac{\gamma}{2\,\delta}\right)^p,
\ee
and
\beno
\log \, M &\geq& p \,\log \, \left(\dfrac{\gamma}{2\,\delta}\right).
\ee
By Proposition 15.12 of \citet{high-dim-stat}, the minimax risk $\mcR_N$ has the lower bound
\beno
\mcR_N &\geq& \delta \, \left[ 1-\dfrac{\mF + \log\,2}{\log\,M}\right],
\ee
where 
\beno
\mF &\coloneqq& \max\limits_{\{i,j\} \subseteq \{1,\ldots,M\}} \; \text{KL}(\nat_i,\nat_j).
\ee
Since $\gamma \in (0,\epsilon)$, the lower bound for $\mcR_N$ can be written as
\beno
\mcR_N &\geq& \delta \, \left[ 1-\dfrac{2 \, \gamma^2\, \mbE\,\norm{\bY}_1\, \widetilde{\lambda}_{\max}^{\epsilon} + \log\,2}{p \,\log \, \left(\gamma / 2\,\delta\right)}\right].
\ee
To obtain the desired lower bound $\mcR_N \ge \delta/2$, we need
\beno
\dfrac{2 \, \gamma^2\, \mbE\,\norm{\bY}_1\, \widetilde{\lambda}_{\max}^{\epsilon} + \log\,2}{p \,\log \, \left(\gamma / 2\,\delta\right)} & \leq & \dfrac{1}{2},
\ee
which implies
\beno
\dfrac{4\,\gamma^2\,\mbE\,\norm{\bY}_1\,\widetilde\lambda_{\max}^{\epsilon}}{p} + \dfrac{2\,\log\,(2)}{p} &\leq& \log(\gamma/2) - \log\,(\delta).
\ee
Exponentiating both sides we have
\beno
\exp\,\left( \dfrac{4\,\gamma^2\,\mbE\,\norm{\bY}_1\,\widetilde\lambda_{\max}^{\epsilon}}{p} +  \dfrac{2\,\log\,(2)}{p}\right) &\leq& \dfrac{\gamma/2}{\delta}.
\ee
This leads us to the following inequality
\beno
\delta &\leq& \dfrac{\gamma}{2} \, \exp\,\left( -\,\dfrac{4\,\gamma^2 \, \mbE\,\norm{\bY}_1\,\widetilde\lambda_{\max}^{\epsilon} }{p} -\dfrac{2\,\log\,(2)}{p} \right).
\ee
Choosing 
\beno
\gamma \= 2 \, C\,\sqrt{\dfrac{p}{\widetilde\lambda_{\max}^{\epsilon} \, \mbE\,\norm{\bY}_1}}
\ee
for some $C>0$, we obtain the bound
\be
\label{delta_upper_bound}
\delta &\leq& C \, \exp \left(- 16 \, C^2 - \dfrac{2\,\log\,(2)}{p} \right) \, \sqrt{\dfrac{p}{\widetilde\lambda_{\max}^{\epsilon} \, \mbE\,\norm{\bY}_1}}.
\ee
As long as $p = O\,(\widetilde\lambda_{\max}^{\epsilon} \,\mbE\,\norm{\bY}_1)$, we can choose $C$ to ensure $\gamma \in (0,\epsilon)$. Finally, for all $\delta>0$ satisfying \eqref{delta_upper_bound}, we have $\mcR_N$ lower bounded by
\beno
\mcR_N &\ge& \dfrac{\delta}{2}.
\ee
Note that as $p \ge 1$,
\beno
\exp \left(- 16 \, C^2 - \dfrac{2\,\log\,(2)}{p} \right) & \ge & \exp \left(- 16 \, C^2 - 2\,\log\,(2) \right),
\ee
we may choose 
\beno
\delta \= C \, \exp \left(- 16 \, C^2 - 2\,\log\,(2) \right) \, \sqrt{\dfrac{p}{\widetilde\lambda_{\max}^{\epsilon} \, \mbE\,\norm{\bY}_1}} \s \\
\= A^\prime \, \sqrt{\dfrac{p}{\widetilde\lambda_{\max}^{\epsilon} \, \mbE\,\norm{\bY}_1}}, 
\ee
where $A^\prime = C \, \exp \left(- 16 \, C^2 - 2\,\log\,(2) \right)$.
Then we obtain the desired lower bound for the minimax risk
\be
\label{minimax_lower}
\mcR_N &\ge& \dfrac{\delta}{2} \= \dfrac{A^\prime}{2} \, \sqrt{\dfrac{p}{\widetilde\lambda_{\max}^{\epsilon} \, \mbE\,\norm{\bY}_1}}.
\ee
Next, we show the lower bound in \eqref{minimax_lower} matches with the upper bound of the $\ell_2$-error of the maximum likelihood estimator $\mle$ provided in Theorem \ref{thm1}. Let $A = A^{\prime}/2$ be an unknown constant independent of $N$ and $p$. 
We have
\beno
\mcR_N &\ge& A \, \sqrt{\dfrac{p}{\widetilde\lambda_{\max}^{\epsilon} \, \mbE\,\norm{\bY}_1}} \s \\
\= A \, \sqrt{\dfrac{\widetilde{\lambda}_{\max}^{\star}}{\widetilde{\lambda}_{\max}^{\star}} } \, \left(\dfrac{\widetilde{\lambda}_{\min}^{\epsilon}}{\widetilde{\lambda}_{\min}^{\epsilon}}\right) \, \dfrac{1}{\sqrt{\widetilde\lambda_{\max}^{\epsilon}}} \, \sqrt{\dfrac{p}{\mbE\,\norm{\bY}_1}}\s \\
\= A \,\dfrac{1}{\sqrt{\widetilde\lambda_{\max}^{\epsilon}}} \, \dfrac{\widetilde{\lambda}_{\min}^{\epsilon}}{\sqrt{\widetilde{\lambda}_{\max}^{\star}}} \, \dfrac{\sqrt{\widetilde{\lambda}_{\max}^{\star}} }{\widetilde{\lambda}_{\min}^{\epsilon}} \; \sqrt{\dfrac{p}{\mbE \norm{\bY}_1}} \s \\
&\geq & A \,\left(\dfrac{\widetilde{\lambda}_{\min}^{\epsilon}}{\widetilde\lambda_{\max}^{\epsilon}} \right)\, \dfrac{\sqrt{\widetilde{\lambda}_{\max}^{\star}} }{\widetilde{\lambda}_{\min}^{\epsilon}} \; \sqrt{\dfrac{p}{\mbE \norm{\bY}_1}},
\ee
where the last inequality holds because $\widetilde\lambda_{\max}^{\epsilon} \,\ge \,\widetilde{\lambda}_{\max}^{\star}$.
Under the assumption that
\beno
\widetilde\lambda_{\max}^{\epsilon}= O\, \left(\widetilde{\lambda}_{\min}^{\epsilon}\right),
\ee
we showed that the lower bound of the minimax risk $\mcR_N$ and the upper bound of the $\ell_2$-error of the maximum likelihood estimator presented in Theorem \ref{thm1} match up to an unknown constant independent of $N$ and $p$.
\qed

%% file: prop2_normal.tex
\begin{proposition}
\label{prop:suff_norm}
Consider a separable multilayer network model following the form of equation \eqref{general_model} and is 
defined on a set of $N \geq 3$ nodes and $K \geq 1$ layers.
Denote by $\bs(\bX) \in \mbR^p$  the sufficient statistic vector of the exponential family $\mbP(\bX = \bx \,|\, \bY = \by)$
as defined in Lemma \ref{lem:s_hetero}. 
Let $\mbE^{\bY}$ be the random conditional expectation operator 
for the distribution  of $\bX$ conditional on $\bY$, and define
\beno
\bS_{\mN} 
&\coloneqq& 
(I(\truth) \, \norm{\bY}_1)^{-1/2} \, (\bs(\bX) - \mbE^{\bY} \bs(\bX)) \s \\ 
\= \dsum_{\{i,j\} \subset \mN} \, (I(\truth) \, \norm{\bY}_1)^{-1/2} \, (\bs_{i,j}(\bX) - \mbE^{\bY} \bs_{i,j}(\bX)). 
\ee
For any measurable convex set $\mA \subset \mbR^p$, 
\beno
\label{ineq normal error homo}
\left| \, \mbP(\bS_\mN\in \mA)-\Phi(\bZ \in \mA)\,  \right| &\leq&
 \dfrac{83}{(\widetilde{\lambda}_{\min}^{\epsilon})^{3/2}} \, 
\sqrt{\dfrac{p^{7/2}}{\mbE \, \norm{\bY}_1}}
+  \dfrac{4}{\mbE \, \norm{\bY}_1} + \dfrac{8 \, \left[D_{g}\right]^{+}}{\left(\mbE \, \norm{\bY}_1 \right)^2}, 
\ee
where $\Phi$ is the standard multivariate normal measure and $\bZ \sim \text{MvtNorm}(\bm{0}_p, \bI_p)$,
where $\bm{0}_p$ is the $p$-dimensional vector of zeros 
and $\bI_p$ is the $p\times p$ identity matrix. 
\end{proposition}

%% file: proof_prop2.tex
Before we prove Proposition \ref{prop:suff_norm}, we introduce a Lyapunov type bound in Lemma \ref{lem:raic} provided by Theorem 1 of Raic \citep{Raic19}. 

\begin{lemma}
\label{lem:raic}
Consider a sequence of $n \geq 1$ independent random vectors $\bW_i \in \mbR^p$. 
Assume that $\mbE \, \bW_i = \bm{0}_p$ and $\sum_{i=1}^{n} \, \var \, \bW_i = \bI_p$  
where $\bm{0}_p$ is the $p$-dimensional vector of zeros and $\bI_p$ is the $p \times p$ identity matrix.
Define 
\beno
\bS_n 
\= \dsum_{i=1}^{n} \, \bW_i
\ee
and let $\bZ$ be the standard multivariate normal random variable, i.e., $\bZ \sim \text{MvtNorm}(\bm{0}_p, \bI_p)$.  
Then,
for all measurable convex sets $\mA \subset \mbR^p$,  
\beno
\left| \mbP(\bS_n \in \mA) - \Phi(\bZ \in \mA) \right|
&\leq& (42 \, p^{1/4} + 16) \, \dsum_{i=1}^{n} \, \mbE \, \norm{\bW_i}_2^3,
\ee 
where $\Phi$ is the standard multivariate normal measure.
\end{lemma}

\s\s

We now turn to proving Proposition \ref{prop:suff_norm}. 

\s

\pproof \ref{prop:suff_norm}. 
By Proposition \ref{prop:inference} and Lemma \ref{lem:s_hetero},   
the conditional distribution of the multilayer network $\bX$ given $\bY$ 
follows an exponential family with sufficient statistic vector that can be decomposed 
into the sum of conditionally independent dyad-based statistics: 
\beno
\bs(\bX) 
\= \dsum_{\{i,j\} \subset \mN} \, \bs_{i,j}(\bX),
\ee
with the precise formula for $\bs_{i,j}(\bX)$ given in Lemma \ref{lem:s_hetero}. 
Define 
\beno
\bS_{\mN} 
&\coloneqq& 
(I(\truth) \, \norm{\bY}_1)^{-1/2} \, (\bs(\bX) - \mbE^{\bY} \bs(\bX)) \s \\ 
\= \dsum_{\{i,j\} \subset \mN} \, (I(\truth) \, \norm{\bY}_1)^{-1/2} \, (\bs_{i,j}(\bX) - \mbE^{\bY} \bs_{i,j}(\bX)), 
\ee
where $I(\truth)$ is the Fisher information matrix of an activated dyad $X_{i,j}$ for $\{i,j\} \subset \mN$ 
satisfying 
$Y_{i,j} = 1$ evaluated at $\truth$ 
per Lemma \ref{lem:min-eig} 
and where 
$\mbE^{\bY}$ is the random conditional expectation operator with respect to the distribution of $\bX$ conditional on $\bY$. 
For $\gamma > 0$ satisfying $\gamma < \mbE \, \norm{\bY}_1$, define the event $\mE(\gamma)$ by
\beno
\mE(\gamma) & \coloneqq &  \left\{\, \by \in \mbY \,:\, 
\norm{\by}_1 \geq \mbE \norm{\bY}_1 - \gamma \right\}. 
\ee
In words, 
$\mE(\gamma)$ is the subset of configurations of the single-layer network $\bY$ 
which have the number of edges equal to at least the  expected number of activated dyads 
$\mbE \, \norm{\bY}_1$ minus $\gamma > 0$. 
The restrictions placed on $\gamma$ ensure that $\mbE \, \norm{\bY}_1 - \gamma > 0$ 
which implies that $\mE(\gamma)$ will not contain the empty graph which has no edges
and that $\mE(\gamma)$ will contain the complete graph with $\binom{N}{2}$ edges 
as $\mbE \, \norm{\bY}_1 < \binom{N}{2}$ (strict inequality following from the fact that 
$g(\by)$, the marginal probability mass function of $\bY$, 
is assumed to be strictly positive on $\mbY$).
Hence, 
$\mbP(\mE(\gamma)) > 0$ and $\mbP(\mE(\gamma)^c) > 0$.   
Let $\mA \subset \mbR^p$ be a measurable convex set.  
By the law of total probability and the triangle inequality, 
we have
\be
\label{tri_ineq}
\left|\, \mbP(\bS_\mN \in \mA) - \Phi(\bZ \in \mA) \, \right| 
&\leq&
\scalebox{0.85}{$\left|\mbP(\bS_n \in \mA \, | \, \mE(\gamma))  - \Phi(\bZ \in \mA)\right| \, \mbP(\mE(\gamma))$} \s \\
 && +  \; \scalebox{0.85}{$\left|\mbP(\bS_n \in \mA \, | \, \mE^c(\gamma))  - \Phi(\bZ \in \mA)\right| \, \mbP(\mE^c(\gamma))$} \s \\ 
 &\leq & \scalebox{0.85}{$\sup\limits_{\by \in \mE(\gamma)} \,\left| \, \mbP(\bS_{\mN} \in \mA \, | \, \bY = \by) - \Phi(\bZ \in \mA)\,  \right|$} \s \\
 && +  \scalebox{0.85}{$\mbP(\mE^c(\gamma))$},
\ee
noting $\left|\mbP(\bS_n \in \mA \, | \, \mE^c(\gamma))  - \Phi(\bZ \in \mA)\right| \leq 1$ 
and $\mbP(\mE(\gamma)) \leq 1$. 
Taking
\beno
\bW_{i,j} & = & (I(\truth) \, \norm{\bY}_1)^{-1/2} \, (\bs_{i,j}(\bX) - \mbE^{\bY} \, \bs_{i,j}(\bX)),
\ee
we have
\beno
\mbE \, [\bW_{i,j} \, | \, \bY = \by ] \= 0, 
\ee
a result of the tower property of conditional expectation, 
and 
\beno 
\var \, \left[\sum_{\{i,j\} \subset \mN} \, \bW_{i,j}  \, | \, \bY = \by \right] \= \bI_p,  
\ee 
which follows from Lemma \ref{lem:min-eig} which establishes that 
$\var[s_{i,j}(\bX) \,|\, \bY = \by] = I(\truth)$ when $Y_{i,j} = 1$, 
recalling the form of the Fisher information matrix of exponential families to be the 
covariance matrix of the sufficient statistic vector 
\citep[e.g., Proposition 3.10, pp. 32,][]{Su19}, 
and due to the fact that 
$\var[s_{i,j}(\bX) \,|\, \bY = \by] = \bm{0}_{p,p}$ when $Y_{i,j} = 0$. 
Applying Lemma \ref{lem:raic} to the first term of the summation of \eqref{tri_ineq}, for any measurable convex set 
$\mA \subset \mbR^p$,
\beno
\left| \mbP(\bS_{\mN} \in \mA \, | \, \bY = \by) - \Phi(\bZ \in \mA) \right|
\ee
is bounded above by
\beno
 (42 \, p^{1/4} + 16) \, \dsum_{\{i,j\} \subset \mN} \, \mbE \, \left[\norm{\bW_{i,j}}_2^3 \, | \, \bY =\by \right].
\ee
Given $\bY = \by$, using standard matrix and vector norm inequalities,
\beno
\norm{\bW_{i,j}}_2
\= \norm{(I(\truth) \, \norm{\by}_1)^{-1/2} \, (\bs_{i,j}(\bX) - \mbE \, \bs_{i,j}(\bX))}_2 \s \\ 
& \leq & \norm{\by}_1^{-1/2} \,\mnorm{I(\truth)^{-1/2}}_2 \, \norm{\bs_{i,j}(\bX) - \mbE \, \bs_{i,j}(\bX)}_2 \s \\ 
& \leq & (\norm{\by}_1 \, \widetilde{\lambda}_{\min}^{\epsilon})^{-1/2} \, p^{1/2} \, y_{i,j}, 
\ee
where $\mnorm{\cdot}_2$ denotes the spectral norm of a $p \times p$ matrix.
From the proof of Lemma \ref{lem:concentration_likelihood}, for all $l \in \{1, \ldots, p\}$, we have
\beno
0 &\leq & s_{l,i,j}(\bx) &\leq& 1,
&&  \{i,j\} \subset \mN, 
\ee
$\mbP$-almost surely. 
Hence, 
\beno
\mbP(\norm{\bs_{i,j}(\bX) - \mbE^{\bY} \bs_{i,j}(\bX)}_\infty \leq y_{i,j} \,|\, \bY = \by) 
\= 1,
\ee 
implying (conditional on $\bY = \by$) 
\beno
\norm{\bs_{i,j}(\bX) - \mbE^{\bY} \, \bs_{i,j}(\bX)}_2 & \leq & (p \, y_{i,j})^{1/2} 
\= p^{1/2} \, y_{i,j},
\ee
$\mbP$-almost surely. 
As a result,
\beno
\mbE \, \left[\norm{\bW_{i,j}}_2^3 \, | \, \bY =\by \right]
&\leq &  (\norm{\by}_1  \, \widetilde{\lambda}_{\min}^{\epsilon})^{-3/2} \, p^{3/2} \, y_{i,j},
\ee
noting that $y_{i,j}^3 = y_{i,j} \in \{0, 1\}$. 
Using the fact that $42 \, p^{1/4} + 16 \leq 58 \, p^{1/4}$ ($p \geq 1$), we have
\beno
(42 \, p^{1/4} + 16) \, \dsum_{\{i,j\} \subset \mN} \,\mbE \, \left[\norm{\bW_{i,j}}_2^3 \, | \, \bY =\by \right] \s \\
\leq \quad 58 \,  p^{7/4} \, \dsum_{\{i,j\} \subset \mN} \, y_{i,j} \, (\norm{\by}_1 \, \widetilde{\lambda}_{\min}^{\epsilon})^{-3/2} \s \\
=  \quad58 \, p^{7/4} \,\norm{\by}_1^{-1/2} \,(\widetilde{\lambda}_{\min}^{\epsilon})^{-3/2} \s \\
\leq \quad 58 \, p^{7/4} \,(\mbE \, \norm{\bY}_1 - \gamma)^{-1/2} \,(\widetilde{\lambda}_{\min}^{\epsilon})^{-3/2},
\ee
as the conditioning event
$\mE(\gamma)$ and choice of $\gamma$
ensure that $\norm{\by}_1 \geq \mbE \norm{\bY}_1 - \gamma > 0$.
We bound the second term in \eqref{tri_ineq} by Chebyshev's inequality using equation \eqref{ineq:var} in the proof of Lemma \ref{lem:concentration_likelihood}:
\beno
\mbP(\mE^c(\gamma))  
&\leq& \dfrac{\mbE \, \norm{\bY}_1 + 2 \, \left[ D_{g} \right]^{+}}{\gamma^2}.  
\ee 
Taking $\gamma = 2^{-1} \, \mbE \, \norm{\bY}_1 > 0$,
we have  
\beno
\mbP(\mE^c(\gamma))  
&\leq&  \dfrac{4}{\mbE \, \norm{\bY}_1} + \dfrac{8 \, \left[D_{g}\right]^{+}}{\left(\mbE \, \norm{\bY}_1 \right)^2}.   
\ee
Combining terms, 
we obtain the bound
\beno 
\left| \mbP(\bS_\mN \in \mA) - \Phi(\bZ \in \mA) \right|
&\leq& \dfrac{83}{(\widetilde{\lambda}_{\min}^{\epsilon})^{3/2}} \, 
\sqrt{\dfrac{p^{7/2}}{\mbE \, \norm{\bY}_1}}
+  \dfrac{4}{\mbE \, \norm{\bY}_1} + \dfrac{8 \, \left[D_{g}\right]^{+}}{\left(\mbE \, \norm{\bY}_1 \right)^2}.
\ee
\qed
\s \s

Note that the asymptotic multivariate normality can be established provided 
\beno
\lim\limits_{N \to \infty} \, \left[
\dfrac{83}{(\widetilde{\lambda}_{\min}^{\epsilon})^{3/2}} \,
\sqrt{\dfrac{p^{7/2}}{\mbE \, \norm{\bY}_1}}
+  \dfrac{4}{\mbE \, \norm{\bY}_1} + \dfrac{8 \, \left[D_{g}\right]^{+}}{\left(\mbE \, \norm{\bY}_1 \right)^2} 
\right]
\= 0,  
\ee
resulting in the following asymptotic convergence in distribution:  
\beno
\bS_\mN &\overset{D}{\longrightarrow}& \bZ &\sim& \text{MvtNorm}\left(\textbf{0}, \, \bI_p \right). 
\ee

%% file: proof_thm2.tex
In order to demonstrate the feasibility of the normal approximation for maximum likelihood estimators $\mle$  of $\truth$,
we first start with a standard Taylor expansion of the negative score equation:
\be
\label{eq:expansion} 
-\nabla_{\nat} \, \ell(\mle; \bx, \by)
\= -\nabla_{\nat} \,  \ell(\truth; \bx, \by) - \nabla_{\nat}^2 \, \ell(\truth; \bx, \by) \, (\mle - \truth) - \bR,
\ee
where $\bR \in \mbR^p$ is the vector of remainders given in the Lagrange form. 
Denoting by $R_i$, $(\mle - \truth)_i$, and $(\nabla_{\nat} \, \ell(\nat;\bx,\by))_i$ the $i^{\text{th}}$ component of $\bR$, $(\mle - \truth)$,
and the score function 
$\nabla_{\nat} \, \ell(\nat;\bx,\by)$,
respectively.
The remainder term $R_i$ ($i = 1, \ldots, p$) is given by 
\beno
R_i 
\=  \dsum_{j=1}^p \, \dfrac{1}{2} \, 
\dfrac{\partial^2\,(\nabla_{\nat} \,\ell(\dot\nat_i;\bx,\by))_i}{\partial \, \theta_j^2} 
(\mle - \truth)_j^2 \s \\
&& + \; \dsum_{1 \le j < k \le p} \, \dfrac{\partial^2 \, 
(\nabla_{\nat} \ell(\dot\nat_i; \bx,\by))_i}{\partial\, \theta_j \, \partial\, \theta_k} \, 
(\mle - \truth)_j \, (\mle - \truth)_k,
\ee
where $\dot\nat_i = t_i \, \mle + (1 - t_i) \, \truth$ (for some $t_i \in [0, 1]$).  
By Proposition \ref{prop:inference}, 
\beno
\ell(\nat; \bx, \by) 
\= \log \, \mbP_{\nat}(\bX = \bx \,|\, \bY = \by) + \log \, g(\by). 
\ee
By Lemma \ref{lem:s_hetero}, 
the probability mass function $\mbP_{\nat}(\bX = \bx \,|\, \bY = \by)$ 
belongs to a minimal exponential family with the sufficient statistic vector $\bs(\bx)$ 
given by equation \eqref{eq:suff} in Lemma \ref{lem:s_hetero}.  
We then have,
\beno
& - \nabla_{\nat} \, \ell(\nat;\bx, \by) 
\= -(\bs(\bx) - \mbE_{\nat}^{\by} \, \bs(\bX)) \s \\ 
& - \nabla_{\nat}^2 \, \ell(\nat;\bx,\by)
\= \var_{\nat}^{\by} \,  \bs(\bX)
\quad = \quad \mcI(\truth) \, \norm{\by}_1,
\ee
where $\mbE_{\nat}^{\by}$ and $\var_{\nat}^{\by}$ are the 
conditional expectation and variance operators,
respectively, 
of the conditional distribution of $\bX$ given $\bY = \by$,
and by using standard formulas for exponential families 
\citep[e.g., Proposition 3.8, pp. 29,][]{Su19}
and the results of Lemma \ref{lem:min-eig}.  
Note $\nabla_{\nat} \,  \ell(\mle; \bx,\by) = 0$, 
as the maximum likelihood estimator $\mle$ solves the score equation by definition. 
We re-arrange \eqref{eq:expansion} and multiply both sides by $(\mcI(\truth) \, \norm{\bY}_1)^{-1/2}$ to obtain
\be
\label{eq:bridge}
&& (\mcI(\truth) \, \norm{\bY}_1)^{1/2} \, (\mle - \truth) - \; (\mcI(\truth) \, \norm{\bY}_1)^{-1/2} \, \bR\s \\ 
\= (\mcI(\truth) \, \norm{\bY}_1)^{-1/2} \, (s(\bX) - \mbE^{\bY} \, s(\bX)). 
\ee
Define $\Delta \coloneqq  (\mcI(\truth) \, \norm{\bY}_1)^{-1/2} \, \bR$.
Let $\mA \subset \mbR^p$ be any measurable convex subset of $\mbR^p$ 
and $\bZ \sim \text{MvtNorm}(\bm{0}_p, \, \bI_p)$.
We are interested in bounding the quantity  
\beno
\left|\mbP((\mcI(\truth) \, \norm{\bY}_1)^{1/2} \, (\mle - \truth) - \Delta \in \mA) - \Phi(\bZ \in \mA)\right|.
\ee 
Then from \eqref{eq:bridge}, 
\beno
\mbP\left((\mcI(\truth) \, \norm{\bY}_1)^{1/2} \, (\mle - \truth) - \Delta \in \mA\right) \s \\
 = \quad \mbP\big((\mcI(\truth) \, \norm{\bY}_1)^{-1/2} \, (s(\bX) - \mbE^{\bY} \, s(\bX)) \in \mA\big).  
\ee
Applying Proposition \ref{prop:suff_norm},
for all measurable convex sets $\mA \subseteq \mbR^p$, 
\beno
&& \left| \mbP\left((\mcI(\truth) \, \norm{\bY}_1)^{-1/2} \, (s(\bX) - \mbE^{\bY} \, s(\bX)) \in \mA\right) 
- \Phi(\bZ \in \mA) \right|
\ee
is bounded above by
\beno
\dfrac{83}{(\widetilde{\lambda}_{\min}^{\epsilon})^{3/2}} \,
\sqrt{\dfrac{p^{7/2}}{\mbE \, \norm{\bY}_1}}
+  \dfrac{4}{\mbE \, \norm{\bY}_1} + \dfrac{8 \, \left[D_{g}\right]^{+}}{\left(\mbE \, \norm{\bY}_1 \right)^2}. 
\ee
Hence, 
\beno
\left|\mbP((\mcI(\truth) \, \norm{\bY}_1)^{1/2} \, (\mle - \truth) - \Delta \in \mA) - \Phi(\bZ \in \mA)\right| 
\ee
is bounded above by 
\beno
\dfrac{83}{(\widetilde{\lambda}_{\min}^{\epsilon})^{3/2}} \,
\sqrt{\dfrac{p^{7/2}}{\mbE \, \norm{\bY}_1}}
+  \dfrac{4}{\mbE \, \norm{\bY}_1} + \dfrac{8 \, \left[D_{g}\right]^{+}}{\left(\mbE \, \norm{\bY}_1 \right)^2}. \s 
\ee
We lastly handle the term $\Delta$
by showing that $\norm{\Delta}_2$ is small with high probability. 
We first use standard vector/matrix norm inequalities to bound 
\beno
\norm{\Delta}_2
\= \norm{(\mcI(\truth) \, \norm{\bY}_1)^{-1/2} \, \bR}_2 
&\leq& \dfrac{\mnorm{\mcI(\truth)^{-1/2}}_2}{\sqrt{\norm{\bY}_1}} \; \norm{\bR}_2
&\leq& \dfrac{\norm{\bR}_2}{\sqrt{\widetilde{\lambda}_{\min}^{\epsilon}  \norm{\bY}_1}},  
\ee
noting that the spectral norm $\mnorm{\mcI(\truth)^{-1/2}}_2$ 
is equal to the largest eigenvalue of $\mcI(\truth)^{-1/2}$ which will be the 
reciprocal of the smallest eigenvalue of $\mcI(\truth)^{1/2}$,
which is bounded below by $\sqrt{\widetilde{\lambda}_{\min}^{\epsilon}}$. 
Using a standard result from the Taylor theorem for functions with multiple variables,
if for each $i = 1, \ldots, p$,
there exists constants $M_i > 0$ such that
\beno
\sup\limits_{\nat \in \mbR^p \,:\, \norm{\nat - \truth}_1 \leq \norm{\mle - \truth}_1} \,
\left| \, \dfrac{\partial^2 (\nabla_{\nat} \, \ell(\nat;\bx,\by))_i}
{\partial\, \theta_j \, \partial\, \theta_k} \, \right|
&\leq& M_i,
&& 1 \leq j \leq k \leq p, 
\ee
then the Lagrange remainder is bounded above by
\beno
R_i & \le & \dfrac{M_i}{2} \, \norm{\mle - \truth}_1^2
\ee
on the set $\{\nat \in \mbR^p : \norm{\nat - \truth}_1 \leq \norm{\mle - \truth}_2\}$. 
By Lemma \ref{lem:Taylor},
conditional on $\bY = \by$,
we have,
for all $i = 1, \ldots, p$, 
the bound $M_i \le 2 \, \norm{\by}_1$. 
Hence,
\be
\label{eq:Rbound1}
\norm{\Delta}_2
&\leq& \scalebox{0.9}{$\dfrac{1}{\sqrt{\widetilde{\lambda}_{\min}^{\epsilon} \, \norm{\by}_1}} \; 
\sqrt{\dsum_{i=1}^{p} \, R_i^2}  
\quad\leq\quad \dfrac{1}{\sqrt{\widetilde{\lambda}_{\min}^{\epsilon} \, \norm{\by}_1}} \, 
\sqrt{\dsum_{i=1}^{p} \, \norm{\by}_1^2 \, \norm{\mle - \truth}_1^4}$} \s\s\\
&\leq& \scalebox{0.9}{$\dfrac{1}{\sqrt{\widetilde{\lambda}_{\min}^{\epsilon} \, \norm{\by}_1}} \,
\sqrt{p \, \norm{\by}_1^2 \, \norm{\mle - \truth}_1^4}  
\quad\leq\quad \dfrac{\sqrt{p} \, \norm{\by}_1 \, \norm{\mle - \truth}_1^2}{\sqrt{\widetilde{\lambda}_{\min}^{\epsilon} \, \norm{\by}_1}}$} \s\s\\
&\leq& \dfrac{\sqrt{p} \, \sqrt{\norm{\by}_1} \, p \, \norm{\mle - \truth}_2^2}{\sqrt{\widetilde{\lambda}_{\min}^{\epsilon}}} 
\quad\leq\quad \dfrac{p^{3/2}\, \sqrt{\norm{\by}_1} \, \norm{\mle - \truth}_2^2}{\sqrt{\widetilde{\lambda}_{\min}^{\epsilon}}}. 
\ee
By Chebyshev's inequality---as in the proof of Lemma \ref{lem:concentration_likelihood}---we can establish that
\be
\label{eq:event1}
\mbP\left(\left|\norm{\bY}_1 - \mbE \, \norm{\bY}_1\right| > \dfrac{1}{2} \, \mbE \, \norm{\bY}_1\right)
&\leq& \dfrac{4}{\mbE \, \norm{\bY}_1} + \dfrac{8 \, [D_{g}]^{+}}{(\mbE \, \norm{\bY}_1)^2}.
\ee
Under Assumptions \ref{assump1}, \ref{assump2} and \ref{assump3}, Theorem \ref{thm1} established that there exist constants $C>0$ and $N_0 \geq 3$ such that,
for all $N \geq N_0$, 
the event 
\be
\label{eq:event2}
\norm{\mle - \truth}_2 &\leq&
C \; \dfrac{\sqrt{\widetilde{\lambda}_{\max}^{\star}} }{\widetilde{\lambda}_{\min}^{\epsilon}} \; \sqrt{\dfrac{p}{\mbE \norm{\bY}_1}}
\ee
occurs with probability at least $1 - \exp\,(-2\,p) \, - \, (\mbE \, \norm{\bY}_1)^{-1}$.
Define $\mE_1$ to be the event
\beno 
|\norm{\bY}_1 - \mbE \, \norm{\bY}_1| 
&\leq& \dfrac{1}{2} \, \mbE \, \norm{\bY}_1
\ee
and $\mE_2$ to be the event in \eqref{eq:event2},
and define $\mcR$ to be the corresponding values of $\Delta$ in
the event $(\bX, \bY) \in \mE_1 \cap \mE_2$,
under which we have the bound
\be
\label{eq:Rbound2}
\norm{\Delta}_2
&\leq& \dfrac{p^{3/2}\, \sqrt{\norm{\by}_1}}{\sqrt{\widetilde{\lambda}_{\min}^{\epsilon}}} \, 
C^2 \; \dfrac{\widetilde{\lambda}_{\max}^{\star} }{(\widetilde{\lambda}_{\min}^{\epsilon})^2} \; \dfrac{p}{\mbE \norm{\bY}_1} \s\s\\
&\leq& \dfrac{C^2 \, p^{5/2} \, \sqrt{2 \, \mbE \, \norm{\bY}_1}}{\mbE \norm{\bY}_1} \, 
\dfrac{\widetilde{\lambda}_{\max}^{\star} }{(\widetilde{\lambda}_{\min}^{\epsilon})^{5/2}} \s \\
\= \dfrac{\sqrt{2} \, C^2 \, p^{5/2}}{\sqrt{\mbE \norm{\bY}_1}} \, 
\dfrac{\widetilde{\lambda}_{\max}^{\star} }{(\widetilde{\lambda}_{\min}^{\epsilon})^{5/2}}.
\ee
The first inequality in \eqref{eq:Rbound2} is obtained by combining the bounds in \eqref{eq:Rbound1} and \eqref{eq:event2}. The second inequality in \eqref{eq:Rbound2} is using the fact that 
\beno
\norm{\by}_1 
&\leq&  \mbE \, \norm{\bY}_1 + \dfrac{1}{2} \, \mbE \, \norm{\bY}_1 
&\leq& 2 \, \mbE \, \norm{\bY}_1
\ee 
in the event $\by \in \mE_1$.  
Moreover,
a union bound shows that
\beno
\mbP(\Delta \not\in \mcR) &\leq&
\mbP(\mE_1^c) \, + \, \mbP(\mE_2^c)
\s \\
&\leq & \exp\,(-2\,p) \, + \, \dfrac{5}{\mbE \, \norm{\bY}_1} + \dfrac{8 \, [D_{g}]^{+}}{(\mbE \, \norm{\bY}_1)^2} \s \\
&\leq& \exp\,(-2\,p) \, + \, \dfrac{5+ 8 \, C_0}{\mbE \, \norm{\bY}_1}, 
\ee
where the constant $C_0$ and the last inequality follow from Assumption \ref{assump1}.
Hence, 
\beno
\label{eq:R_bound}
\mbP\left(\norm{\Delta}_2 \leq 
\dfrac{\sqrt{2} \, C^2 \, p^{5/2}}{\sqrt{\mbE \norm{\bY}_1}} \, 
\dfrac{\widetilde{\lambda}_{\max}^{\star} }{(\widetilde{\lambda}_{\min}^{\epsilon})^{5/2}} \right)
&\geq& 1 - \exp\,(-2\,p) \, - \, \dfrac{5+ 8 \, C_0}{\mbE \, \norm{\bY}_1}. 
\ee
Taken together, 
we have shown under the assumptions of Theorem \ref{thm1} that there exists $N_0 \geq 3$ such that,
for all $N \geq N_0$,
the error of the multivariate normal approximation 
\beno
\left|\mbP((\mcI(\truth) \, \norm{\bY}_1)^{1/2} \, (\mle - \truth) - \Delta \in \mA) - \Phi(\bZ \in \mA)\right|
\ee
is bounded above by
\beno 
\dfrac{83}{(\widetilde{\lambda}_{\min}^{\epsilon})^{3/2}} \,
\sqrt{\dfrac{p^{7/2}}{\mbE \, \norm{\bY}_1}}
+  \dfrac{4}{\mbE \, \norm{\bY}_1} + \dfrac{8 \, \left[D_{g}\right]^{+}}{\left(\mbE \, \norm{\bY}_1 \right)^2}
\ee
where $\Delta$ satisfies 
\beno
\mbP\left(\norm{\Delta}_2 \leq 
\dfrac{\sqrt{2} \, C^2 \, p^{5/2}}{\sqrt{\mbE \norm{\bY}_1}} \, 
\dfrac{\widetilde{\lambda}_{\max}^{\star} }{(\widetilde{\lambda}_{\min}^{\epsilon})^{5/2}} \right)
&\geq& 1 - \exp\,(-2\,p) \, - \, \dfrac{5+ 8 \, C_0}{\mbE \, \norm{\bY}_1}.
\ee

\qed

\subsection{Auxiliary results for proof of Theorem \ref{thm2}} 
\begin{lemma}
\label{lem:Taylor}
Consider a separable multilayer network model following the form of equation \eqref{general_model} and is 
defined on a set of $N \geq 3$ 
and $K \geq 1$ layers and denote by $\ell(\nat; \bx,\by)$ the log-likelihood function.
Then,
for each  $i = 1, \ldots, p$ and $\epsilon > 0$,
\beno
\sup\limits_{\nat \in \mbR^p \,:\, \norm{\nat - \truth}_2 \leq \epsilon} \quad  
\left| \dfrac{\partial^2 \, (\nabla_{\nat} \, \ell(\nat;\bx,\by))_i}{\partial\, \theta_j \, \partial\, \theta_k} \right|
&\leq& 2 \, \norm{\by}_1,
\ee 
where $(\nabla_{\nat} \, \ell(\nat;\bx,\by))_i$ is the $i^{th}$ component of the score function 
$\nabla_{\nat} \, \ell(\nat;\bx,\by)$. 
\end{lemma}

\s 

\llproof \ref{lem:Taylor}. 
By Proposition \ref{prop:inference}, 
given the observation $\bx$ of $\bX$ (i.e., observation of the event $\bX = \bx$), 
$\bY$ is predictable with unique value $\by \in \mbY$ 
given by the formula in Proposition \ref{prop:inference},
and $(\bx, \by)$ is network concordant. 
Further, 
by Proposition \ref{prop:inference}  
\beno
\ell(\nat; \bx, \by)
\= \log \, \mbP_{\nat}(\bX = \bx \mid \bY = \by) 
+ \log \, g(\by), 
\ee
where $\log \, \mbP_{\nat}(\bX = \bx | \bY = \by)$ 
is the log-likelihood of a minimal, full, and regular 
exponential family. 
Thus, 
the second order derivative of 
$\ell(\nat; \bx, \by)$ with respect to the $i^{\text{th}}$ and $j^{\text{th}}$ 
components of $\nat$ correspond to the variance (in the case $i = j$) 
or covariance (in the case of $i \neq j$) 
of corresponding sufficient statistic(s) of the exponential family \citep[e.g., Proposition 3.8, p. 29,][]{Su19},
with sufficient statistics given in Lemma \ref{lem:s_hetero}.  
For $\{i, j\} \subseteq \{1, \ldots,p\}$, 
\beno
\dfrac{\partial \, (\nabla_{\nat} \, \ell(\nat; \bx, \by))_i}{\partial\, \theta_j} 
\= \dfrac{\partial^2 \, \ell(\nat;\bx,\by)}{\partial\, \theta_i \; \partial\, \theta_j}
\= \cov_{\nat}(s_i(\bX),s_j(\bX) \,|\, \bY = \by),
\ee
and when $i = j \in \{1, \ldots, p\}$,
\beno
\dfrac{\partial \, (\nabla_{\nat} \, \ell(\nat;\bx,\by))_i}{\partial\, \theta_i} 
\= \dfrac{\partial^2 \, \ell(\nat; \bx, \by)}{\partial \, \theta_i^2}
\= \var_{\nat}(s_i(\bX) \,|\, \bY = \by).
\ee
As a result, for $\{i,j\} \subseteq \{1,\ldots,p\}$ and $k \in \{1, \ldots, p\}$,
\beno
\left| \, \dfrac{\partial^2 \, (\nabla_{\nat} \, \ell(\nat;\bx,\by))_i}{\partial\, \theta_j \, \partial\, \theta_k} \, \right| 
\=
\left| \, \dfrac{\partial \,\cov_{\nat}(s_i(\bX),s_j(\bX) \,|\, \bY = \by)}{\partial \, \theta_k} \, \right|,
\ee
and when $i = j \in \{1, \ldots, p\}$ and $k \in \{1, \ldots, p\}$,
\beno
\left| \, \dfrac{\partial^2 \, (\nabla_{\nat} \, \ell(\nat;\bx,\by))_i}{\partial\, \theta_i \, \partial\, \theta_k} \, \right| 
\=
\left| \, \dfrac{\partial \,\var_{\nat}(s_i(\bX) \,|\, \bY = \by)}{\partial \, \theta_k} \, \right|.
\ee
By Lemma \ref{lem:s_hetero} equation \eqref{eq:suff}, 
conditional on $\bY = \by$, 
each sufficient statistic $s_i(\bX)$ ($i \in \{1, \ldots, p\}$) 
can be decomposed into the sum of conditionally independent statistics of each dyad 
$\bX_{v,w}$, for $\{ v,w \} \subseteq \mN$. 
We can then write 
\beno
\scalebox{0.95}{$\cov_{\nat}(s_i(\bX),s_j(\bX) \,|\, \bY = \by)$} 
\=  \scalebox{0.95}{$\dsum_{\{v,w\} \subset \mN} \, \cov_{\nat}(s_{i,v,w}(\bX_{v,w}), \, s_{j,v,w}(\bX_{v,w}) \,|\, \bY = \by)$}, 
\ee
noting that by conditional independence 
$\cov_{\nat}(s_{i,v,w}(\bX_{v,w}),s_{j,r,t}(\bX_{r,t}) \,|\, \bY = \by) = 0$
whenever $\{r,t\} \neq \{v,w\}$, 
and when $i = j$, we can write
\beno
\var_{\nat}(s_i(\bX) \,|\, \bY = \by) 
\= \dsum_{\{v,w\} \subset \mN} \,  \var_{\nat}(s_{i,v,w}(\bX_{v,w}) \,|\, \bY = \by),
\ee
again appealing to the conditional independence given $\bY$ of the random variables 
$s_{i,v,w}(\bX_{v,w})$ ($\{v,w\} \subset \mN$). 
As a result, for $k \in \{1, \ldots, p\}$, it suffices to show that, 
\beno
\left|\,\dfrac{\partial \,\cov_{\nat}(s_{i,v,w}(\bX_{v,w}), \, s_{j,v,w}(\bX_{v,w}) \,|\, \bY = \by )}{\partial \, \theta_k} \,\right|
& \leq & 2,
\ee
and 
\beno
\left| \,\dfrac{\partial \,\var_{\nat}(s_{i,v,w}(\bX_{v,w}) \,|\, \bY = \by)}{\partial \, \theta_k} \, \right|
& \leq & 1.
\ee
Recall that the sufficient statistic $s_{i,v,w}(\bX)$ ($i = 1, \ldots ,p$) 
is defined in Lemma \ref{lem:s_hetero} by
\beno
s_{i,v,w}(\bX_{v,w}) \= \dprod_{t=1}^{h} \, X_{v,w}^{(k_{t})},  
&& \{v,w\} \subset \mN, 
\ee
for some $h \in \{1, \ldots, H\}$ 
and $\{k_1, \ldots, k_h\} \subseteq \{1, \ldots, K\}$. 
Define the set $S_{i,v,w}$ of components of the sufficient statistic vector $\bs_{v,w}(\bX)$ for $\{v,w\} \subset \mN$ and $i = 1,\ldots,p$ by
\beno
S_{i,v,w} \; \coloneqq \; \left\{ \dprod_{t=1}^{h^\prime} \, X_{v,w}^{(l_{t})} \; : \; \{l_1, \ldots, l_{h^\prime} \}  \subset \{k_1,\ldots,k_h\}, \; h^\prime < h \right\},
\ee
where $h \in \{1, \ldots, H\}$ and $\{k_1, \ldots, k_h\} \subseteq \{1, \ldots, K\}$.
The set $S_{i,v,w}$ is the set of components of the sufficient statistic vector $\bs_{v,w}(\bX)$ of dyad $\{v,w\} \subset \mN$ that have a value of 1 when $s_{i,v,w}(\bX) = 1$.
For the ease of notation, 
let $I_{S_{i,v,w}}$ be the set of dimension indices whose corresponding components of the sufficient statistic vector $\bs_{v,w}(\bX)$ belong to the set $S_{i,v,w}$:
\beno
I_{S_{i,v,w}} \; \coloneqq \; \{ j \in \{1,\ldots,p\} \; : \; s_{j,v,w}(\bX) \in S_{i,v,w} \}.
\ee
Define the conditional expectation of $s_{i,v,w}(\bX)$ given $\bY = \by$ for any $i \in \{1,\ldots,p\}$ and $\{v,w\} \subset \mN$ by
\beno
P_{i,v,w}(\nat ; \bX, \by) \; \coloneqq  \; \mbP_{\nat} \, (s_{i,v,w}(\bX) = 1 \, | \, \bY = \by). 
\ee
Denote by $L_i$ the set of layer indices $\{k_1, \ldots, k_h\} \subseteq \{1, \ldots, K\}$ that define the $i^{th}$ component $s_{i,v,w}(\bX_{v,w})$ of the sufficient statistic vector $\bs_{v,w}(\bX)$ 
for any $\{v,w\} \subset \mN$, $j \in \{1,\ldots,p\}$, and some $h \in \{1, \ldots, H\}$.
We then define
\beno
\bX_{v,w }^{(L_i)} & \coloneqq & \left\{X_{v,w}^{(k_1)}, \ldots, X_{v,w}^{(k_h)} \right\}, & \; &
\bX_{v,w }^{(-L_i)} \; \coloneqq \; \bX_{v,w} \setminus \bX_{v,w }^{(L_i)}, 
\ee
and the corresponding sample space
\beno
\mbX_{v,w }^{(L_i)} & \coloneqq & \{0,1\}^h, & \; &
\mbX_{v,w }^{(-L_i)} \; \coloneqq \; \{0,1\}^{H-h}, 
\ee
for some $h \in \{1, \ldots, H\}$.
Then we can write
\beno
P_{i,v,w}(\nat ; \bX, \by) 
\= \mbP_{\nat}\left(\dprod_{l \in L_i} \, X_{v,w}^{(l)} = 1 \, | \, \bY = \by \right) \s \\
\= \dfrac{\dsum_{\mbX_{v,w }^{(-L_i)}} \, \exp\,\left(\sum_{j \in I_{S_{i,v,w}}}\,\theta_j + \sum_{j \in I_{S_{i,v,w}}^c}\, \theta_j\, s_{j,v,w}(\bx)\right)}{\dsum_{\mbX_{v,w}} \, \exp\left( \sum_{j = 1}^p \, \theta_j \, s_{j,v,w}(\bx) \right)}.
\ee
Let 
\beno
Z(\nat) & \coloneqq & \dsum_{\mbX_{v,w}} \, \exp\left( \sum_{j = 1}^p \, \theta_j \, s_{j,v,w}(\bx) \right),
\ee
and take the derivative of $P_{i,v,w}(\nat; \bx,\by)$ with respect to $\theta_k$ for $k = 1,\ldots,p$. We have
\be
\label{cov_ineq}
\dfrac{\partial\, P_{i,v,w}(\nat ; \bX, \by) }{\partial\,\theta_k} \s \\
 \leq \quad
 \scalebox{0.9}{$\dfrac{\dsum_{\mbX_{v,w }^{(-L_i)}} \, \exp\,\left(\sum_{j \in I_{S_{i,v,w}}}\,\theta_j + \sum_{j \in I_{S_{i,v,w}}^c}\, \theta_j\, s_{j,v,w}(\bx)\right) \, \left(Z(\nat) - \dfrac{\partial \, Z(\nat) }{ \partial \, \theta_k}\right)}{Z(\nat)^2}$} \s \\
 =  \scalebox{0.7}{$\dfrac{\dsum_{\mbX_{v,w }^{(-L_i)}} \, \exp\,\left(\sum_{j \in I_{S_{i,v,w}}}\,\theta_j + \sum_{j \in I_{S_{i,v,w}}^c}\, \theta_j\, s_{j,v,w}(\bx)\right) \, \left(\dsum_{\mbX_{v,w}} \, \exp\left( \sum_{j = 1}^p \, \theta_j \, s_{j,v,w}(\bx)\right) \, (1-s_{k,v,w}(\bx))\right)}{Z(\nat)^2}$} \s \\
 \leq \quad 
 \scalebox{0.9}{$\dfrac{\dsum_{\mbX_{v,w }^{(-L_i)}} \, \exp\,\left(\sum_{j \in I_{S_{i,v,w}}}\,\theta_j + \sum_{j \in I_{S_{i,v,w}}^c}\, \theta_j\, s_{j,v,w}(\bx)\right) \, Z(\nat)}{Z(\nat)^2}$} \s \\
 \leq \quad
1.
\ee
The first inequality is obtained because $s_{k,v,w} (\bx) \leq 1$, and the last inequality is due to the fact that
\beno
\dsum_{\mbX_{v,w }^{(-L_i)}} \, \exp\,\left(\sum_{j \in I_{S_{i,v,w}}}\,\theta_j + \sum_{j \in I_{S_{i,v,w}}^c}\, \theta_j\, s_{j,v,w}(\bx)\right) 
\ee
is bounded above by 
\beno
\dsum_{\mbX_{v,w}} \, \exp\left( \sum_{j = 1}^p \, \theta_j \, s_{j,v,w}(\bx) \right).
\ee
Now we turn to show the derivative of the conditional variance and covariance of the sufficient statistics of each dyad are bounded. Given $\bY = \by$, for all $\{i\} \subset \{ 1,\ldots, p\}$, $s_{i,v,w} (\bX)$ are conditionally independent across $\{v,w\} \subseteq \mN$.
Then we have 
\beno
& \cov_{\nat}(s_{i,v,w}(\bX_{v,w}), \, s_{j,v,w}(\bX_{v,w}) \,|\, \bY = \by) \s \\
 = & \scalebox{0.95}{$\mbE \,[ s_{i,v,w}(\bX) \, s_{j,v,w}(\bX) \,|\, \bY = \by] - \mbE\, [s_{i,v,w}(\bX)\,|\, \bY = \by] \, \mbE \, [s_{j,v,w}(\bX)\,|\, \bY = \by]$} \s \\
 = & \mbP_{\nat} \, (s_{i,v,w}(\bX)  = 1, s_{j,v,w}(\bX) =1 \,|\, \bY = \by) - P_{i,v,w} (\nat;\bX,\by) \, P_{j,v,w} (\nat;\bX,\by) \s \\
 = & \mbP_{\nat} \left( \dprod_{l \in L_i \cup L_j} \bX_{v,w}^{(l)} = 1  \, | \, \bY = \by \right) - P_{i,v,w} (\nat;\bX,\by) \, P_{j,v,w} (\nat;\bX,\by).
\ee
Using the inequality derived in \eqref{cov_ineq} and suppressing the notation of $\{v,w\}$ and $(\bX, \by)$ in $P_{i,v,w} (\nat;\bX,\by)$, the derivative of the covariance with respect to $\theta_k$, $k = 1,\ldots, p$ is given by
\beno
& \left| \, \dfrac{\partial \,\cov_{\nat}(s_{i,v,w}(\bX_{v,w}), \, s_{j,v,w}(\bX_{v,w}) \,|\, \bY = \by)}{\partial \, \theta_k} \, \right| \s \\
= & \dfrac{\partial \, \mbP_{\nat} \left(\dprod_{l \in L_i \cup L_j} \bX_{v,w}^{(l)} = 1  \, | \, \bY = \by \right)}{\partial \, \theta_k} - \dfrac{\partial \, P_i(\nat)}{\partial \, \theta_k} \, P_j(\nat) - P_i(\nat) \, \dfrac{\partial \, P_j(\nat)}{\partial \, \theta_k}  \s \\
 \leq & 2.
\ee
Using the same inequality and notation in \eqref{cov_ineq}, the derivative of the variance of a Bernoulli random variable $s_{i,v,w} (\bX)$ is given by
\beno
\left| \,\dfrac{\partial \,\var_{\nat}(s_{i,v,w}(\bX_{v,w}) \,|\, \bY = \by)}{\partial \, \theta_k} \, \right| 
\= \left|\,\left(1 - 2 \, P_i(\nat) \right) \, \dfrac{\partial \, P_i(\nat)}{\partial \, \theta_k} \,\right|
& \leq & 
1.
\ee
Finally, for $\{i,j\} \subseteq \{1,\ldots,p\}$ and $k \in \{1, \ldots, p\}$, we obtain
\beno
\left| \, \dfrac{\partial^2 \, (\nabla_{\nat} \, \ell(\nat;\bx,\by))_i}{\partial\, \theta_j \, \partial\, \theta_k} \, \right| 
\=
\left| \, \dfrac{\partial \,\cov_{\nat}(s_i(\bX),s_j(\bX) \,|\, \bY = \by)}{\partial \, \theta_k} \, \right| \s \\
& \leq & \scalebox{0.95}{$\dsum_{\{v,w\} \subset \mN} \, \left| \dfrac{\partial \, \cov_{\nat}(s_{i,v,w}(\bX_{v,w}), \, s_{j,v,w}(\bX_{v,w}) \,|\, \bY = \by)}{\partial \, \theta_k} \, \right|$} \s \\
& \leq & 2 \, \norm{\by}_1
\ee
due to the fact that $\cov_{\nat}(s_{i,v,w}(\bX_{v,w}), \, s_{j,v,w}(\bX_{v,w}) \,|\, \bY = \by) = 0$ when $Y_{v,w} = 0$ for $\{v,w\} \subset \mN$.
Similarly, $\var_{\nat}(s_{i,v,w}(\bX_{i,v,w}) \,|\, \bY = \by) = 0$ when $Y_{v,w} = 0$ for $\{v,w\} \subset \mN$, 
and when $i = j \in \{1, \ldots, p\}$ and $k \in \{1, \ldots, p\}$, we have
\beno
\left| \, \dfrac{\partial^2 \, (\nabla_{\nat} \, \ell(\nat;\bx,\by))_i}{\partial\, \theta_i \, \partial\, \theta_k} \, \right| 
\=
\left| \, \dfrac{\partial \,\var_{\nat}(s_i(\bX) \,|\, \bY = \by)}{\partial \, \theta_k} \, \right| \s \\
& \leq & \norm{\by}_1.
\ee

\qed